\documentclass[reqno]{amsart}
\usepackage{graphicx, amsmath,amsfonts,amsthm,amssymb, paralist, fullpage, esint}
\newtheorem{theorem}{Theorem}[section]
\newtheorem{corollary}[theorem]{Corollary}
\newtheorem{lemma}[theorem]{Lemma}
\newtheorem{proposition}[theorem]{Proposition}

\theoremstyle{definition}
\newtheorem{definition}[theorem]{Definition}

\theoremstyle{remark}
\newtheorem{assumption}[theorem]{Assumption}
\newtheorem{remark}[theorem]{Remark}
\numberwithin{equation}{section}

\newcommand{\g}{\geqslant}

\newcommand{\RR}{\mathbb{R}}
\newcommand{\ZZ}{\mathbb{Z}}
\newcommand{\CC}{\mathbb{C}}
\newcommand{\s}{\mathcal{S}}

\newcommand{\sph}{\mathbb{S}}
\newcommand{\NN}{\mathbb{N}}
\newcommand{\p}{\partial}

\newcommand{\les}{\leqslant}
\newcommand{\lesa}{\lesssim}

\newcommand{\mc}[1]{\mathcal{#1}}
\newcommand{\mb}[1]{\mathbf{#1}}
\newcommand{\eref}[1]{(\ref{#1})}
\newcommand{\lr}[1]{ \left\langle #1 \right\rangle}
\newcommand{\ind}{\mathbbold{1}}

\DeclareSymbolFont{bbold}{U}{bbold}{m}{n}
\DeclareSymbolFontAlphabet{\mathbbold}{bbold}

\DeclareMathOperator*{\supp}{supp}
\DeclareMathOperator*{\diam}{diam}

\DeclareMathOperator*{\dist}{dist}

\DeclareMathOperator*{\conv}{Convex Hull}

\setdefaultenum{(i)}{(a)}{1}{A}

\begin{document}

\title{Multi-scale Bilinear Restriction Estimates for General Phases}%
\author[T.~Candy]{Timothy Candy}%
\address[T.~Candy]{Universit\"at Bielefeld, Fakult\"at f\"ur Mathematik,
  Postfach 100131, 33501 Bielefeld, Germany}%
\email{tcandy@math.uni-bielefeld.de}%
\thanks{Financial support by the DFG through the CRC ``Taming uncertainty and
profiting from randomness and low regularity in analysis, stochastics
and their applications'' is acknowledged.}

\begin{abstract}
  We prove (adjoint) bilinear restriction estimates for general phases at different scales in the full non-endpoint mixed norm range, and give bounds with a sharp and explicit dependence on the phases. These estimates have applications to high-low frequency interactions for solutions to partial differential equations, as well as to the linear restriction problem for surfaces with degenerate curvature. As a consequence, we obtain new bilinear restriction estimates for elliptic phases and wave/Klein-Gordon interactions in the full bilinear range, and give a refined Strichartz inequality for the Klein-Gordon equation. In addition, we extend these bilinear estimates to hold in adapted function spaces by using a  transference type principle which holds for vector valued waves.
\end{abstract}

\maketitle

\section{Introduction}

Let $n \g 2$ and for $j=1,2$ take phases $\Phi_j:\Lambda_j \rightarrow \RR$ with $\Lambda_j \subset \RR^n$. We consider the problem of obtaining (adjoint) bilinear restriction estimates of the form
            \begin{equation}\label{eqn:intro bilinear bound}
                \| e^{it \Phi_1(-i\nabla)}f e^{ i t \Phi_2(-i\nabla)} g \|_{L^q_t L^r_x(\RR^{1+n})} \les \mb{C} \| f\|_{L^2_x(\RR^n)} \| g \|_{L^2_x(\RR^n)}
            \end{equation}
with an explicit dependence of the constant $\mb{C}$ on the phases $\Phi_j$ and sets $\Lambda_j$. Here $1 \les q, r \les 2$, $\supp \widehat{f} \subset \Lambda_1$, $\supp \widehat{g} \subset \Lambda_2$, and $\widehat{f}(\xi) = \int_{\RR^n} f(x) e^{ -i x\cdot \xi} dx$ is the spatial Fourier transform of $f$. This problem is interesting for two key reasons. Firstly, in applications to PDE, the dependence on $\Phi_j$ and $\Lambda_j$ reflects the \emph{derivative} cost required to place the product into $L^q_t L^r_x$. Controlling the number of derivatives required in \eref{eqn:intro bilinear bound} leads to, for instance, sharp null form estimates in the case of the wave equation $\Phi = |\xi|$ \cite{Tao2001b, Lee2008, Lee2008a}, refined Strichartz estimates for the wave and Klein-Gordon equations \cite{Bourgain1998, Killip2012b, Ramos2012}, as well as scattering for the Wave Maps equation \cite{Sterbenz2010}, and Dirac-Klein-Gordon equation \cite{Candy2017, Candy2017a}. Secondly, bilinear restriction estimates of the form \eref{eqn:intro bilinear bound} have played an important role in the  classical linear restriction problem in harmonic analysis, see for instance [15, 20]. In particular, understanding the dependence of the constant $\mb{C}$ on the phases $\Phi_j$, has led to improved linear restriction estimates for surfaces with curvature that may degenerate in certain directions \cite{Buschenhenke2017, Stovall2017}.

We now consider some concrete examples. In the case of the wave equation $\Phi_1 = \Phi_2 = |\xi|$, the estimate \eref{eqn:intro bilinear bound} has a long history and is now essentially well understood. More precisely, if $\Lambda_1 \subset \{ |\xi| \approx 1 \}$ is at scale 1,  $\Lambda_2 \subset \{ |\xi| \approx \lambda \}$ is at scale $\lambda \g 1$, and $\Lambda_1$ and $\Lambda_2$ have angular separation of size $1$, then it is known that \eref{eqn:intro bilinear bound} holds in the range $\frac{1}{q} + \frac{n+1}{2r}  \les \frac{n+1}{2}$ and we have $ \mb{C} \approx \lambda^{\frac{1}{q} - \frac{1}{2} +\epsilon}$ for every $\epsilon>0$. This is sharp up to the $\epsilon$ frequency loss, and in the unit scale ($\lambda =1$), non-endpoint case $ q = r >\frac{n+3}{n+1}$, is due to the breakthrough work of Wolff \cite{Wolff2001}. The endpoint case $q=r = \frac{n+3}{n+1}$ for general scales $\lambda \g 1$ was  obtained by Tao \cite{Tao2001b}. In the mixed exponents case, $q \not = r$, the bilinear estimate for general scales is due to Tataru \cite{Tataru2003} for the non-constant coefficient wave equation, and Lee-Vargas \cite{Lee2008}, Lee-Rogers-Vargas \cite{Lee2008a} in the constant coefficient case. The endpoint estimate for $q\not = r$ is also known, and was proven by  Temur \cite{Temur2013}.

On the other hand, in the case of the Schr\"{o}dinger equation, $\Phi_1 = \Phi_2 = |\xi|^2$, if  $\Lambda_1$ is a ball of radius $1$, and $\Lambda_2$ is a ball of radius $\lambda \g1$ such that the sets $\Lambda_1$ and $\Lambda_2$ are  separated by a distance $\lambda$, then work of Tao \cite{Tao2003a} shows that for $q= r > \frac{n+3}{n+1}$ the bilinear estimate \eref{eqn:intro bilinear bound} holds with $\mb{C} \approx \lambda^{\frac{1}{r}-1 + \epsilon}$ for every $\epsilon>0$.

In the case of general phases which are both at unit scale, under suitable transversality and curvature assumptions, Lee-Vargas \cite{Lee2010} and Bejenaru \cite{Bejenaru2016b} have shown that \eref{eqn:intro bilinear bound} again holds in the non-endpoint range $q=r>\frac{n+3}{n+1}$. For general phases which are \emph{not} at unit scale, only in certain special cases is the dependence of $\mb{C}$ in \eref{eqn:intro bilinear bound} on the phases $\Phi_j$ and sets $\Lambda_j$ known. In particular, if $\Phi_1$ and $\Phi_2$ are elliptic phases with curvature at different scales, then \eref{eqn:intro bilinear bound} was obtained by Stovall \cite{Stovall2017} with an almost sharp dependence on the scale. If $n=2$,  and $\Phi_1 = \Phi_2 = \xi_1^{m_1} + \xi_2^{m_2}$ are surfaces of finite type with rectangular sets $\Lambda_j \subset \{0<\xi_1, \xi_2<1\}$, then the essentially optimal dependence of the constant $ \mb{C}$ on the rectangles $\Lambda_j$ and parameters $m_j \g 2$ has been obtained in recent work of Buschenhenke-M\"uller-Vargas \cite{Buschenhenke2017}. \\

In the current article we unify and improve the above examples, and show that in fact for general phases $\Phi_j$,  under certain transversality and curvature assumptions on the surfaces $S_j = \{ (\Phi_j(\xi), \xi) | \xi \in \Lambda_j\} \subset \RR^{1+n}$, if (for instance) $\Lambda_1$ is contained in a ball of radius $1$, then the bilinear restriction estimate \eref{eqn:intro bilinear bound} holds in the full, non-endpoint, bilinear range $\frac{1}{q} + \frac{n+1}{2r} < \frac{n+1}{2}$ with constant
         \begin{equation}\label{eqn:informal constant} \mb{C} \approx  \mc{V}_{max}^{\frac{1}{r}-1} \,\big(\min\{\mc{H}_{1}, \mc{H}_2\}\big)^{\frac{1}{2}-\frac{1}{q}} \big(\max\{\mc{H}_1, \mc{H}_2\}\big)^{\frac{1}{2} - \frac{1}{r}}\end{equation}
where we introduce the notation
	$$\mc{V}_{max} = \sup_{\substack{\xi \in \Lambda_1 \\ \eta \in \Lambda_2}} |\nabla \Phi_1(\xi) - \nabla \Phi_2(\eta)|, \qquad \mc{H}_j = \| \nabla^2 \Phi_j \|_{L^\infty(\Lambda_j)}, $$
see Theorem \ref{thm-main} below for a more precise statement. Moreover, we show that this dependence on the phases $\Phi_j$ is sharp. The estimate \eref{eqn:informal constant} recovers all examples mentioned above \emph{without} an epsilon loss. For instance, in the case of the wave equation, a computation shows that provided the sets $\Lambda_1$ and  $\Lambda_2$ are at scales $1$ and  $\lambda\g 1$ respectively, and are angularly separated by a distance $1$, then we have $\mc{V}_{max} \approx 1$ and $\mc{H}_1 = 1$, $\mc{H}_2 \approx \lambda^{-1}$. In particular, for the bilinear restriction estimate for the cone, \eref{eqn:informal constant} recovers the optimal dependence on the frequency  $\lambda$ without an $\epsilon$ loss. Similarly \eref{eqn:informal constant} also recovers the correct dependence on the frequency parameter $\lambda$ in the Schr\"odinger case, after observing that if $\Lambda_1$ and $\Lambda_2$ are separated by a distance $\lambda$, then $\mc{V}_{\max} \approx \lambda$ and $\mc{H}_1 = \mc{H}_2 = 1$.

As a concrete application of the general formula \eref{eqn:informal constant}, we obtain new bilinear restriction estimates for elliptic phases, and Klein-Gordon interactions with differing masses, see Theorem \ref{thm-elliptic bilinear} and Theorem \ref{thm-bilinear small scale KG} respectively. In particular, we obtain sharp (non-endpoint) estimates for wave Klein-Gordon interactions. Further, we give two additional consequences of the general bilinear restriction estimate.  The first, stated in Theorem \ref{thm-refined strichartz} below, is a refined Strichartz estimate for the Klein-Gordon equation which shows that, for $n \g 2$ there exists $\theta>0$ such that
    \begin{equation}\label{eqn:intro refined strichartz} \big\| e^{ it \lr{\nabla}} f \big\|_{L^{2\frac{n+1}{n-1}}_{t,x}(\RR^{1+n})} \lesa  \| f \|_{X}^\theta \| f \|_{H^{\frac{1}{2}}}^{1-\theta}, \end{equation}
 where $\| \cdot \|_X$ is a norm which detects concentration of the Fourier support of $f$ to certain rectangular sets, and satisfies $\| f\|_X \les \| f \|_{H^\frac{1}{2}}$. The classical Strichartz estimate for the Klein-Gordon equation corresponds to $\theta=0$. The estimate \eref{eqn:intro refined strichartz} gives the Klein-Gordon version of a refined Strichartz estimate for the wave equation due to Ramos \cite{Ramos2012} and extends a similar estimate of Killip-Stovall-Visan \cite{Killip2012b}. Bounds of the form \eref{eqn:intro refined strichartz} are a closely related to profile decompositions, see for instance \cite{Killip2012b, Ramos2012} for further details.

  The second consequence is an extension of \eref{eqn:intro bilinear bound} from free waves $e^{i t \Phi(-i \nabla)}f$, to more general functions belonging to certain adapted function spaces $U^2_\Phi$, see Corollary \ref{cor-main Up version} below.  This corresponds to a multi-scale version of recent work of the author and Herr \cite{Candy2016}. The new bilinear bounds we obtain here, are applied in \cite{Candy2017, Candy2017a} to obtain conditional scattering results for the Dirac-Klein-Gordon system. The application to the Dirac-Klein-Gordon equation, which includes Klein-Gordon interactions with different masses, formed a key motivation to study general bilinear restriction estimates of the form \eref{eqn:intro bilinear bound}.

\subsection{Main Theorem} Let $\ell^2_c(\ZZ)$ denote the collection of all compactly supported sequences in $\ell^2(\ZZ)$, thus every sequence in $\ell^2_c(\ZZ)$ has only finitely non-zero components. By a slight abuse of notation, the standard norm on $\ell^2$  is denoted by $|\cdot|$. We define the product of two $\ell^2_c(\ZZ)$-valued functions $u$ and $v$ as simply the tensor product $ uv = u \otimes v$. Thus $|uv|=|u| |v|$.
 We say that $u(t,x):\RR^{1+n}\rightarrow \ell^2_c(\ZZ)$ is a \emph{$\Phi_j$-wave} if
		$$ u(t,x) = e^{ it \Phi_j(-i\nabla)} f(x)$$
 with $f: \RR^n \rightarrow \ell^2_c(\ZZ)$ and $\supp \widehat{f} \subset \Lambda_j$. Note that, if $u$ is a $\Phi_j$-wave, then $\supp \widehat{u}$ is independent of time, $\| u \|_{L^\infty_t L^2_x} = \| u(0) \|_{L^2_x}$, and $u$ is a solution to the equation
        $$ i \p_t u + \Phi_j(-i \nabla) u = 0.$$
 Alternatively, up to a Jacobian factor, $u$ can be viewed as the extension operator associated to the surface $S_j=\{(\Phi_j(\xi), \xi) \mid \xi \in \Lambda_j \} \subset \RR^{1+n}$.\\

It is well-known that to obtain bilinear restriction estimates in the full bilinear range $\frac{1}{q} + \frac{n+1}{2r} \les \frac{n+1}{2}$ requires firstly that the surfaces $S_j$ are \emph{transverse}, in the sense that $|\nabla \Phi_1(\xi) - \nabla \Phi_2(\eta)| >0$,  and secondly that there is some curvature along the $n-1$ dimensional surfaces of intersection $S_1 \cap ( \mathfrak{h} - S_2)$ for $\mathfrak{h} \in \RR^{1+n}$. These conditions are sufficient in the case of the wave or Schr\"{o}dinger equation, but in the general case, a stronger condition is required that ensures that the normals to the surface $S_k$ intersect transversally with the conic surface formed by taking the normal directions to the surface $S_j$ at points in $S_j \cap (\mathfrak{h} - S_k)$, see for instance \cite{Lee2010, Bejenaru2016b, Candy2016} for assumptions of this type. In the current article, we use a normalised version of the conditions appearing in \cite{Candy2016}. To describe the precise conditions we impose on the phases $\Phi_j$, we require some additional notation. Given $\mathfrak{h}=(a, h) \in \RR^{1+n}$ and $\{j, k\}=\{1,2\}$ we define the surfaces
        $$ \Sigma_j(\mathfrak{h}) = \{ \xi \in \Lambda_j\cap (h-\Lambda_k) \mid \Phi_j(\xi) + \Phi_k(h-\xi) = a \}.$$
Note that $\Sigma_j(\mathfrak{h}) \subset \Lambda_j \subset \RR^n$ and these surfaces are related to the intersection of the surfaces $S_j \cap (\mathfrak{h} - S_k)$ through the formula
    $$ S_j \cap (\mathfrak{h} - S_k) = \{ (\Phi_j(\xi), \xi) | \xi \in \Sigma_j(\mathfrak{h}) \}.$$
The sets $\Sigma_j(\mathfrak{h})$ play a crucial role in what follows. The key curvature, transversality, and smoothness assumptions we make on the phases $\Phi_j$ are the following.\\

\begin{assumption}\label{assump-main}
Fix constants $\mb{C}_0, \mb{d}_0 >0$. For $j=1, 2$ we let $\Lambda_j \subset \RR^n$ be an open set, and $\Phi_j:\Lambda_j \rightarrow \RR$ satisfy the conditions:

  \begin{enumerate}
      \item[(\textbf{A1})] for every $\{j, k\} = \{1, 2\}$, $\mathfrak{h} \in
            \RR^{1+n}$, $\xi, \xi' \in \Sigma_j(\mathfrak{h})$, and
            $\eta \in \Lambda_k$ we have
            $$ \big|\big(\nabla \Phi_j(\xi) - \nabla \Phi_j(\xi')\big)
            \wedge \big( \nabla \Phi_j(\xi) - \nabla
            \Phi_k(\eta)\big)\big| \g \mb{C}_0 \mc{V}_{max} \mc{H}_j |\xi - \xi'|$$
            and
                $$ \mb{C}_0 \big| \nabla \Phi_j(\xi) - \nabla \Phi_j(\xi')\big| \les  \mc{H}_j |\xi - \xi'|, $$

       \item[(\textbf{A2})] for $j \in \{1,2\}$ we have  $\Phi_j \in C^{5n}(\Lambda_j)$ and moreover
                         $$  \min_{3\les m \les 5n} \bigg( \frac{\mc{H}_j}{\|\nabla^m \Phi_j\|_{L^\infty(\Lambda_j)}}\bigg)^\frac{1}{m-2} \g \mb{d}_0, \qquad \frac{\mc{V}_{max}}{\mc{H}_j} \g \mb{d}_0.$$\\
   \end{enumerate}
\end{assumption}

The key assumption (\textbf{A1}) is essentially equivalent to the conditions appearing previously in the literature \cite{Lee2010, Bejenaru2016b, Candy2016}, except that we require a normalised version to correctly determine the dependence of the constant on the phases $\Phi_j$. To gain a better geometric understanding of (\textbf{A1}), we observe that letting $\xi' \rightarrow \xi$ and taking $\eta = h-\xi$, since $\nabla\Phi_j(\xi) - \nabla \Phi_k(h -\xi)$ is normal to $\Sigma_j(\mathfrak{h})$,  (\textbf{A1}) implies the \emph{local} condition
	\begin{equation}\label{eqn:intro local cond} \big| \big[ \nabla^2 \Phi_j(\xi) v \big] \wedge n \big| \g \mb{C}_0 \mc{H}_j |v| \end{equation}
for all $v \in T_{\xi} \Sigma_j(\mathfrak{h})$, where $n$ is the unit normal to $\Sigma_j(\mathfrak{h})$, and $T_\xi \Sigma_j(\mathfrak{h})$ is the tangent plane at $\xi \in \Sigma_j(\mathfrak{h})$. Under certain conditions, the local condition \eref{eqn:intro local cond} is in fact equivalent to the global assumption (\textbf{A1}), see Lemma \ref{lem:local version of A1} below. In view of \eref{eqn:intro local cond}, (\textbf{A1}) states that the Hessian $\nabla^2\Phi_j(\xi)$  can not rotate tangent vectors in $\Sigma_j(\mathfrak{h})$ to normal vectors.  A similar observation was made to \cite{Bejenaru2016b} where the assumption (\textbf{A1}) was replaced with a condition on the shape operator associated to the surface $S_j$. We should emphasise here that in the general phase case, it is \emph{not} sufficient to simply assume the surfaces $S_j$ intersect transversally, and a stronger condition of the form (\textbf{A1}) is necessary \cite{Lee2006, Vargas2005}.

There are two immediate consequences of (\textbf{A1}) which we wish to highlight. The first, is that for every $\xi \in \Lambda_1$ and $\eta \in \Lambda_2$ we have the transversality bound
\begin{equation}\label{eqn-main trans assump}
	\big|\nabla \Phi_1(\xi) - \nabla \Phi_2(\eta) \big| \g \mb{C}_0 \mc{V}_{max}
\end{equation}
which follows by observing that for every $\xi \in \Lambda_j$ there exists $\mathfrak{h} \in \RR^{1+n}$ such that  $\xi \in \Sigma_j(\mathfrak{h})$. The second, is that for every $j \in  \{1,2\}$,  $\mathfrak{h} \in \RR^{1+n}$, and  $\xi, \xi'\in \Sigma_j(\mathfrak{h})$ we have the (global) curvature type bound
	\begin{equation}\label{eqn-main curv assump}
	 	\big|\nabla \Phi_j(\xi) - \nabla \Phi_j(\xi') \big| \g \mb{C}_0 \mc{H}_j |\xi - \xi'|.
	\end{equation}
The bounds \eref{eqn-main trans assump} and \eref{eqn-main curv assump} play a key role in the proof of Theorem \ref{thm-main}, and in fact the stronger bound contained in (\textbf{A1}) is only directly required in the proof of Lemma \ref{lem - surface C transverse} below.

The smoothness assumption (\textbf{A2}) requires that the phases to lie in $C^{5n}(\Lambda_j)$, this can probably be improved, but we do not consider this issue here. While the constant $\mb{C}_0$ should be thought of as a universal constant, it is important to keep track of the parameter $\mb{d}_0$ as it encodes the size of the Fourier support of the underlying waves, and scales like the frequency $\xi$. For instance, if $\Phi_1(\xi) = |\xi|^s$ for some $s>0$, then (\textbf{A2}) essentially requires $\mb{d}_0 \lesa |\xi|$ on $\Lambda_1$.  More generally, we observe that both the assumptions (\textbf{A1}) and (\textbf{A2}) are invariant under the rescaling
	\begin{equation} \label{eqn-phase rescaling} \big( \Lambda_j, \Phi_j(\xi), \mb{d}_0 \big) \mapsto \big( \mu^{-1} \Lambda_j, \lambda \Phi_j(\mu\xi), \mu^{-1} \mb{d}_0 \big)\end{equation}
and linear translations
	\begin{equation}\label{eqn-translation invariance} \Phi_j \mapsto \Phi_j + \xi \cdot \xi_0. \end{equation}
These invariances significantly simplify the arguments to follow, as we can then reduce to the case where the maximum velocity and curvature are normalised to be one. In fact the precise assumptions (\textbf{A1}) and (\textbf{A2}) were carefully chosen with the scaling \eref{eqn-phase rescaling} and translation invariance \eref{eqn-translation invariance} in mind.\\

  We are now ready to state our first main theorem.

\begin{theorem}\label{thm-main}
Let $n\g2$, $\mb{C}_0>0$,  $1\les q, r \les 2$, and $\frac{1}{q} +\frac{n+1}{2r} < \frac{n+1}{2}$. There exists a constant $C>0$, such that for any $ \mb{d}_0>0$, any open sets $\Lambda_1, \Lambda_2 \subset \RR^n$, any phases $\Phi_1$ and $\Phi_2$ satisfying Assumption \ref{assump-main} with $\mc{H}_2 \les \mc{H}_1$, and any $\Phi_1$-waves $u$, and $\Phi_2$-waves $v$ satisfying the support conditions
    $$\supp \widehat{u} + \mb{d}_0 \subset \Lambda_1, \qquad \supp \widehat{v} + \mb{d}_0 \subset \Lambda_2, \qquad \min\{\diam(\supp \widehat{u}), \diam(\supp \widehat{v})\} \les \mb{d}_0, $$
we have
		\begin{equation}\label{eqn-intro bilinear data in ball} \big\| uv  \big\|_{L^q_t L^r_x(\RR^{1+n})}\les C\mb{d}_0^{n+1 - \frac{n+1}{r} - \frac{2}{q}} \,\mc{V}_{max}^{\frac{1}{r}-1} \,\mc{H}_{1}^{1-\frac{1}{q} - \frac{1}{r}}\, \Big( \frac{\mc{H}_{1}}{\mc{H}_{2}}\Big)^{\frac{1}{q} - \frac{1}{2}}\| u \|_{L^\infty_t L^2_x} \| v \|_{L^\infty_t L^2_x}. \end{equation}
\end{theorem}

The proof of Theorem \ref{thm-main} follows the strategy used by Tao in the proof of the endpoint bilinear restriction estimate for the cone \cite{Tao2001b}. In particular, we replace the combinatorial Kakeya type arguments used by Wolff in the seminal paper \cite{Wolff2001} and further exploited in, for instance, \cite{Tao2003a, Lee2010, Candy2016, Stovall2017}, with energy estimates across transverse hypersurfaces.  A similar argument was used by Bejenaru \cite{Bejenaru2016b} to obtain a bilinear restriction estimate with $q=r$ for general phases at unit scale. For our purposes, arguing via energy estimates has a number of advantages, firstly it avoids the $\epsilon$ loss that occurs when using the pigeon hole type arguments that arise in the combinatorial approach, and secondly it can be adapted without too much effort to handle interactions with waves at very different scales by using the \emph{wave table} construction of Tao.

\begin{remark}
The extension of the bilinear restriction estimates to \emph{vector valued} waves, was first obtained by Tao \cite{Tao2001b} where the vector valued nature of the waves played an important technical role in the induction argument. As the proof of Theorem \ref{thm-main} follows the argument used by Tao, the fact that we may take vector valued waves in Theorem \ref{thm-main} also plays an important technical role here. However, there is an additional gain that comes from working with vector valued waves, which turns out to be very useful for applications of Theorem \ref{thm-main} to problems in nonlinear PDE. This gain comes from the simple observation that if we have an estimate for vector valued waves, then we can immediately deduce that the same bound holds in the atomic Banach space $U^2_{\Phi_j}$. In particular, Theorem \ref{thm-main} also holds in $U^2_{\Phi_j}$, see Corollary \ref{cor-main Up version} below. The fact that Theorem \ref{thm-main} holds for vector valued waves, is a reflection of the fact that the proof essentially only exploits bilinear estimates in $L^2$, for which the theory for scalar valued and vector valued waves is identical. A previous joint work of the author and Herr \cite{Candy2016}, used a similar observation, although phrased in terms of $U^2_{\Phi_j}$ atoms, to prove a unit scale bilinear restriction in the adapted function space $V^2_{\Phi_j}$. The work \cite{Candy2016} provided a strong motivation to consider the case of non-unit scale waves,  however it is important to note that the argument used in \cite{Candy2016} follows the combinatorial Kakeya  type approach to bilinear restriction estimates, in contrast to the energy type approach used here. In particular, it is also possible to prove estimates for vector valued waves via the combinatorial approach.
\end{remark}

Theorem \ref{thm-main} is optimal up to endpoints, in the sense that the restriction $\frac{1}{q} +\frac{n+1}{2r} \les \frac{n+1}{2}$ is necessary, and the dependence on $\Phi_j$ and $\mb{d}_0$ in (\ref{eqn-intro bilinear data in ball}) cannot be improved, see Section \ref{sec:counter example}. In particular, if the Fourier support of $u$ or $v$ is contained in a ball of radius $\mb{d}_0$, then the dependence of (\ref{eqn-intro bilinear data in ball}) on $\mb{d}_0$ is sharp. However, if the support of $\widehat{u}$ (or $\widehat{v}$) is \emph{not} well approximated by a ball, then certain improvements are possible. One possibility is the following. Given sets $\Lambda_j^* \subset \Lambda_j$, we define the quantity
     \begin{equation} \label{eqn-defn of d as surface measure} \mb{d}[\Lambda_1^*, \Lambda_2^*] = \sup_{\mathfrak{h}=(a,h)\in \RR^{1+n}} \Big( \sigma_{\Sigma_1(\mathfrak{h})}\big[ \Sigma_1(\mathfrak{h}) \cap \Lambda_1^* \cap \big( h - \Lambda_2^* \big) \big]\Big)^{\frac{1}{n-1}} \end{equation}
where $\sigma_{\Sigma_j(\mathfrak{h})}$ denotes the (induced) Lebesgue surface measure on the surface $\Sigma_j(\mathfrak{h})\subset \RR^n$. It easy to check that, under the transversality assumption (\ref{eqn-main trans assump}), we have the bound ${\mb{d}[\Lambda_1^*, \Lambda_2^*] \lesa \min\{ \diam(\Lambda_1^*), \diam(\Lambda_2^*)\}}$. However $\mb{d}[\Lambda_1^*, \Lambda_2^*]$ can potentially be much smaller as it only requires the \emph{intersection} of $\supp \widehat{u}$,  the surface $\Sigma_1(\mathfrak{h})$, and translates of $\supp \widehat{v}$ to be small. The quantity $\mb{d}[\Lambda_1^*, \Lambda_2^*]$ arises precisely in the bilinear $L^2_{t,x}$ estimate. In fact, under the transversality assumption (\ref{eqn-main trans assump}) we have for every $\Phi_1$-wave $u$, and $\Phi_2$-wave $v$,
        $$ \| u v \|_{L^2_{t,x}} \les ( \mb{C}_0 \mc{V}_{max})^{-\frac{1}{2}} (\mb{d}[\supp \widehat{u}, \supp \widehat{v}])^{\frac{n-1}{2}} \| u \|_{L^\infty_t L^2_x} \| v \|_{L^\infty_t L^2_x},$$
see Theorem \ref{thm-classical bilinear L2 estimate} below. Exploiting this observation leads to the following improvement of Theorem \ref{thm-main}.

\begin{theorem}\label{thm-main small scale}
Let $n\g2$, $\mb{C}_0>0$,  $1\les q, r \les 2$, and $\frac{1}{q} +\frac{n+1}{2r} < \frac{n+1}{2}$. There exists a constant $C>0$, such that for any $ \mb{d}_0>0$, any open sets $\Lambda_1, \Lambda_2 \subset \RR^n$, any phases $\Phi_1$ and $\Phi_2$ satisfying Assumption \ref{assump-main} with $\mc{H}_2 \les \mc{H}_1$, and any $\Phi_1$-waves $u$, and $\Phi_2$-waves $v$ satisfying the support conditions
    $$\supp \widehat{u} + \mb{d}_0 \subset \Lambda_1, \qquad \supp \widehat{v} + \mb{d}_0 \subset \Lambda_2, \qquad \mb{d}\big[\supp \widehat{u} + \mb{d}_0, \supp \widehat{v}+\mb{d}_0\big] \les \frac{\mb{d}_0}{\mb{C}_0}, $$
we have
		$$ \big\| uv  \big\|_{L^q_t L^r_x(\RR^{1+n})}\les C\mb{d}_0^{n+1 - \frac{n+1}{r} - \frac{2}{q}} \,\mc{V}_{max}^{\frac{1}{r}-1} \,\mc{H}_{1}^{1-\frac{1}{q} - \frac{1}{r}}\, \Big( \frac{\mc{H}_{1}}{\mc{H}_{2}}\Big)^{\frac{1}{q} - \frac{1}{2}}\| u \|_{L^\infty_t L^2_x} \| v \|_{L^\infty_t L^2_x}. $$
\end{theorem}

\begin{remark}
  It is clear that if $\min\{ \diam(\supp \widehat{u}), \diam(\supp \widehat{v}) \} \les \mb{d}_0$, then since one of $\supp \widehat{u} + \mb{d}_0$ or $\supp \widehat{v} + \mb{d}_0$ is again contained in a ball of radius $\mb{d}_0$, a computation gives
         $$ \mb{d}\big[\supp \widehat{u} + \mb{d}_0, \supp \widehat{v}+\mb{d}_0\big] \lesa \mb{d}_0.$$
  In particular, after potentially choosing $\mb{C}_0$ slightly smaller, we see that Theorem \ref{thm-main small scale} implies Theorem \ref{thm-main}. We note here that the $\mb{d_0}$ thickening of the Fourier supports, namely $\supp \widehat{u} + \mb{d}_0$ and $\supp \widehat{v} + \mb{d}_0$, arises as the induction on scales argument used in the proof of Theorem \ref{thm-main small scale} involves numerous decompositions of $u$ and $v$ into wave packets, each of which may have slightly larger Fourier support. Thus we require some room to increase the Fourier support in the induction argument, which eventually manifests itself in the assumptions in Theorem \ref{thm-main small scale} applying to a $\mb{d}_0$ neighbourhood of the Fourier supports of $u$ and $v$.
\end{remark}

We now give a number a number of consequences of Theorem \ref{thm-main small scale}. Namely, we state a version of Theorem \ref{thm-main small scale} in the adapted function space $U^2_{\Phi}$, and give a concrete application to the case of elliptic and Klein-Gordon waves. Moreover, we prove a new refined Strichartz estimate for the Klein-Gordon equation.

\subsection{Bilinear Restriction in Adapted Function spaces} Theorem \ref{thm-main small scale} has a number of applications to linear PDE, for instance, we give a new bilinear restriction estimate for wave/Klein-Gordon interactions, together a refined Strichartz type estimate for the Klein-Gordon equation, see Theorem \ref{thm-bilinear small scale KG} and Theorem \ref{thm-refined strichartz} below. However, in applications to \emph{nonlinear} PDE, it is useful to have a version of Theorem \ref{thm-main small scale} which holds for more general functions than just $\Phi_j$-waves. Often the function spaces which arise in applications satisfy some version of the \emph{transference principle}. This principle roughly asserts that if we have an estimate for $\Phi_j$-waves, then we immediately deduce that the estimate also holds in some more general function space. For instance the well-known $X^{s,b}$ type spaces satisfy the transference principle, essentially since elements of $X^{s, b}$ type spaces can be written as averages of free waves. On the other hand, other function spaces of interest, such as the adapted function spaces $U^p_{\Phi_j}$ and $V^p_{\Phi_j}$, do not generally satisfy such a strong property, and significant additional work can be required to extend estimates for free waves, to estimates in $U^p_{\Phi_j}$, see for instance \cite{Candy2016} for further discussion.

However, the adapted function space $U^2_{\Phi_j}$ \emph{does} satisfy a slightly weaker transference type principle: if an estimate holds for \emph{vector valued} waves, then it also holds in $U^2_{\Phi_j}$. Consequently, the fact that Theorem \ref{thm-main small scale} holds for vector valued waves, immediately implies that bilinear restriction type estimates for functions in $U^2_{\Phi_j}$, see the proof of Corollary \ref{cor-main Up version} below or \cite[Remark 5.2]{Candy2016} for details of this argument. Bilinear restriction type estimates in the adapted function spaces $U^p_{\Phi_j}$ and $V^p_{\Phi_j}$ first appeared in work of Sterbenz-Tataru \cite[Lemma 5.7]{Sterbenz2010}, and the full bilinear range for general phases at unit scale was obtained recently in \cite{Candy2016}. The atomic function spaces $U^p_{\Phi_j}$ have proven to be extremely useful in obtaining endpoint well-posedness results for dispersive PDE, see for instance \cite{Koch2005, Koch2007, Hadac2009} and the references therein. \\

Before we can state our next result, we require the definition of $U^p_{\Phi_j}$. Let $1<p<\infty$. We say that $\phi$ is a $U^p_{\Phi_j}$ \emph{atom}, if there exists a finite partition $\mc{I}$ of $\RR$, and a collection of $\Phi_j$-waves $(\phi_I)_{I\in \mc{I}}$ such that $\phi = \sum_{I \in \mc{I}} \ind_I(t) \phi_I$ and
        $$ \bigg( \sum_{I \in \mc{I}} \| \phi_I \|_{L^\infty_t L^2_x}^p \bigg)^{\frac{1}{p}} \les 1.$$
The atomic Banach space $U^p_{\Phi_j}$ is then defined as
        $$ U^p_{\Phi_j} = \Big\{ \sum_m c_m \phi_m \,\,\Big| \,\, (c_m)_{m \in \NN} \in \ell^1(\NN), \,\, \text{$\phi_m$ are $U^p_{\Phi_j}$ atoms }\Big\}$$
with the induced norm
    $$ \| u \|_{U^p_{\Phi_j}} = \inf_{ u = \sum_m c_m \phi_m} \sum_m |c_m|$$
where the inf is over all representations of $u$ in terms of $U^p_{\Phi_j}$ atoms $\phi_m$. Note that functions in $U^p_{\Phi_j}$ take values in $\ell^2(\ZZ)$, in particular, every $\Phi_j$-wave with norm of size one, is a (trivial) example of a $U^p_{\Phi_j}$ atom.

\begin{corollary}\label{cor-main Up version}
Let $n\g2$, $\mb{C}_0>0$,  $1\les q, r \les 2$, and $\frac{1}{q} +\frac{n+1}{2r} < \frac{n+1}{2}$. There exists a constant $C>0$, such that for any $ \mb{d}_0>0$, any open sets $\Lambda_1, \Lambda_2 \subset \RR^n$, any phases $\Phi_1$ and $\Phi_2$ satisfying Assumption \ref{assump-main} with $\mc{H}_2 \les \mc{H}_1$, and any $u \in U^2_{\Phi_1}$, $v \in U^2_{\Phi_2}$ satisfying the support conditions
    $$\supp \widehat{u} + \mb{d}_0 \subset \Lambda_1, \qquad \supp \widehat{v} + \mb{d}_0 \subset \Lambda_2, \qquad \mb{d}\big[\supp \widehat{u} + \mb{d}_0, \supp \widehat{v}+\mb{d}_0\big] \les \frac{\mb{d}_0}{\mb{C}_0}, $$
we have
		$$ \big\| uv  \big\|_{L^q_t L^r_x(\RR^{1+n})}\les C\mb{d}_0^{n+1 - \frac{n+1}{r} - \frac{2}{q}} \,\mc{V}_{max}^{\frac{1}{r}-1} \,\mc{H}_{1}^{1-\frac{1}{q} - \frac{1}{r}}\, \Big( \frac{\mc{H}_{1}}{\mc{H}_{2}}\Big)^{\frac{1}{q} - \frac{1}{2}}\| u \|_{U^2_{\Phi_1}} \| v \|_{U^2_{\Phi_2}}. $$
\begin{proof}
We exploit the transference type principle mentioned above. It is enough to consider the case where $u= \sum_{I \in \mc{I}} \ind_I(t) u_I$ and $v=\sum_{J \in \mc{J}} \ind_J(t) v_I$ are atoms. Define the vector valued waves $U = (u_I)_{I \in \mc{I}}$ and $V = (v_J)_{J\in \mc{J}}$. Then since $u$ and $v$ are atoms, $U$ is a $\Phi_1$-wave, $V$ is a $\Phi_2$-wave, we have the energy estimate $\| U \|_{L^\infty_t L^2_x} , \| V \|_{L^\infty_t L^2_x} \les 1$, and the pointwise bounds
    \begin{equation}\label{eqn:atom bounded by vector} \sum_{I \in \mc{I} } \ind_I(t) |u_I| \les \sum_{I \in \mc{I}} \ind_I(t) |U| \les |U|, \qquad \sum_{J \in \mc{J}} \ind_J(t) |v_J| \les |V|.\end{equation}
Therefore corollary follows from an application of Theorem \ref{thm-main small scale}.
\end{proof}
\end{corollary}

The bound \eref{eqn:atom bounded by vector} is somewhat crude, and in fact Corollary \ref{cor-main Up version} can be improved by further exploiting the time localisation of the $U^p_{\Phi_j}$ atoms directly in the proof of Theorem \ref{thm-main small scale}. More precisely, by refining the argument used to prove Theorem \ref{thm-main small scale}, and interpolating with an ``linear'' version of the bilinear $L^2_{t,x}$ inside the induction on scales argument, we obtain the following improvement of Corollary \ref{cor-main Up version}.

\begin{theorem}\label{thm-main Up version II}
Let $n\g2$, $\mb{C}_0>0$,  $1\les q, r \les 2$ and $\frac{1}{q} +\frac{n+1}{2r} < \frac{n+1}{2}$. Let $ \frac{2}{(n+1)q}< \frac{1}{b} \les \frac{1}{a} \les \frac{1}{2}$ and $\frac{1}{a} + \frac{1}{b} \g \frac{1}{\min\{q, r\}}$. There exists a constant $C>0$, such that for any $ \mb{d}_0>0$, any open sets $\Lambda_1, \Lambda_2 \subset \RR^n$, any phases $\Phi_1$ and $\Phi_2$ satisfying Assumption \ref{assump-main} with $\mc{H}_2 \les \mc{H}_1$, and any $u \in U^a_{\Phi_1}$, $v \in U^b_{\Phi_2}$ satisfying the support conditions
    $$\supp \widehat{u} + \mb{d}_0 \subset \Lambda_1, \qquad \supp \widehat{v} + \mb{d}_0 \subset \Lambda_2, \qquad \mb{d}\big[\supp \widehat{u} + \mb{d}_0, \supp \widehat{v}+\mb{d}_0\big] \les \frac{\mb{d}_0}{\mb{C}_0}, $$
we have
\begin{equation}\label{eqn:improve Up bilinear bound}\begin{split}
    \big\| uv  &\big\|_{L^q_t L^r_x(\RR^{1+n})} \\
    &\les C\mb{d}_0^{n+1 - \frac{n+1}{r} - \frac{2}{q}} \,\mc{V}_{max}^{\frac{1}{r}-1} \,\mc{H}_{1}^{1-\frac{1}{q} - \frac{1}{r}}\, \Big( \frac{\mc{H}_{1}}{\mc{H}_{2}}\Big)^{\frac{1}{q} - \frac{1}{2} + (n+1)(\frac{1}{2}-\frac{1}{a})} \Big( \frac{\mu^n \mc{V}_{max}}{\mb{d}_0^{n+1} \mc{H}_1}\Big)^{(1-\frac{1}{r}-\frac{1}{b})_+}\| u \|_{U^a_{\Phi_1}} \| v \|_{U^b_{\Phi_2}}
\end{split}
\end{equation}
where $\mu = \min\{ \diam( \supp \widehat{u}), \diam( \supp \widehat{v}), \mb{d}_0\}$ and $s_+ = s$ if $s>0$, and zero otherwise.
\end{theorem}

The bound \eref{eqn:improve Up bilinear bound} has the useful consequence that we may place $v$ into the weaker $V^2_\Phi$ space \emph{without} the standard high-low frequency loss. To make this claim more precise, define the space $V^2_{\Phi}$ as all right continuous functions such that
    $$ \| u \|_{V^2_{\Phi}} =\| u \|_{L^\infty_t L^2_x} +  \sup_{(t_k) \in \mc{Z}} \Big( \sum_{k \in \ZZ} \| e^{-it_k\Phi(-i\nabla)}u(t_k) - e^{-it_{k-1}\Phi(-i\nabla)} u(t_{k-1}) \|_{L^2_x}^2 \Big)^\frac{1}{2}<\infty$$
where $\mc{Z} = \{ (t_j)_{j \in \ZZ} \mid t_j \in \RR \text{ and } t_j<t_{j+1} \}$. Thus functions in $V^2_\Phi$ have finite quadratic variation along the flow $e^{it \Phi(-i\nabla)}$. An application of the continuous embedding $V^2_\Phi \subset U^a_\Phi$ for $a>2$ \cite[Lemma 6.4]{Koch2005}, together with \eref{eqn:improve Up bilinear bound} implies that we then have
         $$\big\| uv  \big\|_{L^q_t L^r_x(\RR^{1+n})}\les C\mb{d}_0^{n+1 - \frac{n+1}{r} - \frac{2}{q}} \,\mc{V}_{max}^{\frac{1}{r}-1} \,\mc{H}_{1}^{1-\frac{1}{q} - \frac{1}{r}}\, \Big( \frac{\mc{H}_{1}}{\mc{H}_{2}}\Big)^{\frac{1}{q} - \frac{1}{2}}\| u \|_{U^2_{\Phi_1}} \| v \|_{V^2_{\Phi_2}}$$
provided that $1<q \les 2$, $1<r<2$, and $\frac{1}{q} + \frac{n+1}{2r}<\frac{n+1}{2}$. This bound is precisely the same as the case of free solutions. In particular, we avoid the high-low  frequency loss which would come from simply deducing a $V^2_\Phi$ bound from Corollary \ref{cor-main Up version} via an interpolation argument. In the special case $q=r=2$, we get a loss, and only obtain
         $$\big\| uv  \big\|_{L^2_{t,x}(\RR^{1+n})}\les C\mb{d}_0^{-\frac{n-1}{2}} \,\mc{V}_{max}^{-\frac{1}{2}} \Big( \frac{\mu^n \mc{V}_{max}}{\mb{d}_0^{n+1} \mc{H}_1}\Big)^{\epsilon}\| u \|_{U^2_{\Phi_1}} \| v \|_{V^2_{\Phi_2}}.$$
This bilinear $L^2_{t,x}$ bound is particularly interesting in the case of the wave or Klein-Gordon equation, as it allows us to place the high frequency wave into the weaker $V^2_{\Phi_2}$ space without losing any high-frequency derivatives.

\subsection{The Elliptic Case}  We say that a phase $\Phi:\Lambda \rightarrow \RR$ is \emph{elliptic} on $\Lambda$, if for all $\xi \in \Lambda$ and $v\in \RR^n$  we have
 $$ \big| \big( \nabla^2 \Phi(\xi) v\big) \cdot v \big|\approx \| \nabla^2 \Phi\|_{L^\infty(\Lambda)} |v|^2.$$
Equivalently, the eigenvalues of $\nabla^2 \Phi$ all have the same sign, and are essentially of the size $\|\nabla^2 \Phi_j\|_{L^\infty(\Lambda)}$. A typical example of an elliptic phase is the Schrodinger phase $\Phi = \frac{1}{2}|\xi|^2$, or the Klein-Gordon phase $(m^2 + |\xi|^2)^\frac{1}{2}$ in the region $|\xi| \ll m$. Bilinear restriction estimates for elliptic phases was recently exploited by Stovall \cite{Stovall2017} to deduce new results for the linear restriction problem. Applying Theorem \ref{thm-main small scale} to the case of elliptic phases, gives the following improvement to \cite[Theorem 2.1]{Stovall2017}.

\begin{theorem}\label{thm-elliptic bilinear}
Let $1\les q, r \les 2$, $\frac{1}{q} +\frac{n+1}{2r} < \frac{n+1}{2}$, and $\mb{d}_0>0$. Let $\Lambda_j \subset \RR^n$ be convex, and $\Phi_j \in C^{5n}(\Lambda_j)$ be elliptic phases such that $\mc{H}_2 \les \mc{H}_1$,
        $$   \diam(\Lambda_j) \lesa  \frac{\mc{V}_{max}}{\mc{H}_j}, \qquad \diam(\Lambda_j) \lesa \min_{2\les m \les 5n} \Big( \frac{\mc{H}_j}{\| \nabla^m \Phi_j \|_{L^\infty(\Lambda_j)} } \Big)^{\frac{1}{m-2}}, $$
and for all $\xi\in \Lambda_1$, $\eta \in \Lambda_2$
        $$ \big| \nabla \Phi_1(\xi) - \nabla \Phi_2(\eta)\big| \gtrsim \mc{V}_{max}.$$
Then for any $\Phi_1$-waves $u$, and $\Phi_2$-waves $v$ satisfying the support conditions
    $$\supp \widehat{u} + \mb{d}_0 \subset \Lambda_1, \qquad \supp \widehat{v} + \mb{d}_0 \subset \Lambda_2, \qquad \mb{d}\big[\supp \widehat{u} + \mb{d}_0, \supp \widehat{v}+\mb{d}_0\big] \les \frac{\mb{d}_0}{\mb{C}_0}, $$
we have
		$$ \big\| uv  \big\|_{L^q_t L^r_x(\RR^{1+n})}\lesa \mb{d}_0^{n+1 - \frac{n+1}{r} - \frac{2}{q}} \,\mc{V}_{max}^{\frac{1}{r}-1} \,\mc{H}_{1}^{1-\frac{1}{q} - \frac{1}{r}}\, \Big( \frac{\mc{H}_{1}}{\mc{H}_{2}}\Big)^{\frac{1}{q} - \frac{1}{2}}\| u \|_{L^\infty_t L^2_x} \| v \|_{L^\infty_t L^2_x}. $$
\end{theorem}
\begin{remark}
The precise assumptions used in the Theorem \ref{thm-elliptic bilinear} can be weakened slightly (in particular the convexity assumption). See the proof of Lemma \ref{lem:local version of A1} below.
\end{remark}

\subsection{Bilinear Restriction for the Klein-Gordon equation} We now turn to the case of the Klein-Gordon equation. Let
    $$ \lr{\xi}_m = (m^2 + |\xi|^2)^\frac{1}{2}.$$
Note that when $m=0$ this is simply the wave phase $|\xi|$. If we have the support assumptions
    $$  \supp \widehat{f} \subset \{  |\xi|  \approx  \lambda, \theta(\xi, \omega_0)\ll  \alpha \},  \qquad \supp \widehat{g} \subset \{ |\xi|  \approx  \lambda, \theta(\xi,\omega_1) \ll \alpha \}$$
with $\theta(\omega_0, \omega_1) \approx \alpha$ (here $\theta(x,y)$ denotes the angle between vectors $x, y \in \RR^n$) then it is known that for $\frac{1}{q} + \frac{n+1}{2r} < \frac{n+1}{2}$, $1\les q, r\les 2$, and $\lambda\g \mu$ we have the bilinear estimate for the wave equation
    $$ \| e^{ it |\nabla|} f e^{it |\nabla|} g \|_{L^q_t L^r_x} \lesa \alpha^{n-1 - \frac{n-1}{r}  - \frac{2}{q}} \mu^{n-\frac{n}{r} - \frac{1}{q}} \Big( \frac{\lambda}{\mu}\Big)^{\frac{1}{q} - \frac{1}{2}} \| f\|_{L^2} \| g \|_{L^2}, $$
see \cite{Tao2001b, Tataru2003, Lee2008, Lee2008a}. An application of Theorem \ref{thm-main small scale} gives the following Klein-Gordon counterpart.

\begin{theorem}[Bilinear extension for K-G]\label{thm-bilinear small scale KG}
Let $1\les q, r \les2$, $\frac{1}{q} + \frac{n+1}{2r}<\frac{n+1}{2}$, $m_1, m_2 \g 0$, and $\lambda \g \mu >0$. Let $0< \alpha \les 1$, and suppose we have $\xi_0, \eta_0 \in \RR^n$ such that $\lr{\xi_0}_{m_1} \approx \lambda$, $\lr{\eta_0}_{m_2} \approx \mu$, and
    $$ \frac{ \big| m_2 |\xi_0| - m_1 |\eta_0| \big|}{\lambda \mu} + \Big( \frac{|\xi_0| |\eta_0| \mp \xi_0 \cdot \eta_0}{\lambda \mu} \Big)^\frac{1}{2} \approx \alpha. $$
Define $\beta = ( \frac{m_1}{\alpha \lambda} + \frac{m_2}{\alpha \mu} + 1)^{-1}$. If
    $$ \supp \widehat{f} \subset \big\{ \big| |\xi| - |\xi_0| \big| \ll  \beta \lambda, \,\,  \big( |\xi| |\xi_0| - \xi \cdot \xi_0 \big)^\frac{1}{2} \ll \alpha  \lambda \big\}, \qquad \supp \widehat{g} \subset \big\{ \big| |\xi| - |\eta_0| \big| \ll  \beta  \mu, \,\,  \big( |\xi| |\eta_0| - \xi \cdot \eta_0 \big)^\frac{1}{2} \ll \alpha  \mu \big\}$$
then we have the bilinear estimate
    $$ \big\| e^{it\lr{-i\nabla}_{m_1}}f e^{\pm it \lr{-i\nabla}_{m_2}} g \big\|_{L^q_t L^r_x} \lesa \alpha^{n-1 - \frac{n-1}{r}  - \frac{2}{q}} \beta^{1-\frac{1}{r}}  \mu^{n-\frac{n}{r} - \frac{1}{q}} \Big( \frac{\lambda}{\mu}\Big)^{\frac{1}{q} - \frac{1}{2}} \|f \|_{L^2_x} \| g \|_{ L^2_x} $$
where the implied constant is independent of $m_1, m_2$.
\end{theorem}

Note that Theorem \ref{thm-bilinear small scale KG} contains both wave like regimes (when $\beta \approx 1$), and Schr\"{o}dinger like regimes (when $\beta \ll 1$). In particular, if $\beta = 1$, then the support assumptions on $f$ and $g$ are the same as in the case of the wave equation, while if $\beta \approx \alpha $, then the support assumptions reduce to simply requiring that
        $$ \supp \widehat{f} \subset \{ |\xi - \xi_0|  \ll  \alpha \lambda \}, \qquad \supp \widehat{g} \subset \{ |\xi - \eta_0| \ll  \alpha  \mu, \}$$
which matches the standard assumptions used in the Schr\"odinger case (see for instance Theorem \ref{thm-elliptic bilinear} above).

\begin{remark}\label{rem:bilinear atomic KG}
If we apply the atomic bilinear restriction estimate, Theorem \ref{thm-main Up version II}, in place of Theorem \ref{thm-main}, under the same assumptions as those in Theorem \ref{thm-bilinear small scale KG}, we arrive at the stronger bounds
        $$ \| u v \|_{L^q_t L^r_x} \lesa \alpha^{n-1 - \frac{n-1}{r}  - \frac{2}{q}} \beta^{1-\frac{1}{r}}  \mu^{n-\frac{n}{r} - \frac{1}{q}} \Big( \frac{\lambda}{\mu}\Big)^{\frac{1}{q} - \frac{1}{2}}  \| u \|_{U^2_{\lr{\nabla}_{m_1}}} \| v \|_{V^2_{\pm \lr{\nabla}_{m_2}}}.$$
In particular, in the special case $q=r=2$, we may place $v \in V^2_{\pm \lr{\nabla}_{m_2}}$ \emph{without} a high frequency loss.
\end{remark}

\subsection{A Refined Strichartz Estimates for the Klein-Gordon Equation} The standard $H^{\frac{1}{2}}$ Strichartz estimate for the Klein-Gordon equation states that
        $$ \| e^{ it \lr{\nabla}} f \|_{L^{2 \frac{n+1}{n-1}}_{t,x}(\RR^{1+n})} \lesa \| f \|_{H^{\frac{1}{2}}}.$$
It is of interest to understand the concentration properties of solutions which come close to maximising this inequality, for instance, understanding these concentrating solutions is a key step in obtaining a \emph{profile decomposition}. In the case of the Schr\"{o}dinger equation, this concentration type property is a consequence of a refined Strichartz estimate introduced by Bourgain \cite{Bourgain1998}  and extended in work of Moyua-Vargas-Vega \cite{Moyua1996, Moyua1999}. A version of the refined Strichartz estimate for the Klein-Gordon equation was proved by Killip-Stovall-Visan \cite{Killip2012b} with data in the Klein-Gordon regime (in other words, the bound included a loss of derivatives when compared to the wave regime).  Here we use Theorem \ref{thm-main} together with an argument used by Ramos \cite{Ramos2012} to obtain a refined Strichartz estimate for the Klein-Gordon equation without any derivative loss.

Before we state the refined Strichartz estimate for the Klein-Gordon equation, we require some notation. Given $\lambda \g 1$, and $0<\alpha <1$, we define $\mc{A}_{\lambda, \alpha}$ to the collection of sets $A \subset \RR^n$ of the form
    $$ A = A(\xi_0) = \big\{ \lr{\xi}\approx \lambda, \,\,  \big| |\xi| - |\xi_0| \big| \ll \tfrac{\alpha \lambda}{1+\alpha \lambda} \lambda, \,\, \big( |\xi| |\xi_0| - \xi \cdot \xi_0 \big)^\frac{1}{2} \ll \alpha \lambda \big\}$$
where the points $\xi_0 \in \RR^n$ satisfy $\lr{\xi_0} \approx \lambda$ and are chosen to ensure that the sets $A \in \mc{A}_{\lambda, \alpha}$  form a finitely overlapping cover of the annulus/ball $\{\lr{\xi}\approx \lambda \}$. In the wave like region $\alpha \gtrsim \frac{1}{\lambda}$, the elements of $\mc{A}_{\lambda, \alpha}$ are radial sectors of the annuli $\lr{\xi} \approx \lambda$, while in the Klein-Gordon (or Schr\"odinger) like case $\alpha \ll \frac{1}{\lambda}$, the elements of $\mc{A}_{\lambda, \alpha}$ are angular sectors with some additional radial restrictions which degenerate to cubes if $\lambda \lesa 1$.

\begin{theorem}[Refined Strichartz]\label{thm-refined strichartz} Let $n \g 2$. There exists $0<\theta<1$ and $1<r<2$ such that
        $$ \big\| e^{ it \lr{\nabla}} f \big\|_{L^{2\frac{n+1}{n-1}}_{t,x}(\RR^{1+n})} \lesa \Big( \sup_{\lambda \in 2^{\NN}, \alpha \in 2^{-\NN}} \sup_{A \in \mc{A}_{\lambda, \alpha}} \Big( \tfrac{\alpha \lambda}{1+\alpha \lambda}\Big)^{\frac{1}{n+1}}  \lambda^\frac{1}{2} |A|^{\frac{1}{2} - \frac{1}{r}} \| \widehat{f} \|_{L^r_\xi(A)} \Big)^\theta \| f \|_{H^{\frac{1}{2}}}^{1-\theta}. $$
\end{theorem}

This estimate is the Klein-Gordon counterpart to the estimate for the wave equation obtain by Ramos \cite{Ramos2012}. Applications of Theorem \ref{thm-refined strichartz} to the $H^{\frac{1}{2}}$ profile decomposition for the Klein-Gordon equation will appear elsewhere.

\subsection{Outline} In Section \ref{sec-applications} we first show that the global condition (\textbf{A1}) is equivalent to the local condition \eref{eqn:intro local cond}, and adapt an argument from \cite{Candy2016} to give a simplification of the transversality/curvarture condition ($\mb{A1}$). We then give the proof of the main consequences of Theorem \ref{thm-main small scale}, namely the elliptic case Theorem \ref{thm-elliptic bilinear}, the Klein-Gordon case, Theorem \ref{thm-bilinear small scale KG}, and the refined Strichartz estimate, Theorem \ref{thm-refined strichartz}.

In Section \ref{sec:counter example} we give a counter example which shows that the dependence on the phases $\Phi_j$, and the frequency parameter $\mb{d}_0$ in Theorem \ref{thm-main} is sharp. The counter example is based on constructing waves $u$ and $v$ which are sums of wave packets concentrating on an $\epsilon^{-1}$ neighbourhood of the  space-time rectangle
    $$ \big\{ \big(s+s', - s \nabla \Phi_1(\xi_0) - s' \nabla \Phi_2(\eta_0)\big) \,\,\big|\,\, 0\les s \les (\epsilon^{2} \mc{H}_1)^{-1}, \,\, 0\les s' \les (\epsilon^2 \mc{H}_2)^{-1} \big\}$$
with $\xi_0 \in \Lambda_1$ and $\eta_0 \in \Lambda_2$.

In Section \ref{sec-ind on scales}, we run the induction on scales argument, which reduces the proof of Theorem \ref{thm-main small scale} to showing that we may control the $L^q_t L^r_x$ norm on a cube of diameter $R$, by a cube of diameter $\frac{R}{2}$. The induction argument applied here is slightly different to that used in Lee-Vargas \cite{Lee2008}, and in particular it is here where we are able to avoid the $\epsilon$ loss in the scale factors that occurred in previous works.

In Section \ref{sec-the localisation argument} we apply localisation arguments to further reduce the proof of Theorem \ref{thm-main small scale} to obtaining a key wave table type decomposition, Theorem \ref{thm-general wave table decomposition}, which decomposes the product $u v$ into a term which is concentrated at scale $\frac{R}{2}$, together with a term that satisfies an improved bilinear estimate. The proof of the decomposition contained in Theorem \ref{thm-general wave table decomposition} is the key step in the proof of Theorem \ref{thm-main small scale}.

 Section \ref{sec-wave packets} contains the general wave packet decomposition used in the present article, while, in Section \ref{sec-geometric consequences of phase}, closely following \cite{Candy2016}, we give the key geometric consequences of (\textbf{A1}). The energy estimates across transverse surfaces we require are given in Section \ref{sec-energy estimates across C_j}. These estimates rely heavily on the geometric consequences of the curvature and transversality assumption (\textbf{A1}) contained in Section \ref{sec-geometric consequences of phase}.

 In Section \ref{sec-wave table construct} we give the main step in the wave table construction, namely, following Wolff \cite{Wolff2001} and Tao \cite{Tao2001b}, we use the wave packet decomposition, together with an energy argument, to decompose $u$ into wave packets such that $v$ restricted to the corresponding tube is concentrated on a small cube.

  In Section \ref{sec-proof of key bilinear est} we use the construction in Section \ref{sec-wave table construct} to give the proof of Theorem \ref{thm-general wave table decomposition} and construct the general wave tables that are required in the proof of the localised bilinear restriction estimate.

  Finally, in Sections \ref{sec:atomic wave tables} and \ref{sec:proof of atomic wave tables}, we give the proof of Theorem \ref{thm-main Up version II}. This relies on developing an ``atomic'' version of the wave table construction used in Section \ref{sec-wave table construct}, together with an additional interpolation argument with a ``linear'' version of the classical bilinear $L^2_{t,x}$ estimate.

\subsection{Notation} Throughout this article, we use the notation $a \lesa b$ if $a \les C B$ for some constant $C>0$ which may depend \emph{only} on $\mb{C}_0$, $q$, $r$, and the dimension $n$. In particular, all implied constants will otherwise be independent of the phases $\Phi_j$. Similarly we write $a \approx b$ if $a\lesa b$ and $b \lesa a$, and $a \gg b$ if $a \g C b$ with  $C \g 100^n$.  For $a \in \RR$, we let $a_+ = a$ if $a>0$, and $a_+=0$ otherwise. \\

Given a set $\Omega \subset \RR^n$ and a vector $h \in \RR^n$, we let $\Omega + h = \{ x + h \mid x \in \Omega \}$ denote the translation of $\Omega$ by $h$. At a risk of causing some minor confusion, for a positive constant $c>0$ we let $\Omega + c = \{ x + y\mid x\in \Omega, \,|y|<c\}$ be the Minkowski sum of $\Omega$ and the ball $\{|x|< c\}$. We also let $\diam(\Omega)= \sup_{x, y \in \Omega} |x-y|$. For a function $f \in L^1(\RR^n)$, we let $\widehat{f}(\xi) = \int_{\RR^n} f(x) e^{ - i x \cdot \xi} dx$ denote the spatial Fourier transform. Similarly, if $u \in L^1(\RR^{1+n})$, we define the space-time Fourier transform of $u$ as $\tilde{u}(\tau, \xi) = \int_{\RR^{1+n}} u(t,x) e^{ - i (t,x) \cdot (\tau, \xi) } \,dt\,dx$ and let $\widehat{u}(t)$ be the Fourier transform of $u$ in the spatial variable $x \in \RR^n$. \\

Let $R \g 1$ and $0<\epsilon\les 1$. The constant $R$ will denote the large space-time scale and $0<\epsilon\les 1$ will be a small fixed parameter used to control the various error terms that arise. Given a scale $r>0$, we decompose phase space into the grid $ \mc{X}_{r, \epsilon} = \frac{r}{\epsilon^2} \ZZ^n \times \frac{1}{r} \ZZ^n$. Given a point $\gamma = (x_0, \xi_0) \in \mc{X}_{r, \epsilon}$ in phase-space, we let $x(\gamma) = x_0$ and $\xi(\gamma)= \xi_0$ denote the projections onto the first and second components respectively.\\

Associated to a phase $\Phi_j$, and a scale $R$, we let $r_j = ( \mc{H}_j R)^\frac{1}{2}$ denote the scale of the corresponding wave packets. Define
        $$\Gamma_j = \big\{ (x_0, \xi_0) \in \mc{X}_{r_j, \epsilon} \,\, \big| \,\, \xi_0 + \tfrac{1}{4}\mb{d}_0 \subset \Lambda_j  \big\}.$$
Thus $\Gamma_j$ contains those phase space points inside $\Lambda_j$, which are at least $\frac{1}{4} \mb{d}_0$ from the boundary of $\Lambda_j$. Given a point $\gamma_j=(x_0, \xi_0) \in \Gamma_j$ we define the associated tubes
    $$T_{\gamma_j} = \{ (t,x) \in \RR^{1+n} \mid |x-x_0 + t \nabla \Phi_j(\xi_0) | \les r_j \}. $$
The geometric condition (\textbf{A1})  plays an important role in the proof of Theorem \ref{thm-main} by restricting the intersections of various tubes to certain transverse hypersurfaces. However, to exploit the curvature hypothesis in (\textbf{A1}), we need to restrict points in $\Gamma_j$ to those  close to $\Sigma_j(\mathfrak{h})$. To this end, given $\mathfrak{h} \in \RR^{1+n}$ and a set $\Lambda \subset \RR^n$, we define
	$$ \Gamma_j(\mathfrak{h}, \Lambda) = \big\{ \gamma_j \in \Gamma_j \,\,\big| \,\, \xi(\gamma_j) \in \Sigma_j(\mathfrak{h})\cap \Lambda  + \tfrac{C}{r_j} \,\, \big\}$$
where the constant $C$ will depend on the dimension $n$ and $\mb{C}_0$, but will otherwise be independent of the phase $\Phi_j$. We also define the conic hypersurface
        $$ \mc{C}_j(\mathfrak{h}) = \big\{ \big(s, - s\nabla \Phi_j(\xi) \big) \,\,\big| \,\,s\in \RR,  \xi \in \Sigma_j(\mathfrak{h}) \, \big\}.$$\\

All cubes in this article are oriented parallel to the coordinate axis. Let $Q$ be a cube side length $R$, and take a subscale $0<r\les R$. We define $\mc{Q}_r(Q)$ to be a collection of disjoint subcubes of width $ 2^{-j_0} R$ which form a cover of $Q$, where $j_0$ is the unique integer such that $2^{-1-j_0} R< r \les 2^{-j_0} R$. Thus all cubes in $Q_{r}(Q)$ have side lengths $\approx r$, and moreover, if $r \les r'\les R$ and $q \in \mc{Q}_r(Q)$, $q' \in \mc{Q}_{r'}(Q)$ with $q\cap q' \not = \varnothing$, then $q \in \mc{Q}_{r}(q')$. To estimate various error terms which arise, we will need to create some separation between cubes. To this end, following Tao, we introduce the following construction. Given $0<\epsilon \ll 1$ and a subscale $0<r\les R$, we let
        $$ I^{\epsilon, r}(Q) = \bigcup_{q \in \mc{Q}_r(Q)} (1-\epsilon) q. $$
Note that we have the crucial property, that if $q \in \mc{Q}_r(Q)$, and $(t,x) \not \in I^{\epsilon, r}(Q) \cap q$, then $ \dist\big( (t,x), q\big) \g \epsilon r$. \\

We require smoothed out bump functions adapted to cubes $q\in \mc{Q}_r(Q) $. To this end, given $q \in \mc{Q}_{r}(Q)$ we let $\chi_q \in C^\infty(\RR^{1+n})$ be a positive function such that $\chi_q \gtrsim 1$ on $q$, $\supp \widetilde{\chi}_q \subset \{ |(\tau, \xi)| \les \frac{1}{r}\}$, and we have the decay bound, for every $N \in \NN$,
	$$ \chi_q(t,x) \lesa_N \Big(1 + \frac{  \dist((t,x), q)}{r} \Big)^{-N}. $$
Given a point $\gamma_j=(x_0, \xi_0) \in \Gamma_j$, and a cube $q \in \mc{Q}_{r_j}(Q)$, we define the weights
	$$ w_{\gamma_j, q} = \Big( 1 + \frac{|x_q - x_0 + t_q \nabla \Phi_j(\xi_0)|}{r_j}\Big)^{5n}$$
where $(t_q, x_q)$ denotes the centre of the cube $q$. Thus the weights $w_{\gamma_j, q}$ are essentially one when $T_{\gamma_j} \cap q \not = \varnothing$, and are very large when $q$ and the tube $T_{\gamma_j}$ are far apart.\\

For a subset $\Omega \subset \RR^{1+n}$ we let $\ind_{\Omega}$ denote the indicator function of $\Omega$. Let $\mc{E}$ be a finite collection of subsets of $\RR^{1+n}$, and suppose we have a collection of ($\ell^2_c$-valued) functions $(u^{(E)})_{E\in \mc{E}}$. We then define the associated \emph{quilt}
            $$ [u^{(\cdot)}](t,x) = \sum_{E \in \mc{E}} \ind_E(t,x) |u^{(E)}(t,x)|.$$
This notation was introduced by Tao \cite{Tao2001b}, and plays a key technical role in localising the product $uv$ into smaller scales.


\section{Applications}\label{sec-applications}


\subsection{Alternative conditions on phases} The key assumption $\mb{(A1)}$ is a global condition, in the sense that it depends on the behaviour of the phase $\Phi_j$ at two points $\xi, \xi' \in \Sigma_j(\mathfrak{h})$. In certain cases, it can instead be convenient to have a \emph{local} version of $\mb{(A1)}$ which only depends on the behaviour of $\Phi_j$ in a neighbourhood of $\xi \in \Sigma_j(\mathfrak{h})$. Clearly a suitable candidate for such a local condition can be found by simply taking $\xi' \to \xi$ in $\mb{(A1)}$. More precisely, we can state a local version of $\mb{(A1)}$ as:\\

\begin{enumerate}
	\item[\textbf{(A1')}] for every $\{j, k\} = \{1,2\}$, $\xi \in \Lambda_j$, $\eta \in \Lambda_k$, and $v\in \RR^n$ with $v \cdot (\nabla \Phi_j(\xi) - \nabla \Phi_k(\eta) ) = 0$, we have
	$$ \big| \big[ \nabla^2 \Phi_j(\xi) v \big] \wedge \big( \nabla \Phi_j(\xi) - \nabla \Phi_k(\eta) \big)\big| \g \mb{C}'_0 \mc{H}_j \mc{V}_{max} |v|. $$\\
\end{enumerate}
Note that if $\mathfrak{h} = ( \Phi_j(\xi) + \Phi_k(\eta), \xi + \eta)$, then $\nabla \Phi_j(\xi) - \nabla \Phi_k(\eta)$ is the normal vector to the tangent space $T_\xi \Sigma_j(\mathfrak{h})$. Thus the condition on $v$ in (\textbf{A1'}) can be written as $v \in T_\xi \Sigma_j(\mathfrak{h})$. In particular, letting $\xi' \to \xi$ in $\Sigma_j(\mathfrak{h})$, we see that $\mb{(A1)}$ implies $\mb{(A1')}$. To prove the converse implication requires making an additional global assumption on the behavior of the phases $\Phi_j$ on the sets $\Lambda_j$. One possibility is the following.

\begin{lemma}[Local implies global]\label{lem:local version of A1}
Let $\mb{C}_0'>0$. Assume that for all $\mathfrak{h} \in \RR^{1+n}$ we have
    \begin{equation}\label{eqn:lem loc ver A1:assump on set}
        \conv[\Sigma_1(\mathfrak{h})] \subset \Lambda_1 \cap ( h - \Lambda_2),
    \end{equation}
and the small variation bounds
    \begin{equation}\label{eqn:lem loc ver A1:small normal variation}
        \sup_{\xi, \xi' \in \Lambda_1} \big| \nabla \Phi_1(\xi) - \nabla \Phi_1(\xi') \big| + \sup_{\eta, \eta' \in \Lambda_2} \big| \nabla \Phi_2(\eta) - \nabla \Phi_2(\eta') \big| \les \frac{1}{4} \mb{C}_0' \mc{V}_{max}
    \end{equation}
and, for all $\xi, \xi' \in \Sigma_j(\mathfrak{h})$,
    \begin{equation}\label{eqn:lem loc ver A1:small Hessian variation}
       \big| \nabla \Phi_j(\xi) - \nabla \Phi_j(\xi') - \nabla^2 \Phi_j(\xi) (\xi-\xi') \big|\les \frac{1}{4} \mb{C}_0' \mc{H}_j |\xi - \xi'|.
    \end{equation}
If $\mb{(A1')}$ holds, then $\mb{(A1)}$ also holds with constant $\mb{C}_0 = \frac{1}{4} \mb{C}_0'$.
\end{lemma}
\begin{proof} Let $\mathfrak{h}=(a,h) \in \RR^{1+n}$, $\xi, \xi' \in \Sigma_j(\mathfrak{h})$, and $\eta \in \Lambda_2$. For ease of notation, we let $N=\nabla \Phi_j(\xi) - \nabla \Phi_k(\eta)$ and $v = {(\xi - \xi') - [(\xi - \xi') \cdot \frac{N}{|N|}] \frac{N}{|N|}}$. Suppose for the moment that we have
    \begin{equation}\label{eq:lem local ver of A1:dot product small} |(\xi-\xi') \cdot N|\les \frac{1}{4} \mb{C}_0' \mc{V}_{max} |\xi - \xi'|.\end{equation}
As $v \cdot N=0$, an application of $\mb{(A1')}$ and \eref{eqn:lem loc ver A1:small Hessian variation} then gives
    \begin{align*}
         \big|\big(\nabla &\Phi_j(\xi) - \nabla \Phi_j(\xi')\big)
            \wedge N \big| \\
                &\g \big| \big[ \nabla^2 \Phi_j(\xi) v \big] \wedge N \big| - \mc{V}_{max}\big| \nabla \Phi_j(\xi) - \nabla \Phi_j(\xi') - \nabla^2 \Phi(\xi) (\xi-\xi')| - \mc{H}_j \big|(\xi-\xi') \cdot N \big|  \\
                &\g \mb{C}'_0 \mc{H}_j \mc{V}_{max} \Big(|v| -  \frac{1}{2} |\xi-\xi'|\Big) \g \frac{1}{4} \mb{C}_0' \mc{H}_j \mc{V}_{max} |\xi- \xi'|
    \end{align*}
and hence $\mb{(A1)}$ follows. Thus it remains to verify the bound \eref{eq:lem local ver of A1:dot product small}. We start by observing that the definition of the surface $\Sigma_j(\mathfrak{h})$ implies that $\Phi_j(\xi) - \Phi_j(\xi') = \Phi_k(h-\xi) - \Phi_k(h-\xi')$. On the other hand, as $\Sigma_j(\mathfrak{h}) = h - \Sigma_k(\mathfrak{h})$, \eref{eqn:lem loc ver A1:assump on set} gives $\xi' + t(\xi-\xi') \in \Lambda_j \cap ( h - \Lambda_k)$ for $0\les t \les 1$. Consequently we have the identity
    $$ \int_0^1 \nabla \Phi_j\big(\xi' + t(\xi-\xi')\big) \cdot (\xi - \xi') dt = \int_0^1 \nabla \Phi_k\big( h-\xi' + t(\xi'-\xi)\big) \cdot (\xi'-\xi) dt $$
and hence, by an application of \eref{eqn:lem loc ver A1:small normal variation}, we have
    \begin{align*}
       \Big|N \cdot \frac{\xi - \xi'}{|\xi - \xi'|}\Big|
              &\les \int_0^1 \big| \nabla \Phi_j(\xi) - \nabla \Phi_j\big(\xi' + t(\xi - \xi') \big) \big| dt + \int_0^1 \big| \nabla \Phi_k(\eta) - \nabla \Phi_k\big( h - \xi' + t (\xi' - \xi)\big)\big| dt \\
              &\les \big(\mc{H}_1 + \mc{H}_2\big) \mb{d}_0 \les \frac{1}{4}  \mb{C}'_0 \mc{V}_{max}.
    \end{align*}
\end{proof}

Strictly speaking, the above argument shows that to move from the local condition $\mb{(A1')}$ to the global condition $\mb{(A1)}$, it is enough to verify the bounds \eref{eqn:lem loc ver A1:small Hessian variation} and \eref{eq:lem local ver of A1:dot product small}. However, in most applications of Lemma \ref{lem:local version of A1}, the conditions \eref{eqn:lem loc ver A1:assump on set}, \eref{eqn:lem loc ver A1:small normal variation}, and \eref{eqn:lem loc ver A1:small Hessian variation} can always be imposed by perhaps shrinking the sets $\Lambda_j$ slightly. Thus it generally suffices to simply check $\mb{(A1')}$.\\

In some applications of Theorem \ref{thm-main}, it can be convenient to replace the assumption (\textbf{A1}) with a stronger sufficient condition.

\begin{lemma}[{\cite[Lemma 2.1]{Candy2016}}]\label{lem-simplified conditions on phases}
Assume that \eref{eqn:lem loc ver A1:assump on set} and \eref{eqn:lem loc ver A1:small normal variation} hold. If for all $\xi, \xi' \in \Sigma_j(\mathfrak{h})$ we have the curvature type bound
        \begin{equation}\label{eqn:lem simp cond on phase:curv cond}
            \big| \big(\nabla\Phi_j(\xi) - \nabla \Phi_j(\xi')\big) \cdot (\xi - \xi') \big| \g \mb{C}_0^{(1)} \mc{H}_j |\xi - \xi'|^2,
        \end{equation}
and for $\xi \in \Lambda_1$ and $\eta \in \Lambda_2$ the transversality bound
        \begin{equation}\label{eqn:lem simp cond on phase:trans}
             \big| \nabla \Phi_1(\xi) - \nabla \Phi_2(\eta) \big| \g \mb{C}_0^{(2)} \mc{V}_{max}
        \end{equation}
then $\mb{(A1)}$ holds with constant $\mathbf{C}_0 = \frac{1}{2}\mb{C}_0^{(1)}\mb{C}_0^{(2)}$.
\end{lemma}
\begin{proof} We adapt the argument in \cite{Candy2016} to the current setup. Let $\eta \in \Lambda_k$, $\xi, \xi' \in \Sigma_j(\mathfrak{h})$ and take $\omega = \frac{\xi - \xi'}{|\xi - \xi'|}$. The convexity condition \eref{eqn:lem loc ver A1:assump on set} implies that $|\nabla \Phi_j(\xi) - \nabla \Phi_j(\xi')|\les \mc{H}_j |\xi - \xi'|$. Hence  the assumptions \eref{eqn:lem simp cond on phase:curv cond} and \eref{eqn:lem simp cond on phase:trans} give
\begin{align*} \big|\big(\nabla \Phi_j(\xi) - &\nabla \Phi_j(\xi')\big)
            \wedge \big( \nabla \Phi_j(\xi) - \nabla
            \Phi_k(\eta)\big)|\\
             &\g |\nabla \Phi_j(\xi) - \nabla \Phi_k(\eta)|   \left|
            \big(\nabla
            \Phi_j(\xi) - \nabla \Phi_j(\xi')\big) \cdot \omega \right| - |\nabla \Phi_j(\xi) - \nabla \Phi_j(\xi')|   \left|
            \big(\nabla
            \Phi_j(\xi) - \nabla \Phi_k(\eta)\big) \cdot \omega \right|\\
              &\g \mc{H}_j \mc{V}_{max} |\xi - \xi'| \left( \mb{C}_0^{(1)}\mb{C}_0^{(2)} -  \frac{ |
            \big(\nabla
            \Phi_j(\xi) - \nabla \Phi_k(\eta)\big) \cdot \omega |}{ \mc{V}_{max} } \right)
\end{align*}
where we used the elementary bound $|x\wedge y| \g |y| |x\cdot \omega| - |x| |y \cdot \omega|$. Consequently $\mb{(A1)}$ holds provided that for every $\xi, \xi' \in \Sigma_j(\mathfrak{h})$, $\eta \in \Lambda_k$ we have
        \begin{equation}\notag \left|
            \big(\nabla
            \Phi_j(\xi) - \nabla \Phi_k(\eta)\big) \cdot \frac{\xi -
              \xi'}{|\xi - \xi'|} \right| \les \frac{ \mb{C}_0^{(1)}\mb{C}_0^{(2)}}{2 } \mc{V}_{max} . \end{equation}
But this follows from the argument used to prove \eref{eq:lem local ver of A1:dot product small} together with the small variation condition \eref{eqn:lem loc ver A1:small normal variation}.
\end{proof}

\begin{remark} To better understand the connection between the conditions in Lemma \ref{lem-simplified conditions on phases}, and the assumptions (\textbf{A1}) and (\textbf{A1'}), we observe that the later conditions are roughly equivalent to
         $$ \big| \big(\nabla^2 \Phi_j(\xi) v \big) \wedge N \big| \gtrsim \mc{H}_j \mc{V}_{max} |v| $$
for all $v \in T_{\xi} \Sigma_j(\mathfrak{h})$, where $N$ is the unit normal to $\Sigma_j(\mathfrak{h})$, and $T_\xi \Sigma_j(\mathfrak{h})$ is the tangent plane at $\xi \in \Sigma_j(\mathfrak{h})$. On the other hand, letting $\xi' \rightarrow \xi$, we see that \eref{eqn:lem simp cond on phase:curv cond} is essentially equivalent to the lower bound
	\begin{equation}\label{eqn:alt curve assump loc} \big| \big( \nabla^2\Phi_j(\xi) v \big) \cdot v \big|  \gtrsim  \mc{H}_j |v|^2 \end{equation}
for all $v \in T_\xi \Sigma_j(\mathfrak{h})$ (more precisely, it is easy to check that under an assumption like \eref{eqn:lem loc ver A1:small Hessian variation}, the global condition \eref{eqn:lem simp cond on phase:curv cond} is equivalent to \eref{eqn:alt curve assump loc}). Consequently, (\textbf{A1}) and (\textbf{A1'}) state that the Hessian $\nabla^2\Phi_j(\xi)$  can not rotate tangent vectors in $T_{\xi} \Sigma_j(\mathfrak{h})$ to normal vectors, while \eref{eqn:lem simp cond on phase:curv cond} imposes the much stronger condition, that $\nabla^2\Phi_j(\xi) v$ remains roughly parallel to $v \in T_\xi \Sigma_j(\mathfrak{h})$. A similar remark, although phrased slightly differently, was made in \cite{Bejenaru2016b}.
\end{remark}

 In general, even with the use of the previous simplifications, it can be somewhat involved to apply Theorem \ref{thm-main} as the dependence on $\Phi_j$ is only sharp when applied phase which behave essentially \emph{uniformly} on $\Lambda_j$. However a rough strategy is as follows. Suppose we are given phases $\Phi^*_j$ on some (large) subset of $\RR^n$. The first step is to restrict the domain until the gradients $\nabla \Phi^*_j$ have ``small'' variation. The second step is to rescale the domain, and hence replace $\Phi_j^*$ with a rescaled version $\Phi_j$,  to ensure that the (nonzero) components of the Hessian $\nabla^2 \Phi_j$ are all of a similar size. The final step is then to restrict the rescaled domain further, to ensure that the transversality and curvature conditions hold. Clearly each of these steps depends heavily on the precise phases under consideration, and thus  some amount of trial and error is needed to find appropriate domains.

\subsection{Elliptic Phases: Proof of Theorem \ref{thm-elliptic bilinear}} A short computation shows that, after potentially partitioning the sets $\Lambda_j$ into smaller sets, the conditions \eref{eqn:lem loc ver A1:small normal variation} and \eref{eqn:lem loc ver A1:small normal variation} in Lemma \ref{lem:local version of A1} hold. Hence, using the ellipticity assumption, since for all $\xi \in \Lambda_j$, $\eta \in \Lambda_k$, and $v \cdot ( \nabla \Phi_j(\xi) - \nabla \Phi_k(\eta) ) =0$ we have
    $$ \big| \big( \nabla^2 \Phi(\xi) v\big) \wedge \big( \nabla \Phi_j(\xi) - \nabla \Phi_k(\eta) \big) \big| \gtrsim \big| \big( \nabla^2 \Phi_j(\xi) v \big) \cdot v \big| \mc{V}_{max} \gtrsim \mc{H}_j \mc{V}_{max} |v|^2 $$
  we see that $\mb{(A1)}$ holds and thus result follows by an application of Theorem \ref{thm-main small scale}.

\subsection{The Wave/Klein-Gordon case: Proof of Theorem \ref{thm-bilinear small scale KG}}
 We only consider the case $\alpha \ll 1$, the case $\alpha \approx 1$ is easier and follows directly from Theorem \ref{thm-main} (via Lemma \ref{lem-simplified conditions on phases}). Applying a rotation, we may assume that $\xi_0=(a, 0, \dots, 0)$, $\eta_0 = (\sqrt{1-\delta^2} b, \delta b, 0, \dots, 0)$ with $0<\delta < 1$, $\lr{a}_{m_1} \approx \lambda$, $\lr{b}_{m_2} \approx \mu$, and
    \begin{equation}\label{eqn-cor small scale KG-separation assump} \frac{|m_2a -  m_1 b|}{\lambda \mu} + \Big( \frac{ a b}{\lambda \mu}\Big)^\frac{1}{2} \delta  \approx \alpha .\end{equation}
After rescaling, and perhaps decomposing $u$ and $v$ into slightly smaller sets if necessary, it is enough to prove that if $\supp \widehat{f} \subset \Lambda_1$, $\supp \widehat{g} \subset \Lambda_2$ we have
  \begin{equation}\label{eqn-cor small scale KG-reduced estimate}\| e^{ it \Phi_1} f e^{ i t \Phi_2} g \|_{L^q_t L^r_x(\RR^{1+n})} \lesa \alpha^{-\frac{2}{q}}  \mu^{ n-\frac{n}{r}-\frac{1}{q}} \Big( \frac{\lambda}{\mu}\Big)^{\frac{1}{q} - \frac{1}{2}}\|f \|_{L^2_x} \| g \|_{L^2_x}\end{equation}
where we define the phases
    $$ \Phi_1(\xi) = ( m_1^2 + \beta^2 \xi_1^2 + \alpha^2 |\xi'|^2 )^\frac{1}{2}, \qquad \Phi_2(\xi) = \pm ( m_2^2 + \beta^2 \xi_1^2 + \alpha^2 |\xi'|^2 )^\frac{1}{2}$$
with $\xi = (\xi_1, \xi') \in \RR \times \RR^{n-1}$, $\xi'=( \xi_2, \xi'') \in \RR \times \RR^{n-2}$, and take $\xi_0=(a, 0, \dots, 0)$, $\eta_0 = (\sqrt{1-\delta^2} b, \delta b, 0, \dots, 0)$ and the sets $\Lambda_j$ to be
    $$ \Lambda_1 = \big\{ \big|\beta \xi_1 - a \big| \ll \beta \lambda, \,\, |\xi'| \ll \lambda\big\}, \qquad \Lambda_2 = \big\{ \big|\beta \xi_1 \mp \sqrt{1-\delta^2}b \big| \ll \beta \mu, \,\, |\alpha \xi_2 \mp \delta b| \ll \alpha \mu, \,\, |\xi''| \ll \mu \}. $$
Note that \eref{eqn-cor small scale KG-separation assump} implies that if $(\frac{ a b}{\lambda \mu})^\frac{1}{2} \delta \ll \alpha$ then the sets $\Lambda_j$ must have some radial separation, while if $(\frac{ a b}{\lambda \mu})^\frac{1}{2} \delta \gtrsim \alpha$ the sets are angularly separated. To simplify the computations to follow, we note that the transversality condition (\ref{eqn-cor small scale KG-separation assump}) implies that $\frac{b}{\mu} \delta \lesa \alpha$, this follows by observing that in the fully elliptic case, when $\lambda \approx  m_1$ and $ \mu \approx  m_2$, (\ref{eqn-cor small scale KG-separation assump}) becomes $ \Big| \frac{\eta_0}{m_1} - \frac{\xi_0}{m_2} \Big| \approx \alpha$.
 Similarly, by the definition of $\beta$, if either $\lambda \approx m_1$ or $\mu \approx m_2$, then we have $\beta \approx \alpha$.

To obtain the bound (\ref{eqn-cor small scale KG-reduced estimate}), we apply Theorem \ref{thm-main}. Thus we have to check the conditions (\textbf{A1}) and (\textbf{A2}). There a number of ways to do this, here we argue via Lemma \ref{lem-simplified conditions on phases}. To this end, we start by noting that
	$$ \nabla \Phi_j = \frac{ (\beta^2 \xi_1, \alpha^2 \xi')}{(m^2_j + \beta^2 \xi_1^2 + \alpha^2 |\xi'|^2)^\frac{1}{2}}$$
and hence a computation using the definitions of $\beta$ and $\Lambda_j$, gives the derivative bounds
        \begin{equation}\label{eqn-cor small scale KG-gradient bound}
        \sup_{\xi \in \Lambda_1} \Big|\nabla \Phi_1(\xi) -  \frac{ (\beta a, 0, \dots, 0)}{\lr{a}_{m_1}}\Big| + \sup_{\xi \in \Lambda_2}  \Big|\nabla \Phi_2(\xi) - \frac{(\beta \sqrt{1-\delta^2}b,\alpha \delta b, 0,\dots, 0)}{\lr{ b }_{m_2}}\Big| \approx  \frac{\alpha^2}{10} .
    \end{equation}
Similarly, to estimate $\nabla^2 \Phi_j$, we observe that
	$$ \p_1^2 \Phi_j =   \frac{ \beta^2 }{(m^2_j + \beta^2 \xi_1^2 + \alpha^2 |\xi'|^2)^\frac{1}{2}} - \frac{\beta^4 \xi^2_1}{(m^2_j + \beta^2 \xi_1^2 + \alpha^2 |\xi'|^2)^\frac{3}{2}} = \alpha^2 \frac{\beta^2 ( (\frac{m_j}{\alpha})^2 + |\xi'|^2)}{(m^2_j + \beta^2 \xi_1^2 + \alpha^2 |\xi'|^2)^{\frac{3}{2}}}$$
and  for $k, k'>1$,
	$$ \p_1 \p_k \Phi_j = \alpha^2 \frac{-\beta^2 \xi_1 \xi_k}{( m^2_j + \beta^2 \xi_1 + \alpha^2 |\xi'|)^\frac{3}{2}}, \qquad \p_k \p_{k'} \Phi_j = \alpha^2 \Big( \frac{\delta_{k, k'}}{(m^2_j+\beta^2 \xi_1^2 + \alpha^2 |\xi'|^2)^\frac{1}{2}} -  \frac{ \alpha^2 \xi_k \xi_{k'}}{(m^2_j+\beta^2 \xi_1^2 + \alpha^2 |\xi'|^2)^\frac{3}{2}}\Big).$$
In particular, applying the definition of $\beta$, we have the Hessian bounds $ \mc{H}_1 \approx \frac{\alpha^2}{\lambda}$, $\mc{H}_2 \approx \frac{\alpha^2}{\mu}$. Thus we see that the condition (\textbf{A2}) holds with $\mb{d}_0 = \mu$. It remains to show that (\textbf{A1}) holds, which will follow by checking the conditions in Lemma \ref{lem-simplified conditions on phases}. To this end, the transversality condition (\ref{eqn-cor small scale KG-separation assump}) together with the definition of $\beta$ gives
 \begin{align*}
    \Big| \frac{ (\beta a, 0, \dots, 0)}{\lr{a}_{m_1}} - \frac{(\beta \sqrt{1-\delta^2} b, \alpha \delta b, 0,\dots, 0)}{\lr{b}_{m_2}}\Big|
        &\approx \alpha^2
 \end{align*} The estimate (\ref{eqn-cor small scale KG-gradient bound}) then gives the transversality bound
    $$ |\nabla \Phi_1(\xi) - \nabla \Phi_2(\eta)| \approx \mc{V}_{max} \approx \alpha^2. $$
Since we clearly have the final condition in Lemma \ref{lem-simplified conditions on phases}, it only remains to show that we have the curvature condition
    \begin{equation}\label{eqn-cor small scale KG-curvature bound} \big|\big(\nabla \Phi_1(\xi) - \nabla \Phi_1(\eta) \big) \cdot (\xi - \eta) \big| \gtrsim \frac{\alpha^2}{ \lambda} |\xi - \eta|^2 \end{equation}
for all $\xi, \eta \in \Sigma_1(\mathfrak{h})$, together with a similar bound for $\nabla \Phi_2$. This follows by adapting the argument in \cite[Corollary 6.4]{Candy2016}). We only consider the case $j=1$, as the remaining case is identical. Suppose that $\xi, \eta \in \Sigma_1(\mathfrak{h})$ with $\mathfrak{h} = (a^*, h)$. Let $x= (m_j, \beta \xi_1, \alpha \xi')$, $y = ( m_j, \beta \eta_1, \alpha \eta')$, and $h^*=(m_2-m_1, \beta h_1, \alpha h')$. Then the condition $\xi \in \Sigma_1(\mathfrak{h})$ becomes $|x| \pm |x - h^*| = a^*$. Consequently, for every $\xi, \eta \in \Sigma_1(\mathfrak{h})$ we have
\begin{align}
	 \big| \big(\nabla \Phi_1(\xi) - &\nabla \Phi_1(\eta) \big) \cdot (\xi - \eta) \big|\notag\\
	 			&=  \frac{ |( m_1, \beta \xi_1, \alpha \xi')|+ |(m_1, \beta \eta_1, \alpha \eta')|}{2}  \Bigg|  \frac{( m_1, \beta \xi_1, \alpha \xi')}{|( m_1, \beta \xi_1, \alpha \xi')|} - \frac{( m_1, \beta \eta_1, \alpha \eta') }{|( m_1, \beta \xi_1, \alpha \xi')|} \Bigg|^2 \notag \\
&\approx \lambda\bigg| \frac{x}{|x|} - \frac{y}{|y|} \bigg|^2 \notag \\
&\approx \frac{1}{\lambda} \Big( \lambda^2 \bigg| \frac{x}{|x|} - \frac{y}{|y|} \bigg|^2 + \mu^2 \bigg| \frac{x-h^*}{|x-h^*|} - \frac{y-h^*}{|y-h^*|} \bigg|^2\Big)
\label{eqn-cor small scale bilinear est-curvature bound (i)}
\end{align}
where the last line follows by observing that $\xi, \eta \in \Sigma_1(\mathfrak{h})$ implies that we have the identity
 \begin{align*} |x| |y| \Big| \frac{x}{|x|} - \frac{y}{|y|} \Big|^2 = \big( |x-y|^2 - \big| |x| - |y| \big|^2 \big) &=  \big(\big|(x-h^*) - (y-h^*)\big|^2 - \big| |x-h^*| - |y-h^*| \big|^2 \big) \\
        &= |x-h^*| |y-h^*| \Big| \frac{x-h^*}{|x-h^*|} - \frac{y-h^*}{|y-h^*|} \Big|^2.
     \end{align*}
If either $\lambda \approx m_1$ or $\mu \approx m_2$, then since $z, w \in \RR^n$ with $\lr{z}_m \approx \lr{w}_m$ implies we have
    $$ \Big| \frac{z}{\lr{z}_m} - \frac{x}{\lr{x}_m} \Big| \approx \Big( \frac{m}{\lr{z}_m}\Big)^2 \frac{\big| |z| - |w| \big|}{\lr{z}_m} + \Big( \frac{ |z| |w| - z\cdot w}{\lr{z}_m^2 } \Big)^{\frac{1}{2}}, $$
the required bound (\ref{eqn-cor small scale KG-curvature bound}) follows directly from (\ref{eqn-cor small scale bilinear est-curvature bound (i)}) together with the fact that $\alpha \approx \beta$ in the case $\lambda \approx m_1$ or $\mu \approx m_2$. On the other hand, if $\lambda \gg m_1$ and $\mu \gg m_2$, then we consider separately the cases $|\xi_1 - \eta_1| \gtrsim |\xi' - \eta'|$ and $|\xi' - \eta'| \ll |\xi_1 - \eta_1|$. In the former case, if $\alpha \gg \frac{m_1}{\lambda} + \frac{m_2}{\mu}$, then from the definition of $\beta$ and the condition (\ref{eqn-cor small scale KG-separation assump}) we must have $\beta \approx \delta \approx 1$. Hence we deduce that as $\xi - h, \eta-h \in \Lambda_2$ we have
    $$ |(x-h^*)\wedge(y-h^*)| \g \alpha | \beta (\xi_1 - \eta_1) (\eta' - h') - \beta (\eta_1 - h_1) ( \eta' - \xi')| \approx \alpha \mu |\xi_1 - \eta_1| \approx \alpha \mu |\xi - \eta|$$
which, together with the elementary bound $|\omega - \omega'| \gtrsim |\omega \wedge \omega'|$ for $\omega, \omega'\in \sph^{n-1}$ gives (\ref{eqn-cor small scale KG-curvature bound}). On the other hand, if $\alpha \lesa \frac{m_1}{\lambda} + \frac{m_2}{\mu}$, we have again using the definition of $\beta$
    $$ \frac{|x \wedge y|}{\lambda} + \frac{|(x-h^*)\wedge (y-h^*)|}{\mu} \g \beta \frac{m_1 |\xi_1 - \eta_1|}{\lambda} + \beta \frac{m_2 |\xi_1 - \eta_1|}{\mu}  \approx \alpha |\xi - \eta|.$$
Thus we have (\ref{eqn-cor small scale KG-curvature bound}) when $|\xi_1 - \eta_1| \gtrsim |\xi' - \eta'|$. Finally, if $|\xi'-\eta'| \gg |\xi_1 - \eta_1|$, then we simply observe that
    $$ |x\wedge y| \g \alpha \big| \beta \xi_1 (\xi' - \eta') - \xi' \beta( \xi_1 - \eta_1) \big| \approx \alpha \lambda |\xi' - \eta'| \approx \alpha \lambda  |\xi - \eta|.$$
Therefore (\ref{eqn-cor small scale KG-curvature bound}) follows as required.

\subsection{A Refined Strichartz estimate for the Klein-Gordon equation: Proof of Theorem \ref{thm-refined strichartz}}
\label{subsec:refined strichartz}
The proof closely follows the argument of Ramos in \cite{Ramos2012}, with the key bilinear estimate Theorem \ref{thm-bilinear small scale KG} replacing the bilinear estimate of Tao \cite{Tao2001b} used in \cite{Ramos2012}, thus we shall be somewhat brief. Let $p = \frac{n+1}{n-1}$ and $u = e^{ i t \lr{\nabla} }f$, and decompose
    $$ \| u \|_{L^{2p}_{t,x}}^2 \les \sum_{\ell \in 2^\NN} \Big\| \sum_{\lambda \in 2^\NN} u_\lambda u_{\ell \lambda} \Big\|_{L^p_{t,x}} $$
where $u = \sum_{\lambda \in 2^{\NN}} u_\lambda$ and $\supp \widehat{u}_\lambda \subset \{ \lr{\xi}  \approx \lambda\}$. Since the sum is essentially orthogonal as $\lambda$ varies, by applying the weak orthogonality in $L^p$ (see for instance \cite[Lemma 2.2]{Ramos2012}) we have
    $$ \Big\| \sum_{\lambda \in 2^\NN} u_\lambda u_{\ell \lambda} \Big\|_{L^p_{t,x}} \lesa \Big( \sum_{\lambda\in 2^{\NN} } \| u_\lambda u_{\ell \lambda} \|_{L^p_{t,x}}^{\min\{p, p'\}} \Big)^{\frac{1}{\min\{p, p'\}}}$$
where $p' = \frac{p}{p-1}$ denotes the dual exponent to $p$. Given $A = A(\xi_0) \in \mc{A}_{\lambda, \alpha}$, and $B = B(\eta_0) \in \mc{A}_{\ell \lambda, \alpha}$, we say $A \sim B$ if we have the transversality type condition
        $$ \frac{ \big| |\xi_0| - |\eta_0| \big|}{ \ell \lambda^2} + \Big( \frac{ |\xi_0| |\eta_0| - \xi_0 \cdot \eta_0}{\ell \lambda^2} \Big)^\frac{1}{2} \approx \alpha. $$
Applying a Whitney type decomposition, we can write
    $$ u_\lambda u_{\ell \lambda} = \sum_{\alpha \in 2^{-\NN}} \sum_{A \in \mc{A}_{\lambda, \alpha}} \sum_{\substack{ B \in \mc{A}_{\ell \lambda, \alpha} \\ A \sim B}} u_A u_B $$
where $\supp \widehat{u}_A \subset A$, and $\| \widehat{u}_A \|_{L^r_\xi} \approx \| \widehat{u} \|_{L^r(A)}$, $\| \widehat{u}_A \|_{L^r_\xi} \approx \| \widehat{u} \|_{L^r(A)}$. An application of Theorem \ref{thm-bilinear small scale KG} gives for all $\frac{n+3}{n+1} <q \les 2$
    $$ \| u_A u_B \|_{L^q_{t,x}} \lesa \alpha^{n-1 - \frac{n+1}{q}} \ell^{\frac{1}{q} - \frac{1}{2}}\Big( \frac{\alpha \lambda}{1 + \alpha \lambda}\Big)^{1-\frac{1}{q}} \lambda^{n-\frac{n+1}{q}} \| \widehat{f} \|_{L^2_\xi(A)} \| \widehat{f} \|_{L^2_\xi(B)}. $$
Interpolating with the trivial bound
    $$ \| u_A u_B \|_{L^\infty_{t,x}} \lesa \| \widehat{f} \|_{L^1_\xi(A)} \| \widehat{f} \|_{L^1_{\xi}(B)} $$
we deduce that for all $\max\{ \frac{1}{2}, \frac{2}{n+1} \} \les \frac{1}{r} < \frac{1}{2} + \frac{2}{(n+1)^2}$ we have
    \begin{align}
        \| u_A u_B \|_{L^p_{t,x}} &\lesa \alpha^{n-1  - \frac{2(n-1)}{r} } \ell^{\frac{n-1}{n+1} + \frac{1}{r} - 1} \Big( \frac{\alpha \lambda}{1 + \alpha \lambda}\Big)^{2 - \frac{2}{r} - \frac{n-1}{n+1}} \lambda^{n + 1 - \frac{2n}{r} } \| \widehat{f} \|_{L^r_\xi(A)} \| \widehat{f} \|_{L^r_{\xi}(B)} \notag \\
        &\approx \ell^{\frac{n-1}{n+1} + \frac{n+1}{r} - \frac{n+3}{2}} \Big( \frac{\alpha \lambda}{1 + \alpha \lambda}\Big)^{\frac{2}{n+1}} \lambda^{\frac{1}{2}} |A|^{\frac{1}{2} - \frac{1}{r}} \| \widehat{f}\|_{L^r_\xi(A)} (\ell \lambda)^\frac{1}{2} |B|^{\frac{1}{2} - \frac{1}{r}} \| \widehat{g} \|_{L^r_{\xi}(B)}.   \label{eqn-thm refined strichartz-main bilinear input}
    \end{align}
The bound (\ref{eqn-thm refined strichartz-main bilinear input}) is the key bilinear input, and replaces the corresponding estimate for the wave equation \cite[Corollary 2.1]{Ramos2012} used in the work of Ramos. We now apply the Whitney decomposition to deduce that
        \begin{align}
          \| u_\lambda u_{\ell \lambda} \|_{L^p_{t,x}} &\lesa \sup_{\alpha \approx 1} \sup_{A \in \mc{A}_{\lambda, \alpha}} \sup_{B \in \mc{A}_{\ell \lambda, \alpha}} \| u_A u_B \|_{L^p_{t,x}}  \notag \\
          &\qquad \qquad \qquad + \Big\|  \sum_{ \frac{1}{\lambda} \lesa \alpha \ll 1} \sum_{A \in \mc{A}_{\lambda, \alpha}} \sum_{\substack{ B \in \mc{A}_{\ell \lambda, \alpha} \\ A \sim B}} u_A u_B \Big\|_{L^p_{t,x}} + \Big\|  \sum_{\alpha \ll \frac{1}{\lambda}} \sum_{A \in \mc{A}_{\lambda, \alpha}} \sum_{\substack{ B \in \mc{A}_{\ell \lambda, \alpha} \\ A \sim B}} u_A u_B \Big\|_{L^p_{t,x}}
          \label{eqn-thm refined strichartz-whitney decomp}
        \end{align}
where the first term follows from simply observing that for $\alpha \approx 1$ we have $\# \mc{A}_{\lambda, \alpha} \lesa 1$. For the second and third term in (\ref{eqn-thm refined strichartz-whitney decomp}), we again follow \cite{Ramos2012} and apply the weak orthogonality in $L^p$. More precisely, we claim that for fixed $\ell$ and $\lambda$,  the sets
    $$ \{ (\lr{\xi} + \lr{\eta}, \xi + \eta) \,\, | \xi \in A, \,\, \eta \in B \, \} \subset \RR^{1+n}.$$
 for $A \in \mc{A}_{\lambda, \alpha}$, $B \in \mc{A}_{\ell \lambda, \alpha}$ with $A \sim B$, are essentially disjoint as $\alpha$ and $A$ vary. This follows by observing that if $A=A(\xi_0) \in \mc{A}_{\lambda, \alpha}$ and $B=B(\eta_0) \in \mc{A}_{\ell \lambda, \alpha}$ with $A \sim B$, then we have for every $\xi \in A$ and $\eta \in B$
    $$ |\lr{\xi + \eta}_2 - \lr{\xi} - \lr{\eta} | \approx \lambda \alpha^2, \qquad \big| |\xi + \eta| - |\xi_0 + \eta_0| \big| \lesa \frac{\alpha \lambda}{1 + \alpha \lambda} \ell \lambda, \qquad ( |\xi + \eta| |\xi_0 + \eta_0| - (\xi + \eta) \cdot (\xi_0+\eta_0) \big)^\frac{1}{2} \lesa \alpha \ell \lambda.$$
 In particular, for the second term in (\ref{eqn-thm refined strichartz-whitney decomp}), as this sum only covers the wave like regime where the sets $A$ are simply angular sectors of the annuli, we can bound this contribution by using (\ref{eqn-thm refined strichartz-main bilinear input}) together with the weak orthogonality and delicate summation arguments used in \cite{Ramos2012}. Note that these arguments only use the fact that we are summing up over angular sectors of the annuli, in particular, after using (\ref{eqn-thm refined strichartz-main bilinear input}), the difference between the wave and Klein-Gordon equation no longer plays any role.

 On the other hand, for the third term in (\ref{eqn-thm refined strichartz-whitney decomp}), which comprises the main Klein-Gordon type contribution, we start by observing that due to the restriction $A \sim B$, this contribution is only non-vanishing if $\ell \approx 1$. An application of \cite[Lemma 2.2]{Ramos2012}, together with the finite overlap observation made above and (\ref{eqn-thm refined strichartz-main bilinear input}) then gives
    \begin{align*}
      \bigg( \sum_{\lambda} \Big\|  &\sum_{\alpha \ll \frac{1}{\lambda}} \sum_{A \in \mc{A}_{\lambda, \alpha}} \sum_{\substack{ B \in \mc{A}_{\ell \lambda, \alpha} \\ A \sim B}} u_A u_B \Big\|_{L^p_{t,x}}^{\min\{p, p'\}} \bigg)^{\frac{1}{\min\{p, p'\}}} \\
            &\lesa \bigg( \sum_{\lambda} \sum_{\alpha \ll \frac{1}{\lambda}} \sum_{A \in \mc{A}_{\lambda, \alpha}} \sum_{\substack{ B \in \mc{A}_{\ell \lambda, \alpha} \\ A \sim B}} \|   u_A u_B \|_{L^p_{t,x}}^{\min\{p, p'\}} \bigg)^{\frac{1}{\min\{p, p'\}}} \\
            &\lesa \bigg( \sum_{\lambda} \sum_{\alpha \ll \frac{1}{\lambda}} \sum_{A \in \mc{A}_{\lambda, \alpha}} \sum_{\substack{ B \in \mc{A}_{\ell \lambda, \alpha} \\ A \sim B}} \Big( (\alpha \lambda)^{\frac{2}{n+1}} \lambda^\frac{1}{2} |A|^{\frac{1}{2} - \frac{1}{r}} \| \widehat{f} \|_{L^r_\xi(A)} \lambda^\frac{1}{2} |B|^{\frac{1}{2} - \frac{1}{r}} \| \widehat{f} \|_{L^r_\xi(B)}\Big)^{\min\{p, p'\}} \bigg)^{\frac{1}{\min\{p, p'\}}} \\
            &\lesa \Big( \sup_{\lambda, \alpha} \sup_{A\in \mc{A}_{\lambda, \alpha}}  \Big( \tfrac{\alpha \lambda}{1+\alpha \lambda}\Big)^{\frac{1}{n+1}} \lambda^\frac{1}{2} |A|^{\frac{1}{2}-\frac{1}{r}} \| \widehat{f} \|_{L^r_\xi(A)} \Big)^{\frac{2\min\{p, p'\} -1}{\min\{p, p'\}}} \| f \|_{H^\frac{1}{2}}^\frac{2}{\min\{p, p'\}}.
    \end{align*}

\section{Optimality of Theorem \ref{thm-main}}\label{sec:counter example}

In this section, our goal is to show that the conclusion of Theorem \ref{thm-main} is sharp.   We first observe that, by exploiting the dilation and translation invariance of $L^q_t L^r_x$, is enough to consider the case
   \begin{equation}\label{eqn-assum on phase after rescaling II}
 1\les \|\nabla \Phi_1 \|_{L^\infty} + \|\nabla \Phi_2 \|_{L^\infty} \les 3, \qquad \mc{V}_{max}=1, \qquad \mc{H}_2 \les \mc{H}_1 = 1.
   \end{equation}
More precisely, since
	$$ \mc{V}_{max} \les \inf_{\xi_0 \in \RR^n} \big( \|\nabla \Phi_1 - \xi_0\|_{L^\infty(\Lambda_1)} + \|\nabla \Phi_2 - \xi_0 \|_{L^\infty(\Lambda_2)}\big) \les 3 \mc{V}_{max}$$
after a linear translation of the phases (this corresponds to the spatial translation $x\mapsto x+t\xi_0$) we may assume that
   $$ \mc{V}_{max} \les \|\nabla \Phi_1\|_{L^\infty(\Lambda_1)}
				+ \|\nabla \Phi_2\|_{L^\infty(\Lambda_2)} \les 3 \mc{V}_{max}. $$
To obtain the claimed bounds in (\ref{eqn-assum on phase after rescaling II}), we now apply the rescaling
	$$ \Phi_j(\xi) \mapsto \frac{\mc{H}_1}{\mc{V}_{max}^2} \Phi_j\bigg( \frac{\mc{V}_{max}}{\mc{H}_1} \xi\bigg), \qquad \mb{d}_0 \mapsto \frac{\mc{H}_1}{\mc{V}_{max}} \mb{d}_0.$$
It is important to note that this rescaling and translation leaves the assumptions (\textbf{A1}) and (\textbf{A2}), and the bilinear estimate in Theorem \ref{thm-main}, unchanged.

Fix $\xi_0\in \Lambda_1$, $\eta_0 \in \Lambda_2$, and $0<\epsilon \ll 1$. It is enough to construct $\Phi_1$-wave $u$, and a $\Phi_2$-wave $v$, such that $\supp \widehat{u} \subset \{ |\xi - \xi_0| \les \epsilon\}$, $\supp \widehat{v} \subset \{ |\xi - \eta_0| \les \epsilon\}$, we have the energy bounds
    $$\| u \|_{L^\infty_t L^2_x} \lesa \mc{H}_2^{-\frac{1}{2}} \epsilon^{-\frac{n+1}{2}}, \qquad \| v \|_{L^\infty_t L^2_x} \lesa \epsilon^{-\frac{n+1}{2}},$$
and both $u$ and $v$ are concentrated on the $\epsilon^{-1}$ thickened space-time rectangle
    $$ \Omega = \big\{ s\big(1, - \nabla\Phi_1(\xi_0) \big) + s' \big( 1, -\nabla \Phi_2(\eta_0) \big) \,\,\big| \,\, 0\les s \les \epsilon^{-2}, \,\, 0\les s' \les (\epsilon^2 \mc{H}_2)^{-1} \,\big\} + \epsilon^{-1} $$
in the sense that $\| u v \|_{L^q_t L^r_x(\Omega)} \approx \| \ind_\Omega \|_{L^q_t L^r_x}$. Assuming the existence of $u$ and $v$ for the moment, we easily deduce that the estimate
        $$ \| u v \|_{L^q_t L^r_x} \les \mb{C} \| u \|_{L^\infty_t L^2_x} \| v \|_{L^\infty_t L^2_x}$$
together with the transversality assumption (\ref{eqn-main trans assump}) implies that
        $$  \mb{C} \gtrsim \epsilon^{n+1 - \frac{n+1}{r} - \frac{2}{q}} \mc{H}_2^{\frac{1}{2} - \frac{1}{q}}. $$
Thus we see that $\mb{d}_0 \gtrsim \epsilon$, we have the correct dependence on $\mc{H}_j$,  and letting $\epsilon \rightarrow 0$, we must have $\frac{1}{q} + \frac{n+1}{2r} \les \frac{n+1}{2}$.

It remains to construct $u$ and $v$ satisfying the required properties. One approach is to adapt the example used by Tao in \cite{Tao2001b}. Let $\rho$ be a Schwartz function such that $\rho(x) \approx  1$ for $|x| \les 1$, and $\supp \widehat{\rho} \subset \{ |\xi| \les 1\}$, and define
    $$\rho_{1,k}(x) = e^{i (x+N k \nabla \Phi_2(\eta_0))\cdot \xi_0} \rho\big( \epsilon x + \epsilon  N k \nabla \Phi_2(\eta_0)\big), \qquad \rho_{2,k} = e^{i (x+N k \nabla \Phi_1(\xi_0))\cdot \eta_0} \rho\big( \epsilon x  + \epsilon N k \nabla \Phi_1(\xi_0)\big) $$
where $N$ is some fixed quantity independent of $\epsilon$, which is needed to create some separation between the resulting wave packets. The waves $u$ and $v$ are then defined as
    $$ u(t,x) = \sum_{0\les k \lesa (\epsilon^2 \mc{H}_2)^{-1}} e^{ i ( t - Nk) \Phi_1(-i\nabla)} \rho_{1, k}(x), \qquad v(t,x) = \sum_{0\les k \lesa \epsilon^{-2}} e^{ i t \Phi_2(-i \nabla)} \rho_{2, k}(x).  $$
It is clear that $u$ and $v$ satisfy the required Fourier support properties by the definition of $\rho$. Roughly speaking, $v$ is defined as the sum of wave packets concentrated on (space-time) tubes of size  $(\epsilon^2 \mc{H}_2)^{-1} \times \epsilon^{-n}$ which are oriented in the direction $(1, -\nabla \Phi_2(\eta_0))$ and centered along the points $(0, -k N \nabla \Phi_1(\xi_0))$, and essentially cover the set $\Omega$. Similarly,  $u$ is the sum of wave packets concentrated on tubes of size $\epsilon^{-2} \times \epsilon^{-n}$ oriented in the direction $(1, - \nabla \Phi_1(\xi_0))$ and centered along the (space-time) points $(kN, -kN\nabla \Phi_2(\eta_0))$, which again essentially covers the set $\Omega$. In particular, we expect that the product $|uv|$ is essentially concentrated on $\Omega$. These heuristics can be made rigorous by a somewhat standard integration by parts argument (and potentially choosing $N$ large enough to absorb any absolute constants which arise), see for instance Section \ref{sec-wave packets} for similar arguments. On the other hand, the claimed $L^2_x$ bounds follow by noting that as the wave packets are essentially disjoint for fixed times, we have
        $$ \| v \|_{L^\infty_t L^2_x} \approx \Big( \sum_{0\les k \lesa \epsilon^{-2}} \| \rho_{2, k} \|_{L^2_x}^2 \Big)^\frac{1}{2} \approx \epsilon^{-\frac{n+1}{2}}. $$
The argument to bound $\|u \|_{L^\infty_t L^2_x}$ is slightly more involved as the wave packets are not centered on a fixed time slice $t=0$, but instead centred along a tube oriented in the direction $(1, -\nabla \Phi_2(\eta_0))$. In general, wave packets centered along space-time lines need not be orthogonal, however, since the wave packets making up $u$ are \emph{transverse} to the direction $(1, -\nabla\Phi_2(\eta_0))$, we retain orthogonality of the wave packets. This can be seen by an application of Plancheral which gives
    $$ \| u(t) \|_{L^2_x}^2 = \epsilon^{-n} \sum_{0\les k, k' \lesa (\epsilon^2 \mc{H}_2)^{-1}} \int_{\RR^n} e^{ i N(k-k') [ \xi\cdot\nabla \Phi_2(\eta_0) - \epsilon^{-1} \Phi_1(\epsilon \xi-\xi_0)]} |\rho(\xi)|^2 d\xi \lesa \epsilon^{-n} \sum_{0\les k \lesa (\epsilon^2 \mc{H}_2)^{-1}} 1 \approx \mc{H}_2^{-1} \epsilon^{-n-1} $$
where we used integration by parts together with (\ref{eqn-main trans assump}) to control the off diagonal terms. It is worth observing that
the above argument already foreshadows a number of the ingredients in the proof of Theorem \ref{thm-main}, namely the orthogonality properties of wave packets together with energy type estimates across transverse surfaces.


\section{Induction on Scales}\label{sec-ind on scales}


In this section we reduce the proof of Theorem \ref{thm-main small scale}  to proving local estimates on cubes $Q \in \RR^{1+n}$. This follows the argument of Tao \cite{Tao2001b} which was adapted by Bejenaru \cite{Bejenaru2016b} and Lee-Vargas \cite{Lee2008}, but we refine the induction somewhat to remove the epsilon derivative loss that appears in the previous bounds of Tao and Lee-Vargas. \\

Fix constants $\mb{d}_0, \mb{C}_0>0$, and open sets $\Lambda_j \subset \RR^n$. Let $\Phi_1$ and $\Phi_2$ be phases satisfying (\textbf{A1}) and (\textbf{A2}). After exploiting the dilation and translation invariance of $L^q_t L^r_x$ as in Section \ref{sec:counter example}, it is enough to consider the case
   \begin{equation}\label{eqn-assum on phase after rescaling}
 1\les \|\nabla \Phi_1 \|_{L^\infty} + \|\nabla \Phi_2 \|_{L^\infty} \les 3, \qquad \mc{V}_{max}=1, \qquad \mc{H}_2 \les \mc{H}_1 = 1.
   \end{equation}
Take sets $\Lambda_j^* + \mb{d}_0 \subset \Lambda_j$ such that $\mb{d}[\Lambda_1^* + \mb{d}_0, \Lambda_2^*+\mb{d}_0] \les \frac{\mb{d}_0}{\mb{C}_0}$, and fix $R_0 \gg \frac{1}{\mb{d}_0^2 \mc{H}_2}$ with
    $$ \mb{d}\big[\Lambda_1^* + \mb{d}_0, \Lambda_2^* + \mb{d}_0\big] \lesa \big( \mc{H}_2 R_0\big)^{-\frac{1}{2}} \les \mb{d}\big[\Lambda_1^*+\mb{d}_0, \Lambda_2^* + \mb{d}_0 \big].$$
The constant $R_0$ will denote the smallest scale of cubes we consider, while the sets $\Lambda_j^*$ will contain the supports of $\widehat{u}$ and $\widehat{v}$. The proof of Theorem \ref{thm-main small scale} is based on an induction on scales argument. This requires the following definition.

\begin{definition}\label{defn of A(R)}
For any $R \g R_0$ and $1\les q, r \les 2$, we define $A_{q,r}(R)$ to the best constant for which the inequality
	$$ \| uv \|_{L^q_t L^r_x(Q)} \les A_{q,r}(R) \|u \|_{L^\infty_t L^2_x} \| v \|_{L^\infty_t L^2_x}$$
holds for all cubes $Q \subset \RR^{1+n}$ of radius $R$, and all  $\Phi_1$-waves $u$ and $\Phi_2$-waves $v$ satisfying the support assumption
	$$ \supp \widehat{u}  \subset \Lambda_1^* + 4(\mc{H}_2 R)^{-\frac{1}{2}}, \qquad \supp \widehat{v}  \subset \Lambda_2^* + 4 (\mc{H}_2 R)^{-\frac{1}{2}}.$$\\
\end{definition}

It is easy to check that $A_{q,r}(R)$ is always finite, and, since $R \g R_0 \gg (\mb{d}_0^2 \mc{H}_2)^{-1}$, that the required support conditions on $ \widehat{v}$ and $ \widehat{v}$ imply that we always have
	$$ \supp \widehat{u}   +  \frac{\mb{d}_0}{2} \subset \Lambda_1, \qquad \supp \widehat{v} + \frac{\mb{d}_0}{2} \subset \Lambda_2.$$
 It is also worth noting that the support condition becomes stricter as $R$ becomes large, in other words, for large $R$, the Fourier supports must be smaller. Roughly this can be explained as follows. Our goal will be to obtain a bound for $A_{q,r}(2R)$, in terms of $A_{q,r}(R)$. At each scale $2R$, we will need to decompose $u$ and $v$ into wave packets, this will enlarge the Fourier support slightly, and thus to apply the bound at scale $R$ to the wave packets of $u$ and $v$, we will need $A_{q,r}(R)$ to apply to functions with slightly larger support. The support conditions we use play the same role as the margin type conditions used in \cite{Tao2001b}.\\

The proof of Theorem \ref{thm-main small scale} will rely on two key propositions. The first gives a good bound for $A_{q,r}(R)$ for $R$ close to $R_0$, and is needed to begin the induction argument.

\begin{proposition}\label{prop-initial induction bound}
Let $\mb{C}_0>0$,  $1\les q, r \les 2$, and $\frac{1}{q} +\frac{n+1}{2r} < \frac{n+1}{2}$. Assume we have $\mb{d}_0 >0$, open sets $\Lambda_j \subset \RR^n$, and phases $\Phi_j$ satisfying Assumption \ref{assump-main} with the normalisation conditions (\ref{eqn-assum on phase after rescaling}). Then for any sets $\Lambda_j^* + \mb{d}_0 \subset \Lambda_j$ with $\mb{d}[\Lambda_1^* + \mb{d}_0, \Lambda_2^*+\mb{d}_0] \les \frac{\mb{d}_0}{\mb{C}_0}$, and every $R \g R_0$ we have
		$$ A_{q,r}(R) \lesa \mb{d}_0^{n+1 - \frac{n+1}{r} - \frac{2}{q}} \mc{H}_2^{\frac{1}{2}-\frac{1}{q}} \Big( \frac{R}{R_0} \Big)^{\frac{1}{q}}.$$
\end{proposition}

The second implies that we can control $A_{q,r}(2R)$ in terms of $A_{q,r}(R)$.

\begin{proposition}\label{prop-induction step}
Let $\mb{C}_0>0$,  $1\les q, r \les 2$, and $\frac{1}{q} +\frac{n+1}{2r} < \frac{n+1}{2}$. There exists a constant $C>0$, such that for any $ \mb{d}_0>0$, any open sets $\Lambda_j \subset \RR^n$, any phases $\Phi_j$ satisfying Assumption \ref{assump-main} with the normalisation (\ref{eqn-assum on phase after rescaling}), and any sets $\Lambda_j^* + \mb{d}_0 \subset \Lambda_j$ with $\mb{d}[\Lambda_1^*+\mb{d}_0, \Lambda_2^*+\mb{d}_0] \les
\frac{\mb{d}_0}{\mb{C}_0}$, we have for every $R \g R_0$ and $0<\epsilon \ll 1$
	$$ A_{q,r}(2R) \les (1 + C\epsilon) A_{q,r}(R) + \epsilon^{-C} \mc{H}_2^{\frac{n+1}{2r} - \frac{n}{2}} R^{\frac{1}{q} + \frac{n+1}{2r} - \frac{n+1}{2}}. $$
\end{proposition}

We leave the proof of Propositions \ref{prop-initial induction bound} and \ref{prop-induction step} till Section \ref{sec-the localisation argument}, and now turn to the proof of Theorem \ref{thm-main small scale}.

\begin{proof}[Proof of Theorem \ref{thm-main small scale}] We begin by observing that an application of Proposition \ref{prop-induction step} with $R= 2^m R_0$ and $\epsilon = \frac{1}{C} 2^{\frac{m}{2C}(\frac{1}{q} + \frac{n+1}{2r} - \frac{n+1}{2} )}$ implies that
	$$ A_{q,r}\big(2^m R_0\big)\les (1+  2^{\frac{m}{4C}(\frac{1}{q} + \frac{n+1}{2r} - \frac{n+1}{2} )})A_{q, r}\big(2^{m-1} R_0\big) + C^C\mc{H}_2^{\frac{n+1}{2r} - \frac{n}{2}} R_0^{\frac{1}{q} + \frac{n+1}{2r} - \frac{n+1}{2} } 2^{ \frac{m}{2}(\frac{1}{q} + \frac{n+1}{2r} - \frac{n+1}{2}  )}. $$
Since $\frac{1}{q} + \frac{n+1}{2r} < \frac{n+1}{2}$ both error terms decay in $m$. In particular, after $m$ applications of Proposition \ref{prop-induction step} we have\footnote{This exploits the bound
	$$ ( 1 + C 2^{-\frac{\alpha}{C}m}) \times  ( 1 + C 2^{-\frac{\alpha}{C}(m-1)})  \times \cdots \times  ( 1 + C)  \lesa 1$$
which follows by taking logs, and recalling the elementary estimate $\log(1+x) \les x$. }
	$$ A_{q,r}\big(2^m R_0\big) \lesa A_{q,r}(R_0) +\mc{H}_2^{\frac{n+1}{2r} - \frac{n}{2}} R_0^{\frac{1}{q} + \frac{n+1}{2r} - \frac{n+1}{2} } $$
where the implied constant depends only on $q$, $r$, the dimension $n$, and $\mb{C}_0$. Applying Proposition \ref{prop-initial induction bound} and using the definition of $R_0$, we conclude that for every $m\in \NN$
	$$ A_{q,r}\big(2^m R_0\big) \lesa \mb{d}_0^{n+1 - \frac{n+1}{r} - \frac{2}{q}} \mc{H}_2^{\frac{1}{2}-\frac{1}{q}}.$$
Unpacking the definition of $A_{q,r}\big(2^m R_0\big)$, and letting $m \rightarrow \infty$ we obtain Theorem \ref{thm-main small scale}.
\end{proof}

\section{The localisation argument}\label{sec-the localisation argument}

In this section we reduce the proof of Propositions \ref{prop-initial induction bound} and \ref{prop-induction step} to obtaining a decomposition of the $\Phi_j$-waves, into waves which are concentrated at smaller scales, in the sense that the $L^q_t L^r_x$ norm on a cube of diameter $R$, can be controlled by the same norm at the smaller scale $\frac{R}{2}$. This decomposition requires a variant of the wave table construction introduced by Tao in \cite{Tao2001b}. Recall that, given a collection of functions $(u^{(B)})_{B\in \mc{Q}_r(Q)}$ (or a wave table in the notation of Tao), we defined the corresponding \emph{quilt} $[u^{(\cdot)}]$ to be the sum
	$$ [u^{(\cdot)}] = \sum_{B \in \mc{Q}_r(Q)} \ind_B |u^{(B)}|.$$
Roughly speaking, we decompose $u = \sum_{B\in \mc{Q}_r(Q)} u^{(B)}$ into waves which are ``concentrated'' in $B$. The portion of $u^{(B)}$ away from the cube $B$ has additional decay, and can be treated as an error term. On the other hand, the quilt $[u^{(\cdot)}]$ contains the part of $u$ which is concentrated on $B$, and does not decay in $R$. However it is concentrated at a smaller scale than $u$, and this allows us to exploit the induction assumption. \\

The key bound, and the part of Theorem \ref{thm-main small scale} which requires the most work, is the following theorem.

\begin{theorem}\label{thm-general wave table decomposition}
Let  $\mb{C}_0>0$,  $1\les q, r \les 2$, and $\frac{1}{q} +\frac{n+1}{2r} < \frac{n+1}{2}$. Assume that we have $\mb{d}_0>0$, open sets $\Lambda_j \subset \RR^n$, and phases $\Phi_j$ satisfying $(\mb{A1})$, $(\mb{A2})$, and the normalisation $(\ref{eqn-assum on phase after rescaling})$. Let $Q_R$ be a cube of diameter $R \g R_0$. Then for any $0<\epsilon\ll 1$, any $\Phi_1$-wave $u$, and any $\Phi_2$-wave $v$ with
		$$ \supp \widehat{u} + \frac{\mb{d}_0}{2}  \subset \Lambda_1, \qquad \supp \widehat{v} + \frac{\mb{d}_0}{2} \subset \Lambda_2$$
there exist a cube $Q$ of diameter $2R$ such that we have a decomposition
	$$ u = \sum_{B \in \mc{Q}_{\frac{2R}{4^{M}}}(Q)} u^{(B)}, \qquad v = \sum_{B' \in \mc{Q}_{\frac{R}{2}}(Q)} v^{(B')}$$
where $M\in \NN$ with  $4^{-M} < \mc{H}_2 \les 4^{1-M}$, and $u^{(B)}$ is a $\Phi_1$-wave, $v^{(B')}$ is a $\Phi_2$-wave,  with the support properties
	$$ \supp \widehat{u}^{(B)} \subset \supp \widehat{u} + 2\big(2 \mc{H}_2 R\big)^{-\frac{1}{2}}, \qquad \supp \widehat{v}^{(B')} \subset \supp \widehat{v} +  2 \big(2\mc{H}_2 R\big)^{-\frac{1}{2}}.$$
Moreover we have the energy bounds
	$$ \Big(\sum_{B \in \mc{Q}_{\frac{2R}{4^{M}}}(Q)} \|u^{(B)}\|_{L^\infty_t L^2_x}^2 \Big)^\frac{1}{2} \les (1 +C\epsilon) \| u \|_{L^\infty_t L^2_x}$$
	 $$\Big(\sum_{B' \in \mc{Q}_{\frac{R}{2}}(Q)} \|v^{(B')}\|_{L^\infty_t L^2_x}^2 \Big)^\frac{1}{2} \les (1 + C \epsilon) \| v \|_{L^\infty_t L^2_x}$$
and bilinear estimate
	$$ \| u v \|_{L^q_t L^r_x(Q_R)} \les (1+C\epsilon) \big\| \big[u^{(\cdot)}\big] \big[v^{(\cdot)}\big] \big\|_{L^q_t L^r_x(Q)} + \epsilon^{-C} \mc{H}_2^{\frac{n+1}{2r} - \frac{n}{2}} R^{\frac{1}{q} + \frac{n+1}{2r} - \frac{n+1}{2}} \|u \|_{L^\infty_t L^2_x} \| v \|_{L^\infty_t L^2_x}$$
where the constant $C$ depends only on $\mb{C}_0$, $q$, $r$, and $n$.
\end{theorem}

The proof of Theorem \ref{thm-general wave table decomposition} will take up a large part of the remaining sections to follow, and is left to Section \ref{sec-proof of key bilinear est}. \\

The proof of the initial induction bound, also requires the following somewhat classical bilinear $L^2_{t,x}$ estimate.

\begin{theorem}\label{thm-classical bilinear L2 estimate}
Let $j\in \{1, 2\}$, $\Lambda_j \subset \RR^n$, and assume that the phases $\Phi_j \in C^1(\Lambda_j)$ satisfy the transversality condition $(\ref{eqn-main trans assump})$ for some $\mb{C}_0>0$. If $u$ is a $\Phi_1$-wave, and $v$ is a $\Phi_2$-wave we have
        $$ \| u v \|_{L^2_{t,x}(\RR^{1+n})} \les \big(C_0 \mc{V}_{max}\big)^{-\frac{1}{2}} \big( \mb{d}[\supp \widehat{u}, \supp \widehat{v}]\big)^{\frac{n-1}{2}} \|u \|_{L^\infty_t L^2_x} \| v \|_{L^\infty_t L^2_x}.$$
\end{theorem}
\begin{proof} A computation gives the identities
    $$ \widetilde{(uv)}(\tau, \xi) = \int_{\Sigma_1(\tau, \xi)} \frac{\widehat{f}(\eta)\widehat{g}(\xi - \eta)}{|\nabla \Phi_1(\eta) - \nabla \Phi_2(\xi - \eta)|} \, d\sigma(\eta) $$
and
    $$ \int_\RR \int_{\Sigma_1(\tau,\xi)} \frac{|\widehat{f}(\eta)\widehat{g}(\xi - \eta)|^2}{|\nabla \Phi_1(\eta) - \nabla \Phi_2(\xi-\eta)|} d\sigma(\eta) \, d\tau = \big\| \widehat{f}(\eta)\widehat{g}(\xi - \eta) \big\|_{L^2_\eta(\RR^n)}^2$$
where $d \sigma(\eta)$ denotes the surface measure\footnote{ Explicitly, in the region $|\p_1\Phi_1 - \p_1 \Phi_2| \approx |\nabla \Phi_1 - \nabla \Phi_2|$, there exists a function $\psi(\tau, \xi, \eta'):\RR \times \RR^n \times \RR^{n-1} \rightarrow \RR$ such that
            $$ \Phi_1\big(\xi - (\psi, \eta')\big) + \Phi_2(\psi, \eta')  = \tau.$$
Thus we can write the surface as a graph $\Sigma_2(\mathfrak{h})= \{ (\psi(\eta'), \eta') \in (h-\Lambda_1) \cap \Lambda_2\}$, and hence the surface measure is then $ d\sigma(\eta) = \sqrt{1 + |\nabla_{\eta'} \psi|^2} d\eta' = \frac{|\nabla\Phi_1 - \nabla \Phi_2| }{|\p_1 \Phi_1 - \p_1 \Phi_2|} d\eta'$.} on $\Sigma_1(\tau, \xi)$. Consequently bound follows by an application of H\"{o}lder's inequality and the transversality condition (\ref{eqn-main trans assump}).
\end{proof}

The proof of the Propositions \ref{prop-initial induction bound} and \ref{prop-induction step} is now a consequence of Theorems \ref{thm-general wave table decomposition} and \ref{thm-classical bilinear L2 estimate}.

\begin{proof}[Proof of Proposition \ref{prop-initial induction bound}]
By definition of $A_{q,r}(R)$ and $R_0$, it is enough to prove that for every cube $Q$ of diameter $R \g R_0$ and any $\Phi_1$-wave $u$, and $\Phi_2$-wave $v$ satisfying the support conditions
	$$ \supp \widehat{u} \subset \Lambda_1^* + \frac{\mb{d}_0}{2}, \qquad \supp \widehat{v} \subset \Lambda_2^* + \frac{\mb{d}_0}{2}$$
we have
	\begin{equation}\label{eqn-prop initial induc bound-main est} \| u v \|_{L^q_t L^r_x(Q)}\lesa \Big( \mb{d}_0^{(n-1)(1-\frac{1}{r})} \mc{H}_2^{\frac{1}{r} - \frac{1}{2}} R^{\frac{1}{q} + \frac{1}{r} -1} + \mc{H}_2^{\frac{n+1}{2r} - \frac{n}{2}} R^{\frac{1}{q} + \frac{n+1}{2r} - \frac{n+1}{2}}\Big)\|u\|_{L^\infty_t L^2_x} \|v \|_{L^\infty_t L^2_x}. \end{equation}
An application of Theorem \ref{thm-general wave table decomposition} gives
	$$ \| u v \|_{L^q_t L^r_x(Q)} \lesa \big\| [u^{(\cdot)}] [v^{(\cdot)}] \big\|_{L^q_t L^r_x(Q')} +   C  \mc{H}_2^{\frac{n+1}{2r}-\frac{n}{2}} R^{\frac{1}{q} + \frac{n+1}{2r} - \frac{n+1}{2}}\| u \|_{L^\infty_t L^2_x} \|v \|_{L^\infty_t L^2_x}$$
where $Q'$ is a cube of diameter $2R$. To control the quilt term, we first observe that Theorem \ref{thm-classical bilinear L2 estimate} together with the energy estimate for $u^{(B)}$ and $v^{(B)}$ implies the $L^2$ bound
\begin{align*}
   \big\| [u^{(\cdot)}] [v^{(\cdot)}] \big\|_{L^2_{t,x}(Q')}^2 &\les \sum_{B\in \mc{Q}_{\frac{\mc{H}_2 R}{2}}(Q')} \sum_{B' \in \mc{Q}_{\frac{R}{2}}(Q')} \big\| u^{(B)} v^{(B')} \big\|_{L^2_{t,x}}^2 \\
   &\lesa \mb{d}_0^{n-1} \|u \|_{L^\infty_t L^2_x}^2 \|v  \|_{L^\infty_t L^2_x}^2
\end{align*}
where we used the Fourier support conditions on $u^{(B)}$ and $ v^{(B')}$ to deduce that
    $$ \mb{d}\big[\supp \widehat{u}^{(B)}, \supp \widehat{v}^{(B')} \big] \les \mb{d}\big[\Lambda_1^* + \mb{d}_0,  \Lambda_2^* + \mb{d}_0\big] \lesa \mb{d}_0.$$
On the other hand, simply applying H\"{o}lder's inequality gives
		\begin{align*} \big\| [u^{(\cdot)}] [v^{(\cdot)}] \big\|_{L^2_t L^1_x(Q')}
			&\les \big\| [u^{(\cdot)}] \big\|_{L^2_{t,x}(Q')} \big\| [v^{(\cdot)}] \big\|_{L^\infty_t L^2_x(Q')} \\
			& \lesa (R \mc{H}_2)^\frac{1}{2} \Big(\sum \| u^{(B)} \|_{L^\infty_t L^2_x}^2 \Big)^\frac{1}{2} \Big(\sum \| v^{(B')} \|_{L^\infty_t L^2_x}^2 \Big)^\frac{1}{2}\\
			&\lesa (R \mc{H}_2)^\frac{1}{2} \|u\|_{L^\infty_t L^2_x} \| v \|_{L^\infty_t L^2_x}.
		\end{align*}
Applying Holder in the $t$ variable, and interpolating between the $L^2_t L^1_x$ and $L^2_{t,x}$ bounds, then gives
    $$ \big\| [u^{(\cdot)}] [v^{(\cdot)}] \big\|_{L^q_t L^r_x(Q')} \lesa R^{\frac{1}{q} - \frac{1}{2}} \big\| [u^{(\cdot)}] [v^{(\cdot)}] \big\|_{L^2_t L^r_x(Q')} \lesa \mb{d}^{(n-1)(1-\frac{1}{r})}_0 \mc{H}_2^{\frac{1}{r} - \frac{1}{2}} R^{\frac{1}{q}+\frac{1}{r} -1} \|u\|_{L^\infty_t L^2_x} \| v \|_{L^\infty_t L^2_x}.  $$
Therefore (\ref{eqn-prop initial induc bound-main est}) follows.
\end{proof}

A similar application of Theorem \ref{thm-general wave table decomposition} gives Proposition \ref{prop-induction step}.

\begin{proof}[Proof of Proposition \ref{prop-induction step}]
Let $R\gg (\mb{d}_0^2 \mc{H}_2)^{-1}$. Let $u$ be a $\Phi_1$-wave and $v$ be a $\Phi_2$-wave with the support condition
	$$ \supp \widehat{u} \subset \Lambda_1^* + 4(\mc{H}_2 2R)^{-\frac{1}{2}}, \qquad \supp \widehat{v}  \subset \Lambda_2^* + 4(\mc{H}_2 2R)^{-\frac{1}{2}} .$$
It is enough to consider the normalised case $\|u\|_{L^\infty_t L^2_x} = \| v\|_{L^\infty_t L^2_x} = 1$.
Our goal is to show that for every cube $Q$ of diameter $2R$ we have
	\begin{equation}\label{eqn-prop induction step-goal} \| uv \|_{L^q_t L^r_x(Q)} \les (1 + C \epsilon) A_{q,r}(R)  +   \epsilon^{-C} \mc{H}_2^{\frac{n+1}{2r} - \frac{n}{2}} R^{\frac{1}{q} + \frac{n+1}{2r} - \frac{n+1}{2}}.
	\end{equation}
An application of Theorem \ref{thm-general wave table decomposition} gives a cube $Q'$ of diameter $4R$, and wave tables $(u^{(B)})_{B \in \mc{Q}_{\mc{H}_2 R}(Q')}$, $(v^{(B')})_{B'\in \mc{Q}_{R}(Q')}$ such that
	$$ \| u v\|_{L^q_t L^r_x(Q)} \les (1 + C\epsilon) \big\| [u^{(\cdot)}] [v^{(\cdot)}] \big\|_{L^q_t L^r_x(Q')} +  \epsilon^{-C}  \mc{H}_2^{\frac{n+1}{2r} - \frac{n}{2}} R^{\frac{1}{q} + \frac{n+1}{2r} - \frac{n+1}{2}}$$
and the support properties
	\begin{equation}\label{eqn-thm bilinear est-support assump on wave table} \supp \widehat{u}^{(B)} \subset \Lambda_1^* + 4 (\mc{H}_2 R) ^{-\frac{1}{2}}, \qquad \qquad \supp \widehat{v}^{(B')} \subset \Lambda_2^* + 4 (\mc{H}_2 R)^{-\frac{1}{2}}  \end{equation}
where we used the support assumptions on $\widehat{u}$ and $\widehat{v}$. Let $B' \in \mc{Q}_{R}(Q')$ and define the vector valued $\Phi_1$-wave $U^{(B')}= ( u^{(B)} )_{B \in \mc{Q}_{\mc{H}_2 R}(B')}$. It is easy to check that $U^{(B')}$ is a $\Phi_1$-wave (since $u^{(B)}$ are maps into $\ell^2_c(\ZZ)$, after relabeling, $U^{(B')}$ is also a map into $\ell^2_c(\ZZ)$). In particular, for every $B'\in \mc{Q}_{R}(Q')$ we have a $\Phi_1$-wave $U^{(B')}$ and a $\Phi_2$-wave $v^{(B')}$ satisfying the correct support assumptions to apply the definition of $A_{q,r}(R)$. Hence, if we observe that
	\begin{align*} [u^{(\cdot)}] &= \sum_{B' \in \mc{Q}_{R}(Q')} \sum_{B \in \mc{Q}_{\mc{H}_2 R}(B')} \ind_B |u^{(B)}| \\
			&\les \sum_{B' \in \mc{Q}_{R}(Q')} \ind_{B'} \Big( \sum_{B \in \mc{Q}_{\mc{H}_2 R}(B')} |u^{(B)}|^2 \Big)^\frac{1}{2} =  \sum_{B' \in \mc{Q}_{R}(Q')} \ind_{B'} |U^{(B')}|
	\end{align*}
then, applying the definition of $A_{q,r}(R)$, we deduce that
	\begin{align*}
	   \big\| [ u^{(\cdot)} ] [ v^{(\cdot)}] \big\|_{L^q_t L^r_x(Q')} &\les \sum_{B' \in \mc{Q}_R(Q')} \| U^{(B')} v^{(B')} \|_{L^q_t L^r_x(B')} \\
	   &\les A_{q,r}(R) \Big( \sum_{B' \in \mc{Q}_R(Q')} \|U^{(B')}\|_{L^\infty_t L^2_x}^2 \Big)^\frac{1}{2}\Big( \sum_{B' \in \mc{Q}_R(Q')} \|v^{(B')}\|_{L^\infty_t L^2_x}^2 \Big)^\frac{1}{2} \\
	   &\les A_{q,r}(R) (1 +   C \epsilon)^2
	\end{align*}
where the last line follows by using the energy inequalities in Theorem \ref{thm-general wave table decomposition}.
\end{proof}

\begin{remark}
Note that proof of Proposition \ref{prop-induction step} exploited the fact that the definition of $A_{q, r}(R)$ applies to \emph{vector} valued waves. In fact, this is essentially the only step in the proof of Theorem \ref{thm-main small scale} in which the freedom to use vector valued waves is crucial.
\end{remark}

We have now reduce the problem of proving Theorem \ref{thm-main small scale} to obtaining the decomposition contained in Theorem \ref{thm-general wave table decomposition}. The proof of Theorem \ref{thm-general wave table decomposition} requires the full range of tools used to prove bilinear restriction estimates, namely a sharp wave packet decomposition, the wave table construction due to Tao, geometric information on the conic surfaces $\mc{C}_j(\mathfrak{h})$, and energy estimates across transverse surfaces (which can be thought of as a replacement for the combinatorial Kakeya type inequalities used in the original argument of Wolff).

\section{Wave Packets}\label{sec-wave packets}

 The standard approach to constructing wave packets is to localise on both the spatial and Fourier side, up to the scale given by the uncertainty principle. However, as we need to carefully control the constants in the energy estimates, we use a more refined construction originally due to Tao \cite{Tao2001b} and extended to the case of general phases by Bejenaru \cite{Bejenaru2016b}.  Let $\nu \in \s(\RR^n)$ be a positive smooth function, such that $\supp \widehat{\nu}  \subset \{ |\xi| \les 1\}$,  and for all $x \in \RR^n$
    $$ \sum_{k \in \ZZ^n} \nu( x - k) = 1.  $$
Let $A = \{ \xi = (\xi_1, \dots , \xi_n) \in \RR^n \mid  |\xi_j| \les \frac{1}{2}, j=1, \dots, n\}$ denote the unit cube centered at the origin, and define
        $$ \rho(\xi) =  \fint_{|z| < 1} \ind_A(\xi - z) \,dz $$
 where we take $\fint_\Omega = \frac{1}{|\Omega|} \int_Q$. To a phase space point $\gamma = (x_0, \xi_0) \in \mc{X}_{r, \epsilon}$ and $f \in L^2_x(\RR^n)$, define the phase-space localisation operator
	\begin{equation}\label{eqn-phase space localisation operator} \big(L_\gamma f\big)(x) = \nu\Big( \frac{\epsilon^2}{r} (x - x_0)\Big) \Big[\rho\big(  r ( -i \nabla - \xi_0) \big) f\Big](x). \end{equation}

We have the following basic properties.

\begin{lemma}[{Properties of $L_\gamma$ \cite[Lemma 15.2]{Tao2001b}, \cite[Lemma 4.1]{Bejenaru2016b}}]\label{lem-prop of Lgamma} Let $r>0$, $0<\epsilon\les 1$, and $f \in L^2(\RR^n)$. Then
    $$ \sum_{\gamma \in \mc{X}_{r, \epsilon}} L_\gamma f = f, \qquad \supp \widehat{L_\gamma f} \subset \big\{ |\xi - \xi(\gamma)| \les \tfrac{n+1}{r} \big\} \bigcap \big( \supp \widehat{f} + \frac{1}{r}\big). $$
Moreover, if $(m_{k, \gamma})_{k, \gamma}$ is a positive sequence with $\sup_\gamma \sum_k m_{k, \gamma} \les 1$, then
    $$ \bigg( \sum_k \Big\| \sum_{\gamma \in \mc{X}_{r,\epsilon}} m_{k, \gamma} L_\gamma f \Big\|_{L^2_x}^2 \bigg)^\frac{1}{2} \les (1 +  C \epsilon)   \| f\|_{L^2} $$
where the constant $C$ depends only on the dimension $n$.
\end{lemma}
\begin{proof} As in \cite{Bejenaru2016b}, we adapt the method of Tao \cite{Tao2001b}. However, as we use a slightly more direct argument, we include the details. Clearly, by construction we have
	$$ \sum_{\gamma \in \mc{X}_{r,\epsilon}} L_\gamma f = c_n \int_{\RR^n}  \fint_{|z|<\frac{1}{2}}  \sum_{\xi_0 \in \frac{1}{r} \ZZ^n} \ind_A\big( r(\xi-\xi_0) - z\big) \,dz\,\, \widehat{f}(\xi) e^{ i \xi \cdot x} d\xi = f.    $$
On the other hand, since $\supp \widehat{\nu} \subset \{ |\xi|\les 1\}$ and $\supp \widehat{\rho}\subset \{|\xi| \les \frac{2+\sqrt{n}}{2}\}$, we have $\supp \widehat{L_\gamma f} \subset \frac{\epsilon^2}{r} + \supp \widehat{f} \cap \{ |\xi - \xi_0| \les \frac{n}{r}\}$.
It remains to prove the more delicate energy type estimate. To simplify the notation somewhat, we take 	
	$$\nu_{x_0} = \nu( \frac{\epsilon^2}{r}(x-x_0)), \qquad \widehat{f}_{\xi_0, z} = \ind_A(r(\xi-\xi_0) - z) \widehat{f}, \qquad \widehat{P^{(\xi_0, z)}_{A^*} g}(\xi) = \ind_{A^*}(r(\xi-\xi_0) - z) \widehat{g}(\xi)$$
where
        $$A^* = \{ \xi = (\xi_1, \dots , \xi_n) \in \RR^n \mid  |\xi_j| \les \tfrac{1}{2}-2\epsilon^2, j=1, \dots, n\}.$$
  Applying Minkowski's inequality and decomposing $f_{(\xi_0, z)} $  gives
    \begin{align}
      \bigg( \sum_k \Big\| \sum_\gamma m_{k, \gamma} L_\gamma f \Big\|_{L^2_x}^2 \bigg)^\frac{1}{2}
        &\les \fint_{|z|<1} \bigg( \sum_k \Big\| \sum_{x_0} \sum_{\xi_0} m_{k, (x_0, \xi_0)} \nu_{x_0}  f_{\xi_0, z} \Big\|_{L^2_x}^2 \bigg)^\frac{1}{2} dz \label{eqn - lem prop of L - initial decomp} \\
        &\les \fint_{|z|<1} \bigg( \sum_k \Big\| \sum_{x_0} \sum_{\xi_0} m_{k, (x_0, \xi_0)} \nu_{x_0} P^{(\xi_0, z)}_{A^*} f_{\xi_0, z} \Big\|_{L^2_x}^2 \bigg)^\frac{1}{2} dz \notag \\
         &\qquad + \fint_{|z|<1} \bigg( \sum_k \Big\| \sum_{x_0} \sum_{\xi_0} m_{k, (x_0, \xi_0)} \nu_{x_0}  \big( 1 - P^{(\xi_0, z)}_{A^*}\big) f_{\xi_0, z} \Big\|_{L^2_x}^2 \bigg)^\frac{1}{2} dz. \notag
    \end{align}
  A computation shows that the Fourier supports of the functions  $\nu_{x_0} P^{(\xi_0, z)}_{A^*} f_{\xi_0, z}$ are disjoint as $\xi_0$ varies.
  Consequently we have the identity
    \begin{align*}
      \bigg( \sum_k \Big\| \sum_{x_0} \sum_{\xi_0} m_{k, (x_0, \xi_0)} \nu_{x_0} P^{(\xi_0, z)}_{A^*} f_{\xi_0, z} \Big\|_{L^2_x}^2 \bigg)^\frac{1}{2}
            &=  \bigg( \sum_k \sum_{\xi_0} \Big\| \sum_{x_0}  m_{k, (x_0, \xi_0)} \nu_{x_0} P^{(\xi_0, z)}_{A^*} f_{\xi_0, z} \Big\|_{L^2_x}^2 \bigg)^\frac{1}{2}.
    \end{align*}
  As $\nu_{x_0}$ and $m_{k, (x_0, \xi_0)}$ are positive, we have
    $$ \sum_k \Big| \sum_{x_0} m_{k, (x_0, \xi_0)} \nu_{x_0} P^{(\xi_0, z)}_{A^*} f_{\xi_0, z}\Big|^2
                \les |P^{(\xi_0, z)}_{A^*} f_{\xi_0, z}| \Big| \sum_{x_0} \nu_{x_0} \sum_k m_{k, (x_0, \xi_0)} \Big|^2 \les |P^{(\xi_0, z)}_{A^*} f_{\xi_0, z}| $$
    and therefore
        $$ \fint_{|z|<1} \bigg( \sum_k \Big\| \sum_{x_0} \sum_{\xi_0} m_{k, (x_0, \xi_0)} \nu_{x_0} P^{(\xi_0, z)}_{A^*} f_{\xi_0, z} \Big\|_{L^2_x}^2 \bigg)^\frac{1}{2} dz \les \| f \|_{L^2}.$$
    It remains to estimate the second term in (\ref{eqn - lem prop of L - initial decomp}). Repeating the above argument, but applying almost orthogonality instead of orthogonality, we see that
        $$\fint_{|z|<1} \bigg( \sum_k \Big\| \sum_{x_0} \sum_{\xi_0} m_{k, (x_0, \xi_0)} \nu_{x_0}  \big( 1 - P^{(\xi_0, z)}_{A^*}\big) f_{\xi_0, z} \Big\|_{L^2_x}^2 \bigg)^\frac{1}{2} dz \les C_n \fint_{|z|<1} \bigg( \sum_{\xi_0} \Big\| \big( 1 - P^{(\xi_0, z)}_{A^*}\big) f_{\xi_0, z} \Big\|_{L^2_x}^2 \bigg)^\frac{1}{2} dz.$$
    An application of Holder's inequality in the $dz$ integral, together with Plancheral reduces the problem to proving that
        $$ \sum_{\xi_0} \fint_{|z|\les 1} \Big[ 1 - \ind_{A^*}\big( r(\xi - \xi_0) - z\big) \Big] \ind_A\big( r(\xi - \xi_0) - z \big) dz \les C_n \epsilon.$$
But this follows from a short computation. Thus lemma follows.
\end{proof}

We can now define the wave packets we use in the proof of Theorem \ref{thm-main small scale}. As is more or less standard, see for instance \cite{Bejenaru2016b, Candy2016}, to define the wave packet associated to a phase space point $\gamma \in \Gamma_j$, we conjugate the phase-space localisation operator $L_\gamma$ with the flow $e^{ i t \Phi_j(- i \nabla)}$. More precisely:

\begin{definition}[Wave Packets]\label{def:wave-packets} Let $j=1, 2$, and $r_j \gg \mb{d}_0^{-1}$. If $u \in L^\infty_t L^2_x$ and $\gamma_j \in \Gamma_j$, we define
	$$ \big(\mc{P}_{\gamma_j} u\big)(t) = e^{i t \Phi_j(-i \nabla)} L_{\gamma_j} \Big( e^{  - i t \Phi_j(-i\nabla)} u(t) \Big).$$
\end{definition}

Lemma \ref{lem-prop of Lgamma} has an immediate extension to the wave packets $\mc{P}_{\gamma_j}$.

\begin{lemma}[Orthogonality Properties of Wave Packets]\label{lem-orthog prop of wave packets}
Let $j=1, 2$, $0<\epsilon \les1$, and $r_j \gg  \mb{d}_0^{-1}$. Let $u \in L^\infty_t L^2_x$ with $\supp \widehat{u}(t) + \frac{\mb{d}_0}{2} \subset \Lambda_j$. Then
	$$ u = \sum_{\gamma_j \in \Gamma_j} \mc{P}_ju, \qquad \supp \widehat{\mc{P}_{\gamma_j}u}(t) \subset \{ |\xi - \xi(\gamma_j)|\les \tfrac{n+1}{r_j} \} \bigcap \big( \supp \widehat{u} + \tfrac{1}{r_j}\big).$$
Moreover, if $(m_{k, \gamma_j})_{k, \gamma_j}$ is a positive sequence with $\sup_{\gamma_j \in \Gamma_j} \sum_k m_{k, \gamma_j} \les 1$, then
    \begin{equation}\label{eqn-prop of wave packets-main energy est} \bigg\| \bigg( \sum_k \Big\| \sum_{\gamma_j \in \Gamma_j} m_{k, \gamma_j} \mc{P}_{\gamma_j} u(t) \Big\|_{L^2_x}^2 \bigg)^\frac{1}{2}\bigg\|_{L^\infty_t}  \les (1 + C \epsilon)   \| u\|_{L^\infty_t L^2_x} \end{equation}
where the constant $C$ depends only on the dimension $n$.
\begin{proof} If we observe that the Fourier multiplier $e^{ i t \Phi_j(-i\nabla)}$ is unitary on $L^2_x$, and clearly leaves the Fourier support invariant, then all the required properties follow from Lemma \ref{lem-prop of Lgamma} together with the identity
    $${u=\sum_{\gamma_j \in \mc{X}_{r_j, \epsilon}}  \mc{P}_j u = \sum_{\gamma_j \in \Gamma_j} \mc{P}_j u}$$
which is a consequence of the definition of $\Gamma_j$ and the support assumption on $u$.
\end{proof}
\end{lemma}

More refined properties of the wave packets $\mc{P}_{\gamma_j}u$ are possible if we make further assumptions on $u$. In particular, if $u=e^{it\Phi(-i\nabla)}f$ is a free solution, then the associated  wave packets are concentrated on the tubes $T_{\gamma_j}$. One possible formulation of this statement, which will prove extremly useful in practise, gives a bound in terms of the Hardy-Littlewood maximal function
	$$ \mc{M}[f](x) = \sup_{s>0} \frac{1}{s^n} \int_{|x-y|<s} |f(y)| dy.$$

\begin{lemma}[Concentration of wave packets I]\label{lem-conc of wave packets I}
Let $N \g 2$, $0<\epsilon\les 1$, and $j\in \{1,2\}$. Assume that the phase $\Phi_j \in C^N(\Lambda_j)$ satisfies for every $2\les m \les N$ the bound
	$$  \|\nabla^m \Phi_j\|_{L^\infty(\Lambda_j)} \les \mb{d}_0^{2-m} \mc{H}_j. $$
Let $r_j \gg \mb{d}_0^{-1}$, $\gamma_j=(x_0, \xi_0) \in \Gamma_j$, and let $u=e^{it\Phi_j(-i\nabla)}f$ be a $\Phi_j$-wave. Then for all $(t,x) \in \RR^{1+n}$ we have
	 \begin{equation}\label{eqn-lem conc of wave packets-spatial localisation of wave packets}
       |\mc{P}_{\gamma_j} u(t,x)| \lesa \epsilon^{-2(N+n)} \Big(1+  \frac{|x-x_0 + t \nabla\Phi_j(\xi_0)|}{r_j} \Big)^{-N} \Big( 1 + \frac{|t|}{R}\Big)^N \mathcal{M}\big[  f_{\xi_0} \big](x_0)
    \end{equation}
where $f_{\xi_0} = \rho( r_j(-i\nabla - \xi_0))f$, and the implied constant depends only on $n$ and $N$ (in particular it is independent of the phase $\Phi_j$ and $\mb{d}_0$).
\end{lemma}
\begin{proof}
We start by writing
	\begin{align*}
		\mc{P}_{\gamma_j} u(t, x) &= \int_{\RR^n} \widehat{(L_{\gamma_j} f)}(\xi) e^{ i t \Phi_j(\xi)} e^{ i x \cdot \xi} d\xi \\
		&= \int_{\RR^n} K_{\xi_0}(t, x-y) (L_{\gamma_j} f)(y) dy
	\end{align*}
where the kernel is given by $ K_{\xi_0}(t, x) = \int_{\RR^n} \zeta\big( r_j(\xi - \xi_0)\big) e^{ i t \Phi_j(\xi)} e^{ i x \cdot \xi} d\xi$ and $\zeta \in C^\infty_0(|\xi|<2)$ is a cutoff satisfying $\zeta(\xi) \rho(\xi)  = \rho(\xi) $. The spatial localisation properties of $L_{\gamma_j} f$, together with the standard maximal function inequality
	$$ \int_{\RR^n} ( 1 + s|x-y|)^{-(n+1)} |g(y)| s^n \,dy \lesa \mc{M}[g](x)$$
implies that it suffices to prove the kernel bound
    \begin{equation}\label{eqn-lem conc of wave packets-kernel bound}
        \big| K_{\xi_0}(t, x)\big| \lesa r_j^{-n} \Big( 1 + \frac{ | x + t \nabla\Phi(\xi_0)|}{r_j} \Big)^{-N} \Big( 1 +\frac{|t|}{R}\Big)^N.
    \end{equation}
This bound is clearly true when $| x + t \nabla\Phi(\xi_0)|\les r_j$ by the Fourier support assumption on $\zeta$. On the other hand, an application of integration by parts gives for every $1\les k \les n$
     \begin{align*} \big| K_{\xi_0}(t, x)\big|&= \bigg| \int_{\RR^n} \zeta\big( r_j(\xi - \xi_0)\big) e^{ it [ \Phi_j(\xi) - \Phi_j(\xi_0) -(\xi - \xi_0) \cdot \nabla \Phi_j(\xi_0)]}  e^{ i (x + t \nabla \Phi_j(\xi_0) ) \cdot \xi} d\xi \bigg| \\
     &\les |x_k + t \p_k \Phi_j(\xi_0) |^{-N}  \int_{\RR^n} \Big| \p_k^N \Big( \zeta\big( r_j (\xi - \xi_0)\big) e^{ it [ \Phi_j(\xi) - \Phi_j(\xi_0) -(\xi - \xi_0) \cdot \nabla \Phi_j(\xi_0)]}\Big) \Big| d\xi.
     \end{align*}
In particular, it suffices to show that for every $|\xi - \xi_0| \les 2r_j^{-1}$ we have
      $$
            \Big| \p_k^N \Big( \zeta\big( r_j(\xi - \xi_0)\big) e^{ it [ \Phi_j(\xi) - \Phi_j(\xi_0) -(\xi - \xi_0) \cdot \nabla \Phi_j(\xi_0)]}\Big) \Big| \lesa  r_j^N \Big( 1 +\frac{|t|}{R}\Big)^N.
        $$
This bound is consequence of a somewhat tedious induction argument. For completeness, we sketch the proof here. Let $F=i t ( \Phi_j(\xi) - \Phi_j(\xi_0) -(\xi - \xi_0) \cdot \nabla \Phi_j(\xi_0))$ and
    $$ I_N = e^{ - F} \p_k^N \Big( \zeta\big( r_j(\xi - \xi_0)\big) e^{F }\Big).$$
To compute $I_N$ explicitly, we observe that $I_N = \p_k I_{N-1} + I_{N-1} \p_k F$, which via an induction argument, implies that
        \begin{equation}\label{eqn-lem conc of wave packets-I_N ident} I_N = \p^N_k I_0 +  \sum_{m=1}^N \p_k^{N-m} I_0 \sum_{ \substack{ \ell \in \NN^m \\ \ell_1 + 2 \ell_2 + \dots + m \ell_m = m }} c_{m, \ell} ( \p_k F )^{\ell_1} \dots (\p_k^m F)^{\ell_m} \end{equation}
for some constants $c_{m, \ell} \in \NN$. The assumption (\textbf{A2}) together with $r_j \gg \mb{d}_0$ implies that
	$$ |\p_k  F| \les |t| |\nabla \Phi_j(\xi) - \nabla \Phi_j(\xi_0)| \les \frac{|t|}{R} r_j, \qquad |\p_k^2 F|\les |t| |\nabla^2 \Phi_j(\xi)| \les \frac{|t|}{R} r_j^2$$
and for any $3\les m \les N$
	$$ |\p_k^m F| \les |t| |\nabla^m \Phi_j(\xi)| \les \frac{|t|}{|R|} \mb{d}_0^{2-m} R \|\nabla^2 \Phi_j\|_{L^\infty} \les \frac{|t|}{R} r_j^m.$$
Therefore the bound (\ref{eqn-lem conc of wave packets-kernel bound}) follows from the identity (\ref{eqn-lem conc of wave packets-I_N ident}).
\end{proof}

\begin{remark}
If we restrict attention to the case $|t|\lesa R$ and $\epsilon=1$, then Lemma \ref{lem-conc of wave packets I} is essentially equivalent to the wave packet type decompositions used frequently in previous works. See for instance \cite{Tao2001b, Lee2010, Bejenaru2016b, Candy2016}. It is also important to note that the assumption (\textbf{A2}) is essentially the weakest condition for the required concentration bound (\ref{eqn-lem conc of wave packets-spatial localisation of wave packets}) to hold. This can easily be seen from the identity (\ref{eqn-lem conc of wave packets-I_N ident}), since if we want to obtain (\ref{eqn-lem conc of wave packets-spatial localisation of wave packets}), then we essentially require
		$$ \|\nabla^m \Phi_j\|_{L^\infty} \lesa r_j^{m-2}$$
which is equivalent to having a uniform bound on the quantity
		$$ \Big( \frac{\|\nabla^m \Phi_j\|_{L^\infty}}{\|\nabla^2 \Phi_j\|_{L^\infty}} \Big)^{\frac{1}{m-2}}.$$
and this is precisely the condition appearing in (\textbf{A2}). 		
\end{remark}

Recall that, given a phase space point $\gamma_j \in \Gamma_j$, we have defined the associated weight $w_{\gamma_j, q}$, which penalises the distance from the cube $q \in \mc{Q}_{r_j}(Q)$ and tube $T_{\gamma_j}$, by
	$$ w_{\gamma_j, q} =  \Big( 1 + \frac{|x_q - x_0 + t_q \nabla \Phi_j(\xi_0)|}{r_j}\Big)^{5n}.$$

We now come to the main property of wave packets which we exploit in the sequel.

\begin{proposition}[Concentration of wave packets II]\label{prop - concentration of wave packets II}
Let $j=1, 2$ and $0<\epsilon \les1$. Assume that the phase $\Phi_j \in C^{5n}(\Lambda_j)$ satisfies ($\mb{A2}$).
Let $R \gtrsim  r_j \gg \mb{d}_0^{-1}$ and assume that $u$ is a $\Phi_j$-wave with $\supp \widehat{u} + \frac{\mb{d}_0}{2} \subset \Lambda_j.$
Then for any cube $Q$ of diameter $R$ we have the concentration/orthgonality type estimate
	\begin{equation}\label{eqn-concentration of wave packets II-main conc bound}\sum_{\gamma_j \in \Gamma_j} \sup_{q \in \mc{Q}_{r_j}(Q)} w_{\gamma_j, q}^2 \|   \chi_q \mc{P}_{\gamma_j} u\|_{L^2_{t, x}}^2\lesa \epsilon^{-14n}  r_j \| u \|_{L^\infty_t L^2_x}^2.
	\end{equation}
\end{proposition}
\begin{proof}
Let $N=5n$ and $Q$ be a cube of diameter $R$. By translation invariance, we may assume that $Q$ is centred at the origin. Write $u = e^{it\Phi_j(-i\nabla)}f$ for some $\ell^2_{c}(\ZZ)$ valued map $f$. The identity
	 $$|x-x_0+t\nabla \Phi_j(\xi_0)|=|(t,x) - (0,x_0) - t(1, -\nabla \Phi_j(\xi_0))|,$$
together with an application of Lemma \ref{lem-conc of wave packets I} implies that for any $q\in \mc{Q}_{r_j}(Q)$ and $\gamma_j = (x_0, \xi_0) \in \Gamma_j$
	\begin{align*}
	 \|   \chi_q &\mc{P}_{\gamma_j} u\|_{L^2_{t, x}}^2 \\
	 &\lesa \epsilon^{-(N+n)} |\mc{M}[f_{\xi_0}](x_0)|^2 \Big\| \chi_q(t,x) \Big(1+ \frac{|x-x_0 + t \nabla\Phi_j(\xi_0)|}{r_j} \Big)^{-N} \Big( 1 + \frac{|t|}{R}\Big)^N    \Big\|_{L^2_{t,x}}^2\\
	  &\lesa \epsilon^{-(N+n)}  |\mc{M}[f_{\xi_0}](x_0)|^2 \Big(1+ \frac{|x_q-x_0 + t_q \nabla\Phi_j(\xi_0)|}{r_j} \Big)^{-2N} \Big\| \chi_q(t,x) \Big(1+ \frac{|(t,x) - (t_q,x_q)|}{r_j} \Big)^{2N}  \Big\|_{L^2_{t,x}}^2\\
      &\lesa  \epsilon^{-(N+n)} r_j^{n+1} |\mc{M}[f_{\xi_0}](x_0)|^2 \Big(1+ \frac{|x_q-x_0 + t_q \nabla\Phi_j(\xi_0)|}{r_j} \Big)^{-2N}
	\end{align*}
where we used (\ref{eqn-assum on phase after rescaling}), the rapid decay of $\chi_q$, and the bound
        $$ \Big( 1 + \frac{|t|}{R}\Big) \les 2 \Big( 1 + \frac{|t-t_q|}{R}\Big) \lesa  \Big( 1 + \frac{|t-t_q|}{r_j}\Big)$$
which follows directly from the conditions $|t_q|\les R$ and $r_j \lesa R$. Therefore, we deduce that
    $$ \sum_{\gamma_j \in \Gamma_j} \sup_{q \in \mc{Q}_{r_j}(Q)} w_{\gamma_j, q}^2 \|   \chi_q \mc{P}_{\gamma_j} u\|_{L^2_{t, x}}^2 \lesa r_j^{n+1} \epsilon^{-2(N+n)} \sum_{(x_0, \xi_0) \in \Gamma_j} |\mc{M}[f_{\xi_0}](x_0)|^2.  $$
To complete the proof of (\ref{eqn-concentration of wave packets II-main conc bound}), we recall that the Fourier localisation of $f_{\xi_0}$ together with the maximal inequality gives
   $$ r_j^n \sum_{(x_0, \xi_0) \in \Gamma_j}  \big( \mathcal{M}[f_{\xi_0}](x_0)\big)^2 \lesa \epsilon^{-2n} \sum_{\xi_0} \|\mathcal{M}[f_{\xi_0}] \|_{L^2_x}^2 \lesa  \epsilon^{-2n} \sum_{\xi_0} \| f_{\xi_0} \|_{L^2_x}^2 \approx \epsilon^{-2n} \| f\|_{L^2_x}^2. $$
\end{proof}

\begin{remark}
  Note that technical condition $r_j \lesa R$ is natural since (\ref{eqn-concentration of wave packets II-main conc bound}) only controls $u$ on (essentially) a cube of size $R$. If $r_j\gg R$ (which corresponds to $R\ll \|\nabla^2\Phi_j\|_{L^\infty}$) then each wave packet is supported in an area \emph{larger} than $Q$. In other words, $u$ cannot be localised to scales $R\ll\|\nabla^2 \Phi_j \|_{L^\infty}$.
\end{remark}

\begin{remark}
It is worth comparing (\ref{eqn-concentration of wave packets II-main conc bound}) to previous results in the literature where a \emph{local} in time version was proven, in the sense that the function $\chi_q$ was replaced by the sharp cutoff $\ind_{Cq}$. The slightly more refined global version (\ref{eqn-concentration of wave packets II-main conc bound}) (which also contains the region $|t| \gg R$, albeit with a rapidly decaying weight $t^{-N}$) allows us to avoid a number of additional error estimates, which would otherwise complicate the analysis, particularly in the general phase case that we consider here.
\end{remark}

\section{Geometric Consequences of Assumptions on Phases}\label{sec-geometric consequences of phase}

In this section, following the approach in \cite{Candy2016}, we give two important consequences of the transversality/curvature assumption (\textbf{A1}). The first concerns the surfaces $\Sigma_j(\mathfrak{h})$ and only requires the transverality condition (\ref{eqn-main trans assump}), while the second gives a crucial transverality condition concerning the conical set $\mc{C}_j(\mathfrak{h})$ and requires the full strength of the assumption (\textbf{A1}). In fact, this is the only place where the full generality of (\textbf{A1}) is required, rather than just the immediate consequences (\ref{eqn-main trans assump}) and (\ref{eqn-main curv assump}). \\

Recall that give $\mathfrak{h} \in \RR^{1+n}$ we have defined the surfaces $ \Sigma_j(\mathfrak{h}) = \{ \xi \in \Lambda_j\cap (h-\Lambda_k) \mid \Phi_j(\xi) + \Phi_k(h-\xi) = a \}$. The transversality assumption (\ref{eqn-main trans assump}) implies that these surfaces are well-behaved.

\begin{lemma}[Foliation Properties of $\Sigma_j(\mathfrak{h})$]\label{lem - folation properties of surface Sigma}
  Let $\mathfrak{h} = (a, h) \in \RR^{1+n}$ and $0<\mb{C}_0<1$. Suppose that $\Lambda_j$ is open, and that $\Phi_j \in C^2(\Lambda_j)$ satisfies the transverality condition (\ref{eqn-main trans assump}) and the normalisation condition (\ref{eqn-assum on phase after rescaling}).  Assume $\xi + \frac{3}{\mb{C}_0 r} \in \Lambda_1$ and $\eta + \frac{3}{\mb{C}_0 r} \in \Lambda_2$ with
        $$|\xi + \eta - h|\les \frac{1}{r}, \qquad\qquad \big| \Phi_1(\xi) + \Phi_2(\eta) - a \big| \les \frac{1}{r}.$$
  Then $\xi \in \Sigma_1(\mathfrak{h}) + \frac{3}{\mb{C}_0r }$ and $\eta \in \Sigma_2(\mathfrak{h}) + \frac{3}{\mb{C}_0r}$.
\end{lemma}
\begin{proof} Let $F(\xi) = \Phi_1(\xi) + \Phi_2(h-\xi) -a$, we may assume that $F(\xi)\g0$ since otherwise we can simply replace $F$ with $-F$ in the argument below. Let $\xi(s)$ be the solution to
				\begin{align*}
						\p_s \xi(s) &= - \frac{\nabla F( \xi(s) )}{|\nabla F(\xi(s))|}\\
						     \xi(0) &= \xi.
				\end{align*}
 The conditions on $\xi$ and $\eta$ imply that $t(h-\xi) + (1-t)\eta \in \Lambda_2$ for $0<t<1$, and hence $F(\xi)= |\Phi_1(\xi) + \Phi_2(\eta) - a + \Phi_2(h-\xi) - \Phi_2(\eta)|\les  \frac{3}{r}$. Therefore the transversality condition (\ref{eqn-main trans assump}) gives the bound
$$ F(\xi(s)) = F(\xi) - \int_0^s |\nabla F(\xi(s'))| ds'\les  \frac{3}{r} - s \textbf{C}_0 $$
and consequently we see that $F(\xi(s^*))=0$ for some $0\les s^*\les \frac{3}{r \textbf{C}_0 }$. Since $|\xi(s^*) - \xi|\les s^* \les \frac{3}{r \textbf{C}_0 }$ and $\xi(s^*) \in \Sigma_1(\mathfrak{h})$ we have $\xi \in \Sigma_1(\mathfrak{h}) + \frac{3}{r\textbf{C}_0 }$ as required. Replacing $\Phi_1$ with $\Phi_2$ in the above argument, gives $\eta \in \Sigma_2(\mathfrak{h}) + \frac{3}{r\textbf{C}_0 }$.
\end{proof}

Recall that, given $\mathfrak{h} \in \RR^{1+n}$, we have defined  the conical hypersurface
     $$ \mc{C}_j(\mathfrak{h}) = \big\{ \big(s, - s\nabla \Phi_j(\xi) \big) \,\,\big| \,\,s\in \RR,  \xi \in \Sigma_j(\mathfrak{h}) \, \big\}.$$
The next result we require shows that the normal direction to the surface $\tau = \nabla \Phi_k$ always intersects the conic surface $\mc{C}_j(\mathfrak{h})$ transversely.

\begin{lemma}[Transversality of $\mc{C}_j(\mathfrak{h})$]\label{lem - surface C transverse}
Let $\mathfrak{h} \in \RR^{1+n}$. Suppose that $\Lambda_j\subset \RR^n$ is an open set, and that $\Phi_j \in C^2(\Lambda_j)$ satisfies $(\mb{A1})$ and the normalisation condition (\ref{eqn-assum on phase after rescaling}). For  $\{j, k\} = \{1, 2\}$,   $p, q \in \mc{C}_j(h)$,  and $\eta \in \Lambda_k$ we have
	$$  \big| (p-q) \wedge \big( 1 , - \nabla \Phi_k(\eta)\big) \big| \g  \frac{\mb{C}_0}{3}|p-q|.$$
      \end{lemma}
      \begin{proof}
We follow/repeat the argument given in \cite[Lemma 2.7]{Candy2016}. Let $w, w', w'' \in \RR^n$. A computation shows that for every $ {v \in \text{span}\{ (1, w), (0, w-w') \}}$ we have
	$$|v \wedge (1, w'')| \g \frac{ |(w - w')\wedge(w - w'') |}{ (1 +|w|) |w-w'|} |v|. $$
Let $\xi, \xi' \in \Sigma_j(\mathfrak{h})$ and $\eta \in \Lambda_k$. If we take $w = - \nabla \Phi_j(\xi)$, $w' = - \nabla \Phi_j(\xi')$, and $w'' =- \nabla \Phi_k(\eta)$ in the above bound, we deduce that
	$$ \big| v \wedge \big( 1, - \nabla \Phi_k(\eta) \big) \big| \g  \frac{|(\nabla \Phi_j(\xi) - \nabla \Phi_j(\xi'))\wedge ( \nabla \Phi_j(\xi) - \nabla \Phi_k(\eta))| }{ (1 + \| \nabla \Phi_j\|_{L^\infty(\Lambda_j)}) |\nabla \Phi_j(\xi) - \nabla \Phi_j(\xi') |}|v| $$
for every $v \in \text{span}\{ (1, - \nabla \Phi_j(\xi)), (0, \nabla \Phi_j(\xi) - \nabla \Phi_j(\xi'))\}$. The claimed bound now follows by applying (\textbf{A1}), and observing that given any $p, q\in \mc{C}_j(\mathfrak{h})$ we can write
	$$(p - q) = (r-r') (1, - \nabla \Phi_j(\xi)) + r' (0, \nabla \Phi_j(\xi) - \nabla \Phi_j(\xi') ),$$
for some $\xi, \xi' \in \Sigma_j(\mathfrak{h})$ and $r, r' \in \RR$.
\end{proof}

\section{Energy Estimates Across $\mc{C}_j(\mathfrak{h})$} \label{sec-energy estimates across C_j}

In this section we prove a key energy bound across the conic surface $\mc{C}_j(\mathfrak{h})$. This is the key step in the proof of Theorem \ref{thm-main small scale} where the full strength of (\textbf{A1}) is required, elsewhere it is essentially enough to instead consider the weaker conditions (\ref{eqn-main trans assump}) and (\ref{eqn-main curv assump}). As noted earlier, the following estimate is roughly a continuous counterpart to the ``bush'' type arguments used in the combinatorial approach to bilinear restriction estimates. \\

Recall that given $\mathfrak{h} \in \RR^{1+n}$ we have defined the conical surfaces
$$\mc{C}_j(\mathfrak{h}) = \{ (s, -s\nabla \Phi_j(\xi) ) \mid s\in \RR, \xi \in \Sigma_j(\mathfrak{h})\}.$$
We define a weight $\chi^*_{j, r}(t,x)$ associated to $\mc{C}_j(\mathfrak{h})$ by
	$$ \chi^*_{j, r}(t,x) = \bigg(1 + \frac{\dist((t,x), \mc{C}_j(\mathfrak{h}))}{r} \bigg)^{-(n+2)}.$$
Thus $\chi^*_{j, r}$ is essentially concentrated on a $r$-neighbourhood of the conic surface $\mc{C}_j(\mathfrak{h})$. The geometric condition (\textbf{A1}) implies that the wave packets associated to the phase $\Phi_k$ intersect the surface $\mc{C}_j(\mathfrak{h})$ transversally. A rigorous version of this heuristic is the following.

\begin{theorem}\label{thm-null energy bound}
Let $\{j,k\} = \{1,2\}$, $\mathfrak{h} \in \RR^{1+n}$ and $ r \gg \mb{d}_0^{-1} \gtrsim \mc{H}_j $. Assume the phases $\Phi_j$ satisfy $(\mb{A1})$, $(\mb{A2})$, and the normalisation condition (\ref{eqn-assum on phase after rescaling}). If $v$ is a $\Phi_k$-wave with $\supp \widehat{v}$ contained in a ball of radius $r^{-1}$ and  $\supp \widehat{v} +  \frac{\mb{d}_0}{8} \subset \Lambda_k$ then we have
    $$ \big\| \chi^*_{j, r} v \big\|_{L^2_{t,x}}^2 \lesa r \| v \|_{L^\infty_t L^2_x}^2. $$
\end{theorem}
\begin{proof}
As in \cite{Tao2001b, Bejenaru2016b} by the $TT^*$ argument, the required bound is equivalent to proving that for every $F \in L^2_{t,x}$, we have
	$$ \int_\RR \int_\RR \lr{ \chi_r^*(t) \mc{U}(t) \mc{U}(-s)\big[ \chi_r^*(s) F(s)\big], F(t) }_{L^2_x} dt\,ds\, \lesa r \| F \|_{L^2_{t,x}}^2$$
where we let
	$$ \mc{U}(t)[f](x) = \int_{\RR^n} \rho(\xi) e^{i(t\Phi_k(\xi) + x \cdot \xi)} \widehat{f}(\xi) d\xi$$
with $\supp \widehat{\rho} \subset \Lambda_j$, $|\supp \widehat{\rho}|\lesa r^{-n}$, and $|\nabla^m \rho| \lesa r^m$. Define
	$$ K(t,x) = \int_{\RR^n} \rho^2(\xi) e^{ i(t\Phi_j(\xi) + x \cdot \xi)} d\xi.$$
By an application of Young's inequality, and using the assumption $|\supp \widehat{\rho}| \lesa r^{-n}$, it is enough to prove the pointwise bound
	\begin{equation}\label{eqn-thm null energy bound-key pointwise est} \big| \chi^*_r(t,x) K(t-s, x-y) \chi_r^*(s,y) \big|
			\lesa \int_{\supp \widehat{\rho}} \bigg( 1 +\frac{  |(t-s,x-y)|}{r} \bigg)^{-(n+2)}  d\xi.
	\end{equation}
To this end, take $p, p^* \in\mc{C}_j(\mathfrak{h})$ such that
	$$ |(t,x) - p| \approx \dist\big((t,x), \mc{C}_j(\mathfrak{h})\big), \qquad |(s,y) - p^*| \approx \dist\big((s,y), \mc{C}_j(\mathfrak{h})\big). $$
If $|(t,x) - p| + |(s,y) - p^*| \gtrsim |p-p^*|$ then (\ref{eqn-thm null energy bound-key pointwise est}) follows from the decay of $\chi^*_r$ away from the surface $\mc{C}_j(\mathfrak{h})$,  together with the the observation that
	$$ |(t-s, x-y)| \les |p-p^*|+ |(t,x) - p| + |(s,y) - p^*| \lesa |(t,x) - p| + |(s,y) - p^*|.$$
On the other hand, if $|(t,x) - p| + |(s,y) - p^*| \ll  |p-p^*|$, then an application of Lemma \ref{lem - surface C transverse} together with (\textbf{A1}) implies that for every $\xi \in \Lambda_k$
	\begin{align*} \big|(t-s, x-y)\wedge (1, - \nabla \Phi_k(\xi) \big) \big| &\g \big|(p-p^*)\wedge (1, - \nabla \Phi_k(\xi) \big) \big| - |(t,x) - p| - |(s,y) - p^*| \\
	&\gtrsim  |p-p^*| \\
	&\gtrsim  |(t-s, x-y)|
	\end{align*}
where the last inequality follows by writing $(t-s, x-y) = p-p^* + (t,x) - p + (s,y) - p^*$. In particular, we have
	\begin{equation}\label{eqn-thm null energy bound-phase nondegenerate}
		|x-y + (t-s) \nabla \Phi_k(\xi)| \approx |(t-s, x-y)|.
	\end{equation}
We now exploit the decay of the kernel $K(t,x)$ in directions orthogonal to $(1, - \nabla \Phi_k(\xi))$. More precisely, integrating by parts in the kernel $K(t,x)$ and applying (\textbf{A2}) implies that for every $N<5n$
    \begin{align*} |K&(t,x)| \\
    &\lesa \int_{\supp \widehat{\rho}} \frac{1}{|x+t\nabla \Phi_k(\xi)|} \sum_{\substack{ \ell \in \NN^{1+N} \\ \ell_0 + \ell_1 + 2\ell_2 + \dots + N\ell_N = N}} | \nabla^{\ell_0} \rho^2| \frac{ (t \|\nabla^2 \Phi_k\|)^{\ell_1}(t \|\nabla^3 \Phi_k\|)^{\ell_2} \dots (t \|\nabla^{N+1}\Phi_k\|)^{\ell_N}}{|x + t \nabla \Phi_k(\xi)|^{\ell_1 + \dots + \ell_N}} d\xi \\
                            &\lesa \int_{\supp \widehat{\rho}}  \bigg( \frac{r}{|x+t\nabla \Phi_k(\xi)|} \bigg)^N  \bigg(1 +  \frac{  |t| }{|x+t\nabla \Phi_k(\xi)|}\bigg)^{N} d\xi
    \end{align*}
where we applied (\textbf{A2}), (\ref{eqn-assum on phase after rescaling}), and the assumption $r\gtrsim \mb{d}_0^{-1} \gtrsim \mc{H}_j$. Therefore (\ref{eqn-thm null energy bound-key pointwise est}) follows from (\ref{eqn-thm null energy bound-phase nondegenerate}).
\end{proof}

Recall that, given a cube $q \in \mc{Q}_{r_j}(Q)$, we have defined $\chi_q \in C^\infty(\RR^{1+n})$ to be a positive function such that $\chi_q \gtrsim 1$ on $q$, $\supp \widetilde{\chi}_q \subset \{ |(\tau, \xi)| \les \frac{1}{r}\}$, and we have the decay bound, for every $N \in \NN$,
	$$ \chi_q(t,x) \lesa_N \Big(1 + \frac{  \dist((t,x), q)}{r} \Big)^{-N}. $$
To exploit the previous theorem, we require the following crucial consequence of the curvature of $\Phi_j$ on $\Sigma_j(\mathfrak{h})$.

\begin{lemma}\label{lem-sum concentrated on C}
Let $\mathfrak{h} \in \RR^{1+n}$ and $r_j \gg \mb{d}_0^{-1} \gtrsim \mc{H}_j$. Assume that $\Phi_j \in C^2(\Lambda_j)$ satisfies $(\mb{A1})$ and  the normalisation condition (\ref{eqn-assum on phase after rescaling}). Let $\{j,k\} = \{1,2\}$, $Q$ be a cube of diameter $R$ and $q_0 \in \mc{Q}_{r_j}(Q)$. Then provided that
    $$ \sigma_{\Sigma_j(\mathfrak{h})}\big( \Sigma_j(\mathfrak{h}) \cap \Lambda \big) \lesa \mb{d}_0^{n-1}$$
we have for all $(t,x) \in \RR^{1+n}$
  $$ \sum_{\substack{q \in \mc{Q}_{r_j}(B)\\|q-q_0| \g \epsilon R}} \sum_{\gamma_j \in \Gamma_j(\mathfrak{h}, \Lambda)} \frac{\chi_q(t,x)}{w_{\gamma_j, q_0} w_{\gamma_j, q} } \lesa \epsilon^{-C} \bigg(1 + \frac{\dist[(t-t_{q_0},x - x_{q_0}), \mc{C}_j(\mathfrak{h})]}{r_j} \bigg)^{-(n+2)}.$$
\end{lemma}
\begin{proof}
  We start by observing that, for every $(x_0, \xi_0) \in \Gamma_j(\mathfrak{h}, \Lambda)$, we have
    \begin{align*}
      \dist[(t-t_{q_0},x - x_{q_0})&, \mc{C}_j(\mathfrak{h})]\\
           &\lesa r_j +  |(t,x) - (t_q, x_q)| + | x_q - x_0 + t_q \nabla \Phi_j(\xi_0)| + |x_{q_0} - x_0 + t_{q_0} \nabla \Phi_j(\xi_0) |.
    \end{align*}
 In particular, as $\sum_{q} ( 1 + \frac{|(t,x) - (t_q, x_q)|}{r_j})^{n+2} \chi_q(t,x) \lesa 1$, it is enough to show that for fixed $q, q_0 \in \mc{Q}_{r_j}(Q)$ with $|q-q_0| \g \epsilon R$, we have
    $$ \sum_{(x_0, \xi_0) \in \Gamma_j(\mathfrak{h}, \Lambda)} \bigg(1 + \frac{|x_q - x_0 + t_q \nabla\Phi_j( \xi_0)| + |x_{q_0} - x_0 + t_{q_0} \nabla \Phi_j(\xi_0)| }{r_j}  \bigg)^{-3n} \lesa 1. $$
  For ease of reading, we let $\sigma(x_0, \xi_0) = |x_q - x_0 + t_q \nabla\Phi_j( \xi_0)| + |x_{q_0} - x_0 + t_{q_0} \nabla \Phi_j(\xi_0)|$, thus $\sigma$ is small when the tube $T_{(x_0, \xi_0)}$ passes through the cubes $q$ and $q_0$, and is large otherwise.  Let $(x_0, \xi_0) \in \Gamma_j(\mathfrak{h}, \Lambda)$ denote the phase space point with the associated tube passing closest to the cubes $q$ and $q'$, in other words, we take
  $$\sigma(x_0, \xi_0)=  \min_{(x_0', \xi_0') \in \Gamma_j(\mathfrak{h}, \Lambda)} \{ \sigma(x_0', \xi_0')\}.$$
Suppose for the moment that $\sigma(x_0, \xi_0)\ll \epsilon R$, thus the tube $T_{(x_0, \xi_0)}$ passes ``close'' to the cubes $q$ and $q_0$. In particular, writing
 $$(t_q - t_{q_0}) \big( 1, - \nabla \Phi_j(\xi_0) \big) =(t_q, x_q) - (t_{q_0}, x_{q_0}) - \big(0, x_q - x_0 + t_q \nabla \Phi_j(\xi_0)\big) +\big(0, x_{q_0} - x_0 + t_{q_0} \nabla \Phi_j(\xi_0)\big)$$
the separation of the cubes $q$ and $q_0$ implies that $|t_q - t_{q_0}| \approx \epsilon R$. Consequently, as
	$$ x_0 - x_0' = \big( x_q - x_0' + t_q \nabla \Phi_j(\xi_0')\big) - \big( x_q -x_0 + t_q \nabla \Phi_j(\xi_0) \big) + t_q \big( \nabla\Phi_j(\xi_0) - \nabla \Phi_j(\xi_0') \big), $$
and $|t_q| \lesa R$ we have by definition of $\sigma(x_0, \xi_0)$,
\begin{align*}
	|x_0 - x_0'| + R |\nabla \Phi_j(\xi_0) - \nabla \Phi_j(\xi_0')| &\les \sigma(x_0, \xi_0) + \sigma(x_0', \xi_0') + 2R |\nabla \Phi_j(\xi_0) - \nabla \Phi_j(\xi_0')| \\
	&\lesa \sigma(x_0, \xi_0) + \sigma(x_0', \xi_0') + \frac{1}{\epsilon}\big| (t_q - t_{q_0})\big( \nabla \Phi_j(\xi_0) - \nabla \Phi_j(\xi_0')\big| \\
	&\lesa \frac{1}{\epsilon}  \sigma(x_0', \xi_0').
\end{align*}
We now use the fact that $\xi_0, \xi_0' \in \Sigma_j(\mathfrak{h}) + \frac{C}{r_j}$, together with (\textbf{A1})  to deduce that
	$$ |\xi_0 - \xi_0'| \lesa \frac{|\nabla \Phi_j(\xi_0) - \nabla \Phi_j(\xi_0') |}{\mc{H}_j}  + \frac{1}{r_j}.$$
Therefore, since $\sigma(x_0, \xi_0) \les \sigma(x_0', \xi_0')$, we obtain
    \begin{align*}  \sum_{(x_0', \xi_0') \in \Gamma_j(\mathfrak{h}, \Lambda)} \bigg(1 + \frac{\sigma(x_0', \xi_0') }{r_j}  \bigg)^{-3n}
        &\lesa   \sum_{(x_0, \xi_0) \in \Gamma_j(\mathfrak{h}, \Lambda)}  \bigg(1 + \frac{\epsilon}{r_j}|x_0 - x_0'| + \epsilon r_j |\xi_0 - \xi_0'| \bigg)^{-3n} \\
        &\lesa \epsilon^{-C}.
    \end{align*}
On the other hand, if $\sigma(x_0, \xi_0) \gtrsim \epsilon R$, then we gain a power of $R$, and can simply discard the $\xi_0$ sum, and sum up over the points $x_0$. More precisely, the assumption on the set $\Lambda$, together with the condition $\mb{d}_0 \mc{H}_j \lesa 1$,  implies that
     $$ \#\Big\{ \xi_0 \in \tfrac{1}{r_j} \ZZ^n \,\, \Big| \,\, \xi_0 \in \Sigma_j(\mathfrak{h})\cap \Lambda + \tfrac{C}{r_j} \,\Big\} \lesa \frac{\sigma_{\Sigma_j}\big( \Sigma_j(\mathfrak{h})\cap \Lambda\big) r_j^{-1} + r_j^{-n}}{r_j^{-n}} \lesa ( 1 + \mb{d}_0 r_j)^{n-1} \lesa \Big( 1 + \frac{ R}{r_j} \Big)^{n-1} $$
and hence, as $\sigma(x_0, \xi_0) \gtrsim \epsilon R$, we obtain
    \begin{align*} \sum_{(x_0', \xi_0') \in \Gamma_j(\mathfrak{h}, \Lambda)} \bigg(1 + \frac{\sigma(x_0', \xi_0') }{r_j}  \bigg)^{-3n}
                \lesa \epsilon^{-C} \sup_{\xi_0'} \sum_{x_0'} \bigg(1 + \frac{|x_q - x_0' + t_q \nabla\Phi_j( \xi_0')|  }{r_j}  \bigg)^{-2n}
                \lesa \epsilon^{-C}.
    \end{align*}
\end{proof}

\section{The Wave Table Construction}\label{sec-wave table construct}

As implicitly used in the work of Wolff \cite{Wolff2001}, and more explicitly in the work of Tao \cite{Tao2001b}, we need to decompose $u$ depending on $v$. Roughly speaking we will decompose $u$ into wave packets, and keep the wave packets where $v$ and $\mc{P}_{\gamma_j} u$ simultaneously concentrate on the same cube $B \in \mc{Q}_{\frac{R}{4}}(Q)$. We adapt the following construction from Tao, but change notation slightly.

\begin{definition}[$\Phi_j$-Wave Tables]\label{def:wave tables}
Let $0<\epsilon<1$. Let $0<\| F \|_{L^\infty_t L^2_x} < \infty$, and $Q$ be a cube of diameter $R \gg (\mb{d}_0^2 \mc{H}_j)^{-1}$. Let $u$ be a $\Phi_j$-wave. The $\Phi_j$-wave table of $u$, relative to $F$ and $Q$, is the collection of (vector valued) functions $(\mc{W}_{j,\epsilon}^{(B)})_{B \in \mc{Q}_{\frac{R}{4}}(Q)}$ where
	$$ \mc{W}_{j,\epsilon}^{(B)}=\mc{W}_{j, \epsilon}^{(B)}(u;F, Q) = \sum_{\gamma_j \in \Gamma_{j}} \frac{ m_{\gamma_j, B}(F)}{m_{\gamma_j}(F)} \mc{P}_{\gamma_j} u $$
where the coefficients are defined as
  $$ m_{\gamma_j, B}(F) = \sum_{q \in \mc{Q}_{r_j}(B) } \frac{1}{w_{\gamma_j,q}}\| \chi_q F \|_{L^2_{t,x}(\RR^{1+n})}^2, \qquad m_{\gamma_j}(F) = \sum_{B \in \mc{Q}_{\frac{R}{4}}(Q)} m_{\gamma_1, B}(F).$$
\end{definition}

\begin{remark} The assumption $0<\| F \|_{L^\infty_t L^2_x}< \infty$ implies that\footnote{Since $\chi_q$ has compact Fourier support, it can only have countable number of zeros. In particular, $\chi_q>0$ almost everywhere.} $0<\| \chi_q F \|_{L^2_{t,x}}<\infty$, and hence the coefficients $m_{\gamma_j, B}(F)$ are well-defined and satisfy $0<m_{\gamma_j, B}(F)<\infty$. Note that we essentially have $m_{\gamma_j, B}(F) \approx \| F \|_{L^2_{t,x}(T_{\gamma_j} \cap B)}^2$ and $m_{\gamma_j}(F) \approx \| F \|_{L^2_{t,x}( T_{\gamma_j}\cap Q)}^2$. Thus $\mc{W}^{(B)}_{j, \epsilon}$ contains all wave packets $\mc{P}_{\gamma_j} u$ such that $ F|_{T_{\gamma_j}}$ is concentrated on $B$.
\end{remark}

Recall that if $Q$ is a cube of diameter $R$, $0<r\les R$, and $(u^{(B)})_{B\in \mc{Q}_r(Q)}$ is a collection of ($\ell^2_c(\ZZ)$ valued) functions, we defined the quilt $[u^{(\cdot)}]$ as
	$$ [u^{(\cdot)}] = \sum_{B \in \mc{Q}_{r}(Q)}  \ind_{B} |u^{(B)}|$$
and the $\epsilon$ separated cubes
	$$ I^{\epsilon, r}(Q) = \bigcup_{q \in \mc{Q}_r(Q)} (1-\epsilon) q.$$
The definition of the wave tables $\mc{W}^{(B)}_{j, \epsilon}$ gives the following key bilinear estimate.

\begin{theorem}\label{thm-prop of wave tables}
Let $ \{j, k\} = \{1,2\}$, $0<\epsilon \ll 1$,  and $R \gg ( \mb{d}_0^2 \mc{H}_j)^{-1}$. Assume that the phases $\Phi_1$ and $\Phi_2$ satisfy Assumption \ref{assump-main}. Let $u$ be a $\Phi_j$-wave and $v$ be a $\Phi_k$-wave with
    $$\supp \widehat{u} + \tfrac{\mb{d}_0}{2} \subset \Lambda_j, \qquad \supp \widehat{v} + \tfrac{\mb{d}_0}{2} \subset \Lambda_k, \qquad \mb{d}[\supp\widehat{u}, \supp \widehat{v}]\lesa \mb{d}_0.$$
Let $Q$ be a cube of diameter $R$, and for $B \in \mc{Q}_{\frac{R}{4}}(Q)$, take $\mc{W}^{(B)}_{j, \epsilon} = \mc{W}^{(B)}_{j, \epsilon}(u;v, Q)$. Then $\mc{W}^{(B)}_{j, \epsilon}$ is again a $\Phi_j$-wave, is linear in $u$, satisfies the support condition
	$$\supp \widehat{\mc{W}}^{(B)}_{j, \epsilon} \subset \supp \widehat{u} +  \tfrac{1}{r_j},$$
and we have the energy estimate
	$$ \bigg(\sum_{B \in \mc{Q}_{\frac{R}{4}}(Q)} \big\| \mc{W}^{(B)}_{j, \epsilon} \big\|_{L^\infty_t L^2_x}^2\bigg)^\frac{1}{2} \les (1 +  C\epsilon) \| u \|_{L^\infty_t L^2_x}. $$
Moreover, we have the bilinear estimate
     $$ \Big\| \big( |u| - \big[\mc{W}^{(\cdot)}_{j, \epsilon}\big] \big) v \Big\|_{L^2_{t,x}(I^{\epsilon, \frac{R}{4}}(Q))}\\
        \lesa \epsilon^{-C} r_j^{-\frac{n-1}{2}} \| u\|_{L^\infty_t L^2_x} \|v \|_{L^\infty_t L^2_x}. $$
\end{theorem}

\begin{proof}  For ease of notation, for the remainder of the proof, we let $u^{(B)}=\mc{W}^{(B)}_{j, \epsilon}(u; v, Q)$. By construction, we have $ \sum_{B \in \mc{Q}_{\frac{R}{4}}(Q)} \frac{m_{\gamma_j, B}(v)}{m_{\gamma_j}(v)} = 1$ and hence
          $u = \sum_{B \in \mc{Q}_{\frac{R}{4}}(Q)} u^{(B)}$.
The initial claims follow immediately from the definition of the wave table $(\mc{W}^{(B)}_{j, \epsilon})_{B \in \mc{Q}_{\frac{R}{4}}(Q)}$, together with Lemma \ref{lem-orthog prop of wave packets} and Proposition \ref{prop - concentration of wave packets II} (note that, since $u$ is a $\Phi_j$-wave, we have $\| u \|_{L^\infty_t L^2_x} = \|u(0)\|_{L^2_x}$). We now turn to the proof of the bilinear estimate. The identity $u=\sum_{B \in \mc{Q}_{\frac{R}{4}}(Q)} u^{(B)}$ implies that
    $$ \big| |u| \ind_Q - [u^{(\cdot)}]\big| \les \sum_{B \in \mc{Q}_{\frac{R}{4}}(Q)} |u^{(B)}| \ind_{Q \setminus B}.$$
Consequently, as there are only $4^{n+1}$ cubes in $\mc{Q}_{\frac{R}{4}}(Q)$, the separation of cubes inside $I^{\epsilon, \frac{R}{4}}$ implies that
    $$  \Big\| \big( |u| - \big[u^{(\cdot)}\big]\big)v \Big\|_{L^2_{t,x}(I^{\epsilon, \frac{R}{4}}(Q))} \lesa \sup_{B \in \mc{Q}_{\frac{R}{4}}(Q)} \big\| u^{(B)} v  \big\|_{L^2_{t,x}(I^{\epsilon, \frac{R}{4}}(Q)\setminus B )}\lesa \sum_{\substack{q\in \mc{Q}_{r_j}(Q) \\ |q-B| \g \epsilon R}} \| u^{(B)}  v \|_{L^2_{t,x}(q)}^2. $$
For $\mathfrak{h}=(a,h) \in \RR^{1+n}$, we let $H_{\mathfrak{h}}$ be a Fourier projection onto the set $\{ |\tau - a| \les \frac{1}{r_j}, |\xi - h |\les \frac{1}{r_j} \}$ such that
        $$ \| F \|_{L^2_{t,x}}^2 \approx \sum_{\mathfrak{h} \in \frac{1}{r_j}\ZZ^{1+n}}\| H_{\mathfrak{h}} F\|_{L^2_{t,x}}^2.$$
Observe that given $\gamma_j = (x_0, \xi_0) \in \Gamma_j$ and $\mathfrak{h}=(a,h) \in \RR^{1+n}$, we only have $H_{\mathfrak{h}}(\chi_q \mc{P}_{\gamma_j} u \chi_q v ) \not= 0$ if there exists $\xi \in \Lambda_j$, $\eta \in \Lambda_k$ such that
	$$ |\xi_0 - \xi| \lesa \frac{1}{r_j}, \qquad |\xi + \eta - h|\lesa \frac{1}{r_j}, \qquad |\Phi_j(\xi) + \Phi_k(\eta) - a| \lesa \frac{1}{r_j}.$$
Therefore, an application of Lemma \ref{lem - folation properties of surface Sigma} implies that if $H_{\mathfrak{h}}(\chi_q \mc{P}_{\gamma_j} u \chi_q v ) \not= 0$  we must have $\xi_0 \in \Sigma_j(\mathfrak{h})\cap \Lambda + \frac{ C}{r_j}$ with $\Lambda = \supp \widehat{u} \cap ( h - \supp \widehat{v})$. Consequently we have
	\begin{align*} \| u^{(B)} v \|_{L^2_{t,x}(q)}^2 &\lesa \| \chi_q u^{(B)} \chi_q v \|_{L^2_{t,x}}^2 \\
				&\lesa \sum_{\mathfrak{h} \in \frac{1}{r_j} \ZZ^{1+n}}  \Big\| \sum_{\gamma_j \in \Gamma_j(\mathfrak{h}, \Lambda)} \frac{m_{\gamma_j, B}}{m_{\gamma_j}} H_{\mathfrak{h}}\big( \chi_q \mc{P}_{\gamma_j} u \chi_q v \big)\Big\|_{L^2_{t, x}}^2.
	\end{align*}
Since $\frac{m_{\gamma_j,B}(v)}{m_{\gamma_j}(v)} \les (\frac{m_{\gamma_j, B}(v)}{m_{\gamma_j}(v)})^\frac{1}{2} $, an application of H\"older's inequality gives
    \begin{align*}
     \Big\| \sum_{\gamma_j \in \Gamma_j(\mathfrak{h}, \Lambda)} \frac{m_{\gamma_j, B}(v)}{m_{\gamma_j}(v)} H_{\mathfrak{h}}\big( \chi_q \mc{P}_{\gamma_j} u \chi_q v \big)\Big\|_{L^2_{t, x}}^2
   &\lesa   \sum_{\gamma_j \in \Gamma_j} \frac{w_{\gamma_j, q}}{m_{\gamma_j}(v)} \Big\|   H_{\mathfrak{h}}\big( \chi_q \mc{P}_{\gamma_1} u \chi_q v \big)\Big\|_{L^2_{t, x}}^2 \sup_{\mathfrak{h} \in \RR^{1+n}} \sum_{\gamma_j \in \Gamma_j(\mathfrak{h}, \Lambda)} \frac{m_{\gamma_j, B}(v)}{w_{\gamma_j, q}}.
   \end{align*}
 Hence summing up over $\mathfrak{h} \in \frac{1}{r_j} \ZZ^{1+n}$, and noting that the product $\chi_q \mc{P}_{\gamma_j} u$ has Fourier support in a set of size $r_j^{-(n+1)},$ we deduce that
	\begin{align} &\sum_{\substack{q\in \mc{Q}_{r_j} \\ |q-B| \g \epsilon R}} \| u^{(B)} v \|_{L^2_{t,x}(q)}^2\notag\\
		&\lesa \bigg(\sum_{\substack{q\in \mc{Q}_{r_j}(Q) }} \sum_{\gamma_j \in \Gamma_j} \frac{w_{\gamma_j, q}}{m_{\gamma_j}(v)}  \| \chi_q \mc{P}_{\gamma_j} u \chi_q v \|_{L^2_{t,x}}^2 \bigg) \times \bigg( \sup_{\substack{ \mathfrak{h} \in \RR^{1+n} \\ q \in \mc{Q}_{r_j}(Q),\\ |q-B|\g \epsilon R}} \sum_{\gamma_j \in \Gamma_j(\mathfrak{h}, \Lambda)} \frac{m_{\gamma_j, B}(v)}{w_{\gamma_j, q}} \bigg)\notag\\
        &\lesa r_j^{-(n+1)} \bigg(\sum_{\substack{q\in \mc{Q}_{r_j}(Q) }} \sum_{\gamma_j \in \Gamma_j} \frac{w_{\gamma_j, q}}{m_{\gamma_j}(v)}  \| \chi_q \mc{P}_{\gamma_j} u\|_{L^2_{t,x}}^2 \| \chi_q v \|_{L^2_{t,x}}^2 \bigg) \times \bigg( \sup_{\substack{ \mathfrak{h} \in \RR^{1+n} \\ q \in \mc{Q}_{r_j}(Q),\\ |q-B|\g \epsilon R}} \sum_{\gamma_j \in \Gamma_j(\mathfrak{h}, \Lambda)} \frac{m_{\gamma_j, B}(v)}{w_{\gamma_j, q}} \bigg).\label{eqn-thm prop of wave tables-general wave table coefficients}
        \end{align}
We now estimate the first term in (\ref{eqn-thm prop of wave tables-general wave table coefficients}). The definition of the coefficients $m_{\gamma_j}(v)$ together with Proposition \ref{prop - concentration of wave packets II} gives
    \begin{align*}
      \sum_{\substack{q\in \mc{Q}_{r_j}(Q)}} &\sum_{\gamma_j \in \Gamma_j} \frac{w_{\gamma_j, q}}{m_{\gamma_j}(v)}  \| \chi_q \mc{P}_{\gamma_j} u\|_{L^2_{t,x}}^2 \| \chi_q v \|_{L^2_{t,x}}^2\\
            &\lesa  \bigg(\sum_{\gamma_j \in \Gamma_j} \sup_{q \in \mc{Q}_{r_j}(Q)} w_{\gamma_j, q}^2  \| \chi_q \mc{P}_{\gamma_j} u\|_{L^2_{t,x}}^2\bigg) \times \bigg( \sup_{\gamma_j \in \Gamma_j} \frac{1}{m_{\gamma_j}(v)} \sum_{q\in \mc{Q}_{r_j}(Q)} \frac{1}{w_{\gamma_j, q}} \| \chi_q v \|_{L^2_{t,x}}^2\bigg)\\
            &\lesa r_j \epsilon^{-C} \| u \|_{L^\infty_t L^2_x}^2.
    \end{align*}
 On the other hand, to estimate the second term in \eref{eqn-thm prop of wave tables-general wave table coefficients}, expanding the sum over the $m_{\gamma_j, B}(v)$, we have
  $$ \sum_{\gamma_j \in \Gamma_j(\mathfrak{h}, \Lambda)} \frac{m_{\gamma_j, B}(v)}{w_{\gamma_j, q}} = \sum_{q_0 \in \mc{Q}_{r_j}(B)} \sum_{\gamma_j \in \Gamma_j(\mathfrak{h}, \Lambda)} \frac{1}{w_{\gamma_j, q} w_{\gamma_j,q_0}}\| \chi_{q_0} v \|_{L^2_{t,x}(\RR^{1+n})}^2.$$
Since $\chi_{q_0}$ has Fourier support in a ball of radius $r_j^{-1}$, by orthogonality, we may assume that $\widehat{v}$ is also supported in a ball of radius $r_j^{-1}$. Consequently, an application of Lemma \ref{lem-sum concentrated on C}, translation invariance,  and Theorem \ref{thm-null energy bound} implies that
	$$ \sum_{\gamma_j \in \Gamma_j(\mathfrak{h}, \Lambda)} \frac{m_{\gamma_j, B}(v)}{w_{\gamma_j, q}} \lesa r_j \|v\|_{L^\infty_tL^2_x}. $$
Therefore the required bilinear estimate follows.
\end{proof}

\section{Proof of Theorem \ref{thm-general wave table decomposition}} \label{sec-proof of key bilinear est}

We are now ready to give the proof of Theorem \ref{thm-general wave table decomposition}. Roughly the idea is that, given $R \gg (\mb{d}_0^2 \mc{H}_j)^{-1}$,  Theorem \ref{thm-prop of wave tables} allows us to essentially replace the (low frequency) term $u$ with pieces concentrated on $\frac{R}{4}$ cubes. However, since Theorem \ref{thm-prop of wave tables} applied to $u$ only requires scales $\gg \mb{d}_0^{-2}$ (as $\mc{H}_1=1$) we can continue decomposing $u$ until we get terms concentrated on the much smaller $\mc{H}_2 R$ cubes. We then decompose $v$, but as Theorem \ref{thm-prop of wave tables} applied to $v$ requires $R \gg (\mb{d}_0^2 \mc{H}_j)^{-1}$, we can only apply it once before we start to lose logarithmic factors. The immediate obstruction to applying the above strategy, is that Theorem \ref{thm-prop of wave tables} only applies to the separated cubes $I^{\epsilon, r}(Q)$. To create the required separation, we use an averaging argument due to Tao. \\

Suppose $Q$ is a cube of radius $R$, and let $(\epsilon_m)$ and $(r_m)$ be strictly positive sequences with $r_m\les R$. We then define
    $$X[Q]= \cap_{m=1}^M I^{\epsilon_m, r_m}(Q)$$
where we recall that $I^{\epsilon, r}(Q)= \cup_{q \in \mc{Q}_r(Q)} (1-\epsilon)q$. Thus cubes inside $X[Q]$ are separated at multiple scales (which is needed as we apply Theorem \ref{thm-prop of wave tables} at multiple scales). We will need to move from integrating over $Q$, to integrating over the smaller $X[Q]$. The key tool is the following averaging lemma due to Tao.

\begin{lemma}[{\cite[Lemma 4.1]{Lee2008}, \cite[Lemma 6.1]{Tao2001b}}]\label{lem - cube averaging}
Let $1\les q, r \les \infty$, $R>0$ and let $Q_R$ be a cube of diameter $R$. If $\epsilon = \sum_{m=1}^M \epsilon_j \les 2^{-(n+2)}$ then there exists a cube $Q \subset 4 Q_R$ of side lengths $2R$ such that
    \begin{equation}\label{eqn - lem cube averaging - main ineq} \| F \|_{L^q_t L^r_x(Q_R)} \les \Big(1 + 2^{n+2} \epsilon \Big) \| F \|_{L^q_t L^r_x(Q_R \cap X[Q])}.  \end{equation}
\begin{proof} We start by considering the case $q=r=1$. Let $Q(t,x)$ denote the cube of side lengths $2R$ centered at the point $(t,x) \in Q_R$. It is enough to prove that
    \begin{equation}\label{eqn-lem cube avg-pigeon hole version}
        \| F \|_{L^1(Q_R)} \les (1+ 2^{n+2} \epsilon)  \frac{1}{|Q_R|} \int_{Q_R} \| F \|_{L^1(Q_R \cap X[Q(t,x)])} \,dt\,dx.
    \end{equation}
To this end, we begin by observing that since\footnote{Just use the identities $X[Q(t,x)] = X[Q(0)] + (t,x)$ and $-X[Q(0)] = X[-Q(0)] = X[Q(0)]$. }
        $$ (t',x') \in X[ Q(t,x)] \Longleftrightarrow (t,x) \in X[ Q(t',x')]$$
we have for every $(t',x') \in Q_R$
    $$\int_{Q_R} \ind_{X[Q(t,x)]}(t',x') \,dt \,dx = |\{ (t,x) \in Q_R \mid (t',x') \in X[Q(t,x)]\}| =   | Q_R \cap X[Q(t',x')]|. $$
Thus by an application of Fubini we have
        \begin{align*}
          \int_{Q_R} \| F \|_{L^1(Q_R \cap X[Q(t,x)])} \,dt\,dx &= \int_{Q_R} \int_{Q_R} \ind_{X[Q(t, x)]}(t',x') \,dt\,dx |F(t',x')| \,dt'\,dx' \\
                        &= \int_{Q_R} |Q_R \cap X[Q(t',x')]|\, |F(t',x')| \,dt' \,dx'
        \end{align*}
Consequently the required bound (\ref{eqn-lem cube avg-pigeon hole version}) would follow from the inequality
    \begin{equation}\label{eqn-lem cube avg-bound on Q_R}
      |Q_R| \les (1 + 2^{n+2} \epsilon) \big|Q_R \cap X[ Q(t,x) ]\big|.
    \end{equation}
Let $X_k= \cap_{m=k}^M I^{\epsilon_j, r_j}(Q(t,x))$. Noting the inclusion $Q_R \subset Q(t,x)$, we deduce that
    \begin{align*}|Q_R \cap X_1 | &= |Q_R \cap I^{\epsilon_1, r_1}(Q(t,x)) \cap X_{2}| \\
    &= \sum_{q \in \mc{Q}_{r_1}(Q(t,x))} | (1-\epsilon_1) q\cap Q_R \cap X_2|\\
     &\g \sum_{q \in \mc{Q}_{r_1}(Q(t,x))} | q\cap Q_R \cap X_2|  - \sum_{q\in\mc{Q}_{r_1}(Q(t,x))} |q \setminus (1-\epsilon_1)q| \\
     &= |Q_R\cap X_2|  - \big( 1 - (1-\epsilon_1)^n \big) (2R)^n.
    \end{align*}
Since $1 - (1-\epsilon_1)^n \les n \epsilon_1$, repeating the above argument eventually gives
    $$ |Q_R \cap X[Q(t,x)]| \g |Q_R| - n 2^n R^n \sum_{m=1}^M \epsilon_m \g (1 -  2^{n+1} \epsilon ) |Q_R|$$
and so, provided that $\epsilon\les 2^{-(n+2)}$ we obtain (\ref{eqn-lem cube avg-bound on Q_R}). Thus the case $q=r=1$ follows. The general case follows by observing that there exists $G \in L^{q'}_t L^{r'}_x(Q_R)$ with $\| G \|_{L^{q'}_t L^{r'}_x(Q_R)} \les 1$ such that
            $$ \| F \|_{L^q_r L^r_x(Q_R)} \les \| F G \|_{L^1_{t,x}(Q_R)}$$
together with an application of the $L^1_{t,x}$ case obtained above.
\end{proof}
\end{lemma}

We now come to the proof of Theorem \ref{thm-general wave table decomposition}.

\begin{proof} It is enough to consider the case $\| u \|_{L^\infty_t L^2_x} = \| v \|_{L^\infty_t L^2_x} = 1$. Fix a cube $Q_R$ of diameter $R \gg (\mb{d}_0^2 \mc{H}_2)^{-1}$ and let $M\in \NN$ such that
		$$ 4^{-M} < \mc{H}_2 \les 4^{1-M}.$$
An application of Lemma \ref{lem - cube averaging} implies that there exists a cube $Q$ of radius $2R$ such that
        $$ \| uv \|_{L^q_t L^r_x(Q_R)} \les ( 1  + C \epsilon) \| uv \|_{L^q_t L^r_x(X[Q])}$$
where we take
	$$ X[Q] = \bigcap_{m=1,\dots, M} I^{\epsilon_m, 4^{-m} 2R}(Q), \qquad \epsilon_m = 4^{\delta(m-M)} \epsilon$$
and $\delta>0$ is some fixed constant to be chosen later (which will depend only on the dimension $n$, and the constant $C$ appearing in Theorem \ref{thm-prop of wave tables}). We start by decomposing $u$. Given $B_1 \in \mc{Q}_{\frac{R}{2}}(Q)$ we define the $\Phi_1$-wave
	$$ u^{(B_1)}_1 = \mc{W}^{(B_1)}_{1, \epsilon_1}(u; v, Q)$$
and assuming we have constructed $u_{m}^{(B_{m})}$, $B_{m} \in \mc{Q}_{\frac{2R}{4^m}}(Q)$, we define for $B_{m+1} \in \mc{Q}_{\frac{2R}{4^{m+1}}}(B_m)$
		$$ u_{m+1}^{(B_{m+1})} = \mc{W}^{(B_{m+1})}_{1, \epsilon_{m+1}}\big( u^{(B_{m})}_{m}; v, B_{m}\big)$$
and finally take $u^{(B)} = u_M^{(B)}$ for $B \in \mc{Q}_{\frac{2R}{4^M}}(Q) $. Note that each $u_m^{(B_m)}$ is well-defined, since the diameter of the cube $B_{m-1}$ is $4^{-m} 2R \g 4^{-M} 2R \gg \mb{d}_0^2$. To construct $v$, we start by defining $U=(u^{(B)})_{B\in\mc{Q}_{\frac{2R}{4^{M}}}(Q)}$. Clearly, as each $u^{(B)}$ is a $\Phi_1$-wave, after relabeling, $U$ is also a $\Phi_1$-wave. We now decompose $v$ relative to $U$ and the cube $Q$, in other words for $B' \in \mc{Q}_{\frac{R}{2}}(Q)$ we define
		$$ v^{(B')} = \mc{W}^{(B')}_{2, \epsilon}\big( v; U, Q\big).$$
It is easy to check that, by construction, Theorem \ref{thm-prop of wave tables} implies that we have the identities
	$$ u = \sum_{B \in \mc{Q}_{\frac{2R}{4^M}}(Q)} u^{(B)}, \qquad \qquad v = \sum_{B' \in \mc{Q}_{\frac{R}{2}}(Q)} v^{(B')}$$
together with the support conditions
	$$ \supp \widehat{u}_m^{B_m} \subset \supp \widehat{u}_{m-1}^{B_{m-1}} + \Big( \frac{2R}{4^{m-1}} \Big)^{-\frac{1}{2}}, \qquad \supp \widehat{v}^{(B')} \subset \supp \widehat{v} + \big( 2 \mc{H}_2 R\big)^{-\frac{1}{2}}.$$
Moreover, $v^{(B')}$ satisfies the required energy estimate, and by iterating the estimate
	\begin{align*} \sum_{B_m \in \mc{Q}_{\frac{2R}{4^m}}(Q)} \|u^{(B_m)}_m\|_{L^\infty_t L^2_x}^2 &= \sum_{B_{m-1} \in \mc{Q}_{\frac{2R}{4^{m-1}}}(Q)} \sum_{B_m \in \mc{Q}_{\frac{2R}{4^{m}}}(B_{m-1})} \|u^{(B_m)}_m\|_{L^\infty_t L^2_x}^2 \\
	&\les (1 + \epsilon_m C)^2 \sum_{B_{m-1} \in \mc{Q}_{\frac{2R}{4^{m-1}}}(Q)} \| u^{(B_{m-1})}_{m-1}\|_{L^\infty_t L^2_x}^2.
	\end{align*}
and using the elementary inequality $ \Pi_{m=1}^M (1 + \epsilon_m C) \les 1 + e^{ \sum \epsilon_m C} C \sum \epsilon_m $ and the definition of $\epsilon_m$, a short computation shows that $u^{(B)}$ satisfies the correct energy bound.

We now turn to the proof of the bilinear estimate. Writing
    $$ |u v| = [u^{(\cdot)}] [v^{(\cdot)}] + \big( |u| - [u^{(\cdot)}\big) |v| + [u^{(\cdot)}] \big( |v| - [v^{(\cdot)}]\big),$$
after an application of Holder's inequality in the $t$ variable, it is enough to show that for every $2 \g r >\frac{n+1}{n}$ we have
    \begin{equation}\label{eqn-thm main bilinear est-temp bilinear est}
        \big\| \big( |u| - [u^{(\cdot)}]\big) v \big\|_{L^2_t L^r_x(X[Q])} + \big\| [u^{(\cdot)}] \big( |v| - [v^{(\cdot)}]\big) \big\|_{L^2_t L^r_x(X[Q])}
                \lesa \epsilon^{-C} \big(\mc{H}_2 R\big)^{\frac{n+1}{2r} - \frac{n}{2}}
    \end{equation}
We start by estimating the first term. An application of Theorem \ref{thm-prop of wave tables} gives
	\begin{align*}
		\Big\| \Big( \big[ u^{(\cdot)}_{m-1}\big] -& \big[ u^{(\cdot)}_{m} \big] \Big) v \Big\|_{L^2_{t,x}(X[Q])}^2   \\
		&\les \sum_{B_{m-1} \in \mc{Q}_{\frac{2R}{4^{m-1}}}(Q)} \Big\| \Big( |u^{(B_{m-1})}_{m-1}| - \big[ \mc{W}^{(\cdot)}_{1, \epsilon_m}(u^{(B_{m-1})}_{m-1})\big] \Big) v \Big\|_{L^2_{t,x}(I^{\epsilon_m, \frac{2R}{4^{m}}}(B_{m-1}))}^2\\
		&\lesa \epsilon_m^{-2C} \Big( \frac{4^{m-1}}{2R}\Big)^{\frac{n-1}{2}}  \sum_{B_{m-1} \in \mc{Q}_{\frac{2R}{4^{m-1}}}(Q)} \| u^{(B_{m-1})}_{m-1} \|_{L^\infty_t L^2_x}^2 \\
		&\lesa \epsilon_m^{-2C} \big( 4^{-m} R \big)^{-\frac{n-1}{2}}
	\end{align*}
where we used the definition of $M$. On the other hand, an application of Holder's inequality together with the energy bound for $u^{(\cdot)}_m$ and $u_{m-1}^{(\cdot)}$ implies that
	\begin{align*}
	\Big\| \Big( \big[ u^{(\cdot)}_{m-1}\big] -\big[ u^{(\cdot)}_{m} \big]\Big) v \Big\|_{L^2_t L^1_x(X[Q])}
	&\les \Big\| \Big( \big[ u^{(\cdot)}_{m-1}\big]- \big[ u^{(\cdot)}_{m} \big] \Big)\Big\|_{L^2_{t,x}(Q)}  \| v\|_{L^\infty_t L^2_x(Q)}\\
	&\lesa \Bigg(\sum_{B_{m-1}} \| u^{(B_{m-1})}_{m-1}\|_{L^2_{t,x}(B_{m-1})}^2 + \sum_{B_m} \| u^{(B_m)}_m \|_{L^2_{t,x}(B_m)}^2 \Bigg)^\frac{1}{2}  \\
	&\lesa  \big( 4^{-m} R\big)^{\frac{1}{2}}
	\end{align*}
Hence interpolating gives for any $1\les r \les 2$
	$$ \Big\| \Big( \big[ u^{(\cdot)}_{m-1}\big] - \big[ u^{(\cdot)}_{m} \big] \Big) v \Big\|_{L^2_t L^r_x(X[Q])}
			\lesa \epsilon_m^{-2 C (1-\frac{1}{r})} \big( 4^{-m} R \big)^{\frac{n+1}{2r} - \frac{n}{2}} .$$
The definition of $\epsilon_m$ and $M$ implies that we can write
		$$ \epsilon_m^{-C_n(1-\frac{1}{r})} 4^{\frac{m}{2}(n-\frac{n+1}{r})} \les \epsilon^{-C} \mc{H}_2^{ \frac{n+1}{2r} - \frac{n}{2}} 4^{\delta^* (m-M)}$$
 where $\delta^* = \frac{n}{2} - \frac{n+1}{2r} - C\delta ( 1 - \frac{1}{r})>0$ provided $r>\frac{n+1}{n}$ and we choose $\delta$ sufficiently small. Therefore, telescoping the sum over $m$ and letting $ u_0^{(Q)} = u$, we deduce that
	\begin{align*} \big\| \big( |u| - [u^{(\cdot)}]\big) v \big\|_{L^2_t L^r_x(X[Q])} &\les \sum_{m=1}^M \big\| \big( [ u^{(\cdot)}_{m-1}] - [ u^{(\cdot)}_{m} ] \big) v \big\|_{L^2_t L^r_x(X[Q])} \lesa \epsilon^{-C} \big( \mc{H}_2 R \big)^{\frac{n+1}{2r}- \frac{n}{2} }.
	\end{align*}
Thus it only remains to estimate the second term in (\ref{eqn-thm main bilinear est-temp bilinear est}). To this end, applying the definition of $v^{(B')}$ together with Theorem \ref{thm-prop of wave tables}, we have
		\begin{align*} \big\| [u^{(\cdot)}] \big( |v| - [v^{(\cdot)}] \big)\big\|_{L^2_{t,x} (X[Q])} &\les\big\| U \big( |v| - [v^{(\cdot)}] \big)\big\|_{L^2_{t,x} (I^{\epsilon_1, \frac{R}{2}})} \lesa  \epsilon^{-C}  \big( \mc{H}_2 R \big)^{-\frac{n-1}{4}}
        \end{align*}
where we used the inequality $[u^{(\cdot)}]\les (\sum_B |u^{(B)}|^2)^\frac{1}{2}=  |U|$. On the other hand an application of Holder's inequality together with the energy estimates gives
	 \begin{align*}  \big\| [u^{(\cdot)}] \big( |v| - [v^{(\cdot)}] \big)\big\|_{L^2_t L^1_x (X[Q])} &\les \bigg( \sum_{B \in \mc{Q}_{\frac{2R}{4^M}}} \big\| u^{(B)} \|_{L^2_{t,x}(B)}^2 \bigg)^\frac{1}{2} \Big( \|v \|_{L^\infty_t L^2_x} + \| [v^{(\cdot)}] \|_{L^\infty_t L^2_x}\Big) \lesa  \big( \mc{H}_2 R\big)^{\frac{1}{2}}.
	 \end{align*}
Therefore, interpolating between these bounds gives
	$$ \big\| [u^{(\cdot)}] \big( |v| - [v^{(\cdot)}] \big)\big\|_{L^2_t L^r_x (X[Q])} \lesa \epsilon^{-C} \big(\mc{H}_2 R \big)^{\frac{n+1}{2r} - \frac{n}{2}}$$
and hence (\ref{eqn-thm main bilinear est-temp bilinear est}) follows.
\end{proof}

\section{Atomic Wave Tables and the proof of Theorem \ref{thm-main Up version II}}\label{sec:atomic wave tables}

In this section we give the proof of Theorem \ref{thm-main Up version II}, namely the extension of Corollary \ref{cor-main Up version} to the atomic spaces $U^a_{\Phi_1}$ and $U^b_{\Phi_2}$ with $a, b \g 2$. As the proof closely follows the argument used to prove Theorem \ref{thm-main small scale}, we will be somewhat brief. \\

 We start by defining an \emph{atomic $\Phi_j$-wave} to be a function of the form
        $$ u(t, x) =\sum_{I \in \mc{I}} \ind_I(t) u_I(t,x)$$
where $\mc{I}$ is a finite partition of $\RR$, $\ind_I(t)$ is the indicator function adapted to the interval $I \subset \RR$, and $u_I$ is a $\Phi_j$-wave. Given an atomic $\Phi_j$-wave $u$, we use the shorthand
    $$ \| u \|_{\ell^a L^2_x} = \bigg( \sum_{I\in \mc{I}} \| u_I \|_{L^\infty_t L^2_x}^a \bigg)^\frac{1}{a}.$$
In particular, if $\| u \|_{\ell^a L^2_x} \les 1$, then $u$ is simply a $U^a_{\Phi_j}$ atom. Note that an atomic $\Phi_j$-wave is essentially a special case of a $\Phi_j$-wave, since if we let $U=(u_I)_{I\in \mc{I}}$ then perhaps after relabeling, $U$ is clearly a $\Phi_j$-wave, and moreover, as in the proof of Corollary \ref{cor-main Up version} we have the pointwise bound
    $$ |u(t,x)| \les \bigg( \sum_{I\in \mc{I}} |u_I(t,x)|^2 \bigg)^\frac{1}{2}=|U(t,x)|.$$
However, by exploiting the time localisation, atomic $\Phi_j$-waves satisfy slightly stronger bilinear estimates than those that hold for $\Phi_j$-waves. More precisely, we may replace Theorem \ref{thm-general wave table decomposition} with the following.

\begin{theorem}\label{thm-main bilinear estimate Up case}
Let  $\mb{C}_0>0$,  $\frac{1}{2}\les \frac{1}{q} \les 1$, $\frac{1}{2} \les \frac{1}{r} < \frac{n}{n+1}$, and  $\frac{1}{n+1}<\frac{1}{b}\les \frac{1}{a} \les \frac{1}{2}$. Assume that we have $\mb{d}_0>0$, open sets $\Lambda_j \subset \RR^n$, and phases $\Phi_j$ satisfying $(\mb{A1})$, $(\mb{A2})$, and the normalisation $(\ref{eqn-assum on phase after rescaling})$. Let $Q_R$ be a cube of diameter $R \gg ( \mb{d}_0^2 \mc{H}_2)^{-1}$. Then for any $0<\epsilon \ll 1$, any atomic $\Phi_1$-wave $u = \sum_{I \in \mc{I}} \ind_I u_I$, and any atomic  $\Phi_2$-wave $v = \sum_{J \in \mc{J}} \ind_J v_J$ with
		$$\supp \widehat{u} + \frac{\mb{d}_0}{2} \subset \Lambda_1, \qquad \supp \widehat{v} + \frac{\mb{d}_0}{2} \subset \Lambda_2$$
there exist a cube $Q$ of diameter $2R$ such that we have a decomposition
	$$ u = \sum_{B \in \mc{Q}_{\frac{2R}{4^M}}(Q)} u^{(B)}, \qquad v = \sum_{B' \in \mc{Q}_{\frac{R}{2}}(Q)} v^{(B')}$$
where $M\in \NN$ with $4^{-M} < \mc{H}_2 \les 4^{ 1- M}$, and $u^{(B)}= \sum_{I \in \mc{I}} \ind_I u_I^{(B)}$ is an atomic $\Phi_1$-wave, $v^{(B')}=\sum_{J \in \mc{J}} \ind_J v_J^{(B')}$ is an atomic $\Phi_2$-wave,  with the support properties
	$$ \supp \widehat{u}^{(B)} \subset \supp \widehat{u} + 2\big(2 \mc{H}_2 R\big)^{-\frac{1}{2}}, \qquad \supp \widehat{v}^{(B')} \subset \supp \widehat{v} +  2 \big(2\mc{H}_2 R\big)^{-\frac{1}{2}}.$$
Moreover, for any $a_0, b_0 \g 2$ we have the energy bounds
	$$ \Big(\sum_{B \in \mc{Q}_{\frac{2R}{4^M}}(Q)} \|u^{(B)}\|_{\ell^{a_0} L^2_x}^{a_0} \Big)^\frac{1}{a_0} \les (1 +C\epsilon) \| u \|_{\ell^{a_0} L^2_x}$$
	 $$\Big(\sum_{B' \in \mc{Q}_{\frac{R}{2}}(Q)} \|v^{(B')}\|_{\ell^{b_0} L^2_x}^{b_0} \Big)^\frac{1}{b_0} \les (1 + C \epsilon) \| v \|_{\ell^{b_0} L^2_x}$$
and the bilinear estimate
	\begin{align*}  \| u v \|_{L^q_t L^r_x(Q_R)} &\les (1+C\epsilon)  \big\| \big[u^{(\cdot)}\big] \big[v^{(\cdot)}\big] \big\|_{L^q_t L^r_x(Q)} \\
    &+ C \epsilon^{-C} R^{\frac{1}{q} + \frac{n+1}{2r} - \frac{n+1}{2}} \mc{H}_2^{\frac{n+1}{2r}-\frac{n}{2} + (n+1)(\frac{1}{a} - \frac{1}{2})}\Big( ( \mc{H}_2 R)^{\frac{n+1}{2}} ( \mu + \mb{d}_0)^{n}\Big)^{(1-\frac{1}{r}-\frac{1}{b})_+} \|u \|_{\ell^a L^2_x} \| v \|_{\ell^b L^2_x}
    \end{align*}
where the constant $C$ depends only on $\mb{C}_0$, $q$, $r$, and $n$ and we let $ \mu = \min\{ \diam(\supp \widehat{u}), \diam(\supp \widehat{v})\} $.
\end{theorem}

Note that by taking $\frac{1}{a}=\frac{1}{2}$ and $\frac{1}{b}\g  1 - \frac{1}{r}$, we essentially recover Theorem \ref{thm-general wave table decomposition}. \\

We leave the proof of Theorem \ref{thm-main bilinear estimate Up case} till Section \ref{sec:proof of atomic wave tables}, and now turn to the proof of Theorem \ref{thm-main Up version II}. Fix constants $\mb{d}_0, \mb{C}_0 >0$, and open sets $\Lambda_1, \Lambda_2 \subset \RR^n$. Let $\Phi_j$ be phases satisfying $(\mb{A1})$ and $(\mb{A2})$. As previously, by exploiting dilation and translation invariance, we may assume $\mc{H}_2 \les \mc{H}_1$ and the normalisation conditions \eref{eqn-assum on phase after rescaling}. Take sets $\Lambda_j^* + \mb{d}_0 \subset \Lambda_j$ such that $\mb{d}[\Lambda_1^* + \mb{d}_0, \Lambda_2^*+\mb{d}_0] \les \frac{\mb{d}_0}{\mb{C}_0}$, and fix $R_0 \gg \frac{1}{\mb{d}_0^2 \mc{H}_2}$ with $  (R_0 \mc{H}_2)^{-\frac{1}{2}} \approx \mb{d}_0$. The constant $R_0$ will denote the smallest scale of cubes we consider, while the sets $\Lambda_1^*$ and $\Lambda_2^*$ will contain the support of $\widehat{u}$ and $\widehat{v}$ respectively.

\begin{definition}\label{defn of A*(R)}
For any $R \g R_0$ and $1\les a,b, q, r \les 2$, we define $A^*(R)$ to the best constant for which the inequality
	$$ \| uv \|_{L^q_t L^r_x(Q)} \les A^*(R) \|u \|_{\ell^a L^2_x} \| v \|_{\ell^b L^2_x}$$
holds for all cubes $Q \subset \RR^{1+n}$ of radius $R$, and all  atomic $\Phi_1$-waves $u$ and atomic $\Phi_2$-waves $v$ satisfying the support assumption
	$$ \supp \widehat{u}  \subset \Lambda_1^* + 4(\mc{H}_2 R)^{-\frac{1}{2}}, \qquad \supp \widehat{v}  \subset \Lambda_2^* + 4 (\mc{H}_2 R)^{-\frac{1}{2}}.$$
\end{definition}

As in the proof of Theorem \ref{thm-main small scale} in Section \ref{sec-ind on scales}, the key step in the proof of Theorem \ref{thm-main Up version II} is to prove the following bounds on $A^*(R)$.

\begin{proposition}\label{prop-induction bound Up version}
Let $\mb{C}_0>0$,  $\frac{1}{2}\les \frac{1}{q} \les 1$, $\frac{1}{2} \les \frac{1}{r} < \frac{n}{n+1}$, $\frac{1}{n+1}<\frac{1}{b}\les \frac{1}{a} \les \frac{1}{2}$ and $\frac{1}{a} + \frac{1}{b} \g \frac{1}{\min\{q, r\}}$. There exists a constant $C>0$, such that for any $ \mb{d}_0>0$, any open sets $\Lambda_j \subset \RR^n$, any phases $\Phi_j$ satisfying $(\mb{A1})$ and $(\mb{A2})$ with the normalisation (\ref{eqn-assum on phase after rescaling}), and any sets $\Lambda_j^* + \mb{d}_0 \subset \Lambda_j$ with $\mb{d}[\Lambda_1^* + \mb{d}_0, \Lambda_2^* + \mb{d}_0] \les \frac{ \mb{d}_0}{\mb{C}_0}$, we have for every $R \g R_0$ and $0<\epsilon \ll 1$
    \begin{equation}\label{eqn-prop induction bound Up version-induction bound}
        \begin{split}  A^*(2R) \les ( 1 +& C \epsilon) A^*(R) \\
        &+ C \epsilon^{-C} R^{\frac{1}{q} + \frac{n+1}{2r} - \frac{n+1}{2}} \mc{H}_2^{\frac{n+1}{2r}-\frac{n}{2} + (n+1)(\frac{1}{a} - \frac{1}{2})}\Big( ( \mc{H}_2 R)^{\frac{n+1}{2}} ( \mu + \mb{d}_0)^{n}\Big)^{(1-\frac{1}{r}-\frac{1}{b})_+}
        \end{split}
    \end{equation}
and
 \begin{equation}\label{eqn-prop induction bound Up version-initial bound}A^*(2R) \les C \mb{d}_0^{n+1 - \frac{n+1}{r} - \frac{2}{q}}\mc{H}_2^{\frac{1}{2} - \frac{1}{q} + (n+1)(\frac{1}{a} - \frac{1}{2})} \Big( \mb{d}_0^{-(n+1)} (\mu + \mb{d}_0)^n \Big)^{(1-\frac{1}{r}  - \frac{1}{b})_+} \Big( \frac{R}{R_0}\Big)^{\frac{1}{q}}\end{equation}
 where $\mu = \min\{ \diam(\Lambda_1^*), \diam(\Lambda_2^*)\}$.
\end{proposition}
\begin{proof}
Let $Q$ be a cube of diameter $2R$, and let $u = \sum_{I \in \mc{I}} \ind_I(t) u_I$ be a $U^a_{\Phi_1}$ atom, and $v = \sum_{J \in \mc{J}} \ind_J(t) v_J$ be a $U^b_{\Phi_2}$ atom satisfying the support conditions
	$$ \supp \widehat{u} \subset \Lambda_1^* + 4(\mc{H}_2 2R)^{-\frac{1}{2}}, \qquad \supp \widehat{v}  \subset \Lambda_2^* + 4(\mc{H}_2 2R)^{-\frac{1}{2}} .$$
An application of Theorem \ref{thm-main bilinear estimate Up case} gives a cube $Q'$ of diameter $4R$, and atomic waves $(u^{(B)})_{B \in \mc{Q}_{\frac{R}{4^{M-1}}}(Q)}$, $(v^{(B')})_{B'\in \mc{Q}_{R}(Q)}$ such that
	\begin{equation}\label{eqn-prop induction bound Up version-bounded by quilt}
        \begin{split} \| u v\|_{L^q_t L^r_x(Q)} \les (1 + &C\epsilon) \big\| [u^{(\cdot)}] [v^{(\cdot)}] \big\|_{L^q_t L^r_x(Q')} \\
        &+  C \epsilon^{-C} R^{\frac{1}{q} + \frac{n+1}{2r} - \frac{n+1}{2}} \mc{H}_2^{\frac{n+1}{2r}-\frac{n}{2} + (n+1)(\frac{1}{a} - \frac{1}{2})}\Big( ( \mc{H}_2 R)^{\frac{n+1}{2}} ( \mu + \mb{d}_0)^{n}\Big)^{(1-\frac{1}{r}-\frac{1}{b})_+}
        \end{split}
\end{equation}

and the support properties
	$$ \supp \widehat{u}^{(B)} \subset \Lambda_1^* + 4 (\mc{H}_2 R) ^{-\frac{1}{2}}, \qquad \qquad \supp \widehat{v}^{(B')} \subset \Lambda_2^* + 4(\mc{H}_2 R)^{-\frac{1}{2}}$$
where we used the support assumptions on $\widehat{u}$ and $\widehat{v}$.

To prove (\ref{eqn-prop induction bound Up version-induction bound}), we let $B' \in \mc{Q}_{R}(Q')$ and define the atomic $\Phi_1$-wave $U^{(B')}= \sum_{I\in \mc{I}} \ind_I(t) U^{(B')}_I$ with $U^{(B')}_I = ( u^{(B)}_I )_{B \in \mc{Q}_{\frac{R}{4^{M-1}}}(B')}$. Then for every $B'\in \mc{Q}_{R}(Q)$ we have an atomic $\Phi_1$-wave $U^{(B')}$ and an atomic $\Phi_2$-wave $v^{(B')}$ satisfying the correct support assumptions to apply the definition of $A^*(R)$. Thus
	\begin{align*}
	   \big\| [ u^{(\cdot)} ] [ v^{(\cdot)}] \big\|_{L^q_t L^r_x(Q')} &\les \bigg( \sum_{B' \in \mc{Q}_R(Q')} \| U^{(B')} v^{(B')} \|_{L^q_t L^r_x(B')}^{\min\{q, r\}} \bigg)^{\frac{1}{\min\{q, r\}}} \\
	   &\les A^*(R) \Big( \sum_{B' \in \mc{Q}_R(Q')} \|U^{(B')}\|_{\ell^a L^2_x}^{a} \Big)^\frac{1}{a}\Big( \sum_{B' \in \mc{Q}_R(Q')} \|v^{(B')}\|_{\ell^b L^2_x}^{b} \Big)^\frac{1}{b} \\
	   &\les (1+C \epsilon) A^*(R)
	\end{align*}
where the second line used the assumption $\frac{1}{a} + \frac{1}{b} \g  \frac{1}{\min\{q, r\}}$ and the last applied the energy inequalities in Theorem \ref{thm-main bilinear estimate Up case}. Therefore the induction bound \eref{eqn-prop induction bound Up version-induction bound} follows from an application of (\ref{eqn-prop induction bound Up version-bounded by quilt}).

We now turn to the proof of (\ref{eqn-prop induction bound Up version-initial bound}). We begin by observing that again using the bound (\ref{eqn-prop induction bound Up version-bounded by quilt}) with $\epsilon\approx 1$,  the definition of $R_0$, and using the nesting properties of the $\ell^p$ spaces, it is enough to prove that for every $\frac{1}{b} \les \frac{1}{a} \les \frac{1}{2}$ and $\frac{1}{b}\les 1-\frac{1}{r}$ we have the quilt bound
     \begin{equation}\label{eqn-prop induction bound Up version-initial quilt bound}\big\| [u^{(\cdot)}] [v^{(\cdot)}] \big\|_{L^2_t L^r_x(Q')} \lesa  \mb{d}_0^{\frac{n-1}{b}} ( R \mc{H}_2)^{\frac{1}{2} - \frac{1}{b}} ( \mu + \mb{d}_0)^{n(1-\frac{1}{r}-\frac{1}{b})} \mc{H}_2^{(n+1)(\frac{1}{a} - \frac{1}{2})} \| u \|_{\ell^a L^2} \| v \|_{\ell^b L^2}.\end{equation}
To this end, we first note that applying H\"{o}lder's inequality gives for any $a_0 \g 2$
		\begin{align} \big\| [u^{(\cdot)}] [v^{(\cdot)}] \big\|_{L^2_t L^1_x(Q')}
			&\les \big\| [u^{(\cdot)}] \big\|_{L^2_{t,x}(Q')} \big\| [v^{(\cdot)}] \big\|_{L^\infty_t L^2_x(Q')} \notag\\
			& \lesa (R \mc{H}_2)^\frac{1}{2} \mc{H}_2^{(n+1)(\frac{1}{a_0} - \frac{1}{2})} \Big(\sum_{B \in \mc{Q}_{\frac{2R}{4^M}}(Q')} \| u^{(B)} \|_{L^\infty_t L^2_x}^{a_0} \Big)^\frac{1}{a_0} \sup_{B' \in \mc{Q}_{\frac{R}{2}(Q')}} \| v^{(B')} \|_{L^\infty_t L^2_x}\notag\\
			&\lesa (R \mc{H}_2)^\frac{1}{2} \mc{H}_2^{(n+1)(\frac{1}{a_0} - \frac{1}{2})} \|u\|_{\ell^{a_0} L^2} \| v \|_{\ell^{\infty} L^2}.
\label{eqn-prop induction bound Up version-initial L2L1 bound}
		\end{align}
Similarly, we have the $L^2_{t,x}$ bound
    \begin{align}
      \big\| [u^{(\cdot)}] [v^{(\cdot)}] \big\|_{L^2_{t,x}(Q')} &\lesa \mc{H}_2^{(n+1)(\frac{1}{a_0} - \frac{1}{2})} (R \mc{H}_2)^\frac{1}{2} \sup_{B'} \bigg( \sum_{B} \| u^{(B)} v^{(B')}\|_{L^\infty_t L^2_x(B)}^{a_0}\bigg)^\frac{1}{a_0} \notag \\
      &\lesa \mc{H}_2^{(n+1)(\frac{1}{a_0} - \frac{1}{2})} ( R \mc{H}_2)^{\frac{1}{2}} (\mu + \mb{d}_0)^{ \frac{n}{2}} \bigg( \sum_B \| u^{(B)} \|_{L^\infty_t L^2_x}^{a_0}\bigg)^\frac{1}{a_0} \sup_{B'} \| v^{(B')} \|_{L^\infty_t L^2_x} \notag \\
      &\lesa \mc{H}_2^{(n+1)(\frac{1}{a_0} - \frac{1}{2})} ( R \mc{H}_2)^{\frac{1}{2}} ( \mu + \mb{d}_0)^{\frac{n}{2}} \| u \|_{\ell^{a_0} L^2} \| v \|_{\ell^\infty L^2}. \label{eqn-prop induction bound Up version-linear L2 bound}
    \end{align}
On the other hand, as in the proof of Proposition \ref{prop-initial induction bound}, an application of Theorem \ref{thm-classical bilinear L2 estimate} gives
\begin{align}
   \big\| [u^{(\cdot)}] [v^{(\cdot)}] \big\|_{L^2_{t,x}(Q')} &\les \bigg( \sum_{I \in \mc{I}} \sum_{J \in \mc{J}} \sum_{B\in \mc{Q}_{\frac{2R}{4^M}}(Q')} \sum_{B' \in \mc{Q}_{\frac{R}{2}}(Q')} \big\| u^{(B)}_I v^{(B')}_J \big\|_{L^2_{t,x}}^2\bigg)^\frac{1}{2}
   \lesa \mb{d}^{\frac{n-1}{2}}_0 \| u \|_{\ell^2 L^2} \| v \|_{\ell^2 L^2}. \label{eqn-prop induction bound Up version-bilinear L2}
\end{align}
Therefore the required bound \eref{eqn-prop induction bound Up version-initial quilt bound} follows by interpolating\footnote{This is simply an application of complex interpolation. In more detail, we just repeat the proof of the Riesz-Thorin interpolation theorem, thus given a functions $g \in L^2_t(\RR)$ and  $G \in L^\infty_t L^{r'}_x(\RR^{1+n})$, we define the function $\rho(z)$ for $z\in \CC$ as
    $$\rho(z) = \int_{X[Q]} \bigg( \sum_{I\in \mc{I}} \ind_I(t) \| u_I\|_{L^\infty_t L^2_x}^{a(\frac{z}{a_0} + \frac{1-z}{2}) -1} [u_{I,m}^{(\cdot)}]_{\frac{2R}{4^m}} \bigg) \bigg( \sum_{J \in \mc{J}}  \ind_J(t) \| v_J \|_{L^\infty_t L^2_x}^{b(\frac{z}{b_0} + \frac{1-z}{2}) -1} |v_J|\bigg) g(t) |G|^{r'( \frac{1}{2}+\frac{z}{2}) -1}   G \,dt\, dx $$
where $\frac{1}{a_0} = \frac{ \frac{1}{r} + \frac{1}{a}-1}{\frac{2}{r} -1}$ and $\frac{1}{b_0} = \frac{ \frac{1}{r} + \frac{1}{b}-1}{\frac{2}{r} -1}$. It is easy to check that $\rho$ is complex analytic when $0\les \Re(z)\les 1$, is at most of exponential growth, and $\rho(0)$ can be bounded by the $L^2_{t,x}$ estimate, $\rho(1)$ by the $L^2_t L^1_x$ estimate. Hence interpolated bound follows from the Hadamard Three-Lines Theorem or Lindel\"{o}f's Theorem. } between \eref{eqn-prop induction bound Up version-initial L2L1 bound}, \eref{eqn-prop induction bound Up version-linear L2 bound}, and the bilinear estimate \eref{eqn-prop induction bound Up version-bilinear L2}.
\end{proof}

Finally we come to the proof of Theorem \ref{thm-main Up version II}.

\begin{proof}[Proof of Theorem \ref{thm-main Up version II}]
After rescaling, we may assume the normalisation \eref{eqn-assum on phase after rescaling}.  The atomic definition of the $U^a_\Phi$ spaces, together with the definition of $A^*(R)$, implies that it is enough to show that for every $m \in \NN$ we have
        \begin{equation}\label{eqn-thm main Up ver-main ineq} A^*(2^m R_0) \lesa \mb{d}_0^{n+1 - \frac{n+1}{r} - \frac{2}{q}} \mc{H}_2^{\frac{1}{2} - \frac{1}{q} + (n+1)(\frac{1}{a} - \frac{1}{2})} \Big( \mb{d}_0^{-(n+1)} \mu^n \Big)^{(1-\frac{1}{r}-\frac{1}{b})_+},\end{equation}
with the implied constant independent of $m$. We now essentially repeat the proof of Theorem \ref{thm-main small scale}, but use Proposition \ref{prop-induction bound Up version} in place of Proposition \ref{prop-induction step} and Proposition \ref{prop-initial induction bound}. Let $\delta = \frac{n+1}{2} - \frac{1}{q} - \frac{n+1}{2r} -  \frac{n+1}{2}(1-\frac{1}{r} - \frac{1}{b})_+$, note that by assumption, we have $\delta>0$. An application of Proposition \ref{prop-induction bound Up version} with $R=2^m R_0$ and $\epsilon = \frac{1}{C} 2^{-\frac{m}{2C}\delta}$ implies that
    $$ A^*(R)( 2^m R_0) \les (1 + 2^{ -\frac{m}{2C}\delta}) A^*(2^{m-1} R_0) + C^{C+1} 2^{ -\frac{m}{2}\delta}\mb{d}_0^{n+1 - \frac{n+1}{r} - \frac{2}{q}} \mc{H}_2^{\frac{1}{2} - \frac{1}{q} + (n+1)(\frac{1}{a} - \frac{1}{2})} \Big( \mb{d}_0^{-(n+1)} \mu^n \Big)^{(1-\frac{1}{r}-\frac{1}{b})_+}.$$
In particular, since $\delta>0$, after $m$ applications of the bound \eref{eqn-prop induction bound Up version-induction bound} we have
	$$ A^*(2^m R_0) \lesa A^*(R_0) +\mb{d}_0^{n+1 - \frac{n+1}{r} - \frac{2}{q}} \mc{H}_2^{\frac{1}{2} - \frac{1}{q} + (n+1)(\frac{1}{a} - \frac{1}{2})} \Big( \mb{d}_0^{-(n+1)} \mu^n \Big)^{(1-\frac{1}{r}-\frac{1}{b})_+}$$
where the implied constant depends only on the exponents $q, r, a, b,$ the dimension $n$, and $\mb{C}_0$. The required inequality \eref{eqn-thm main Up ver-main ineq} now follows from another application of Proposition \ref{prop-induction bound Up version}.
\end{proof}

\section{Proof of Theorem \ref{thm-main bilinear estimate Up case}}\label{sec:proof of atomic wave tables}

Here we give the proof of the atomic wave table decomposition in Theorem \ref{thm-main bilinear estimate Up case}. The argument follows  that used to prove Theorem \ref{thm-general wave table decomposition} but, as in the proof of Proposition \ref{prop-induction bound Up version}, we interpolate with an improved $L^2_t L^1_x$ estimate, and an additional $L^2_{t,x}$ estimate to gain the $\ell^a$ and $\ell^b$ sums over intervals. Since $\ell^{b_1} \subset \ell^{b_2}$ for $b_1 \les b_2$, it is enough to consider the case $\frac{1}{n+1}< \frac{1}{b} \les 1- \frac{1}{r}$. Fix a cube $Q_R$ of diameter $R \gg (\mb{d}_0^2 \mc{H}_2)^{-1}$ and let $u=\sum_{I\in \mc{I}} \ind_I(t) u_I$ be an atomic $\Phi_1$-wave, and $v=\sum_{J \in \mc{J}} \ind_J(t) v_J$ be an atomic $\Phi_2$-wave with the support conditions
	 $$\supp \widehat{u} + \frac{\mb{d}_0}{2} \subset \Lambda_1, \qquad \supp \widehat{v} + \frac{\mb{d}_0}{2} \subset \Lambda_2. $$
 An application of Lemma \ref{lem - cube averaging} implies that there exists a cube $Q$ of radius $2R$ such that
        $$ \| uv \|_{L^q_t L^r_x(Q_R)} \les ( 1  + C \epsilon) \| uv \|_{L^q_t L^r_x(X[Q])}$$
where as in Theorem \ref{thm-general wave table decomposition} we take
	$$ X[Q] = \bigcap_{m=1,\dots, M} I^{\epsilon_m, 4^{-m}2R}(Q), \qquad \epsilon_m = 4^{\delta(m-M)} \epsilon$$
and $\delta>0$ is some fixed constant to be chosen later (which will depend only on the dimension $n$, and the constant $C_n$ appearing in Theorem \ref{thm-prop of wave tables}) and we take $M\in \NN$ such that $4^{-(M+1)} < \mc{H}_2 \les 4^{-M}$.

We start by decomposing the components $u_I$ of $u$. Let $V=(v_J)_{J\in \mc{J}}$, thus $V$ is a $\Phi_2$-wave such that $|v|\les |V|$. Given $B_1 \in \mc{Q}_{\frac{R}{2}}(Q)$ we let
	$$ u^{(B_1)}_{I, 1} = \mc{W}^{(B_1)}_{1, \epsilon_1}(u_I; V, Q)$$
(where $\mc{W}$ is as in Definition \ref{def:wave tables}) and assuming we have constructed $u_{I, m}^{(B_{m})}$, $B_{m} \in \mc{Q}_{\frac{2R}{4^{m}}}(Q)$, we define for $B_{m+1} \in \mc{Q}_{\frac{2R}{4^{m+1}}}(B_m)$
		$$ u_{I, m+1}^{(B_{m+1})} = \mc{W}^{(B_{m+1})}_{1, \epsilon_{m+1}}\big( u^{(B_{m})}_{I, m}; V, B_{m}\big).$$
To extend this to the atomic waves, we simply take $u_m^{(B_m)} = \sum_{I \in \mc{I}} \ind_I(t) u_{I, m}^{(B_m)}$.
Finally, for $B \in \mc{Q}_{\frac{2R}{4^M}}(Q)$, we let $ u^{(B)}=u_M^{(B)}$. Clearly $u^{(B)}$ is again an atomic $\Phi_1$-wave, and as in the proof Theorem \ref{thm-general wave table decomposition}, from an application of Theorem \ref{thm-prop of wave tables}, the pointwise decomposition and support properties of $u^{(B)}$ follow immediately. On the other hand, the energy inequality follows by exchanging the order of summation, using the fact that $a\g 2$, and applying the energy estimate in Theorem \ref{thm-prop of wave tables}
	\begin{align*}
  \bigg( \sum_{B \in \mc{Q}_{\frac{2R}{4^M}}(Q)} \| u^{(B)} \|_{\ell^a L^2_x}^a \bigg)^\frac{1}{a}
            &\les \bigg( \sum_{I \in \mc{I}} \bigg( \sum_{B \in \mc{Q}_{\frac{2R}{4^M}}(Q)} \big\| u^{(B)}_{I} \big\|_{L^\infty_t L^2_x}^2 \bigg)^\frac{a}{2}\bigg)^{\frac{1}{a}} \\
            &\les (1 + C \epsilon) \| u \|_{\ell^a L^2_x}.
\end{align*}
The next step is to decompose $v = \sum_{J \in \mc{J}} \ind_J(t) v_J$. Let $U = ( u_I^{(B)})_{I \in \mc{I}, B \in \mc{Q}_{\frac{2R}{4^M}}}$. Then $U$ is a $\Phi_1$-wave with the pointwise bound $[u^{(\cdot)}] \les |U|$ and the energy bound $\| U \|_{L^\infty_t L^2_x} \lesa \| u \|_{\ell^2 L^2_x}$.
We now decompose each $v_J$ relative to $U$ and the cube $Q$, in other words for every $B' \in \mc{Q}_{\frac{2R}{4}}(Q) $ we take
		$$ v^{(B')}_J = \mc{W}^{(B')}_{2, \epsilon}\big( v_J; U, Q\big)$$
and finally define $v^{(B')} = \sum_{J \in \mc{J}} \ind_J(t) v_J^{(B')}$. An application of Theorem \ref{thm-prop of wave tables} implies that $v^{(B')}$ is an atomic $\Phi_2$-wave and that $v^{(B')}$ satisfies the correct Fourier support conditions. Furthermore, by a similar argument to the $u^{(B)}$ case,  the required energy inequality also holds.

We now turn to the proof of the bilinear estimate. After observing that
        $$ \| uv \|_{L^q_t L^r_x(X[Q])} \les \big\| [u^{(\cdot)}] [v^{(\cdot)}] \big\|_{L^q_t L^r_x(X[Q])} + \big\| \big( |u| - [u^{(\cdot)}]\big) v \big\|_{L^q_t L^r_x(X[Q])} + \big\| [u^{(\cdot)}] \big( |v| - [v^{(\cdot)}]\big) \big\|_{L^q_t L^r_x(X[Q])}, $$
an application of Holder's inequality implies that it is enough to show that for $\frac{1}{n+1} < \frac{1}{b} \les 1-\frac{1}{r}$, and $ \frac{1}{b} \les \frac{1}{a} \les \frac{1}{2}$  we have
    \begin{equation}\label{eqn-thm main bilinear Up est-temp bilinear est}
    		\begin{split}
        \big\| \big( |u| - [u^{(\cdot)}]\big) v \big\|_{L^2_t L^r_x(X[Q])} +& \big\| [u^{(\cdot)}] \big( |v| - [v^{(\cdot)}]\big) \big\|_{L^2_t L^r_x(X[Q])}\\
                &\lesa \epsilon^{-C}\big( \mc{H}_2 R \big)^{\frac{1}{2} - \frac{n+1}{2b}} (\mu+\mb{d}_0)^{n(1-\frac{1}{r} - \frac{1}{b})}\mc{H}_2^{(n+1)(\frac{1}{a} - \frac{1}{2})} \| u \|_{\ell^a L^2} \| v \|_{\ell^b L^2}
         \end{split}
    \end{equation}
with $\mu = \min\{ \diam(\supp \widehat{u}), \diam(\supp \widehat{v})\}$. We start by estimating the first term. The point is to interpolate between a ``bilinear'' $L^2_{t,x}$ estimate which decays in $R$, and  ``linear'' $L^2_t L^1_x$ and $L^2_{t,x}$ estimates which can lose powers of $R$. The bilinear $L^2_{t,x}$ bound follows from an application of Theorem \ref{thm-prop of wave tables}
	\begin{align}
		\Big\| \Big( \big[ u^{(\cdot)}_{m-1}\big]-& \big[ u^{(\cdot)}_{m} \big] \Big) v \Big\|_{L^2_{t,x}(X[Q])}^2   \notag \\
		&\les \sum_{I \in \mc{I}} \sum_{B_{m-1} \in \mc{Q}_{\frac{2R}{4^{m-1}}}(Q)} \Big\| \Big( |u^{(B_{m-1})}_{I, m-1}| - \big[ \mc{W}^{(\cdot)}_{1, \epsilon_m}(u^{(B_{m-1})}_{I, m-1}; V, Q)\big] \Big) V \Big\|_{L^2_{t,x}(I^{\epsilon_m, \frac{2R}{4^{m}}}(B_{m-1}))}^2\notag \\
		&\lesa \epsilon_m^{-2C_n} \Big( \frac{4^{m-1}}{2R}\Big)^{\frac{n-1}{2}}  \sum_{I \in \mc{I}}\sum_{B_{m-1} \in \mc{Q}_{\frac{2R}{4^{m-1}}}(Q)} \| u^{(B_{m-1})}_{I, m-1} \|_{L^\infty_t L^2_x}^2 \| V \|_{L^\infty_t L^2_x}^2  \notag \\
		&\lesa \epsilon^{-2C_n}  4^{-2(M-m)(\frac{n-1}{4} - \delta C_n)} \big( \mc{H}_2 R \big)^{-\frac{n-1}{2}} \| u \|_{\ell^2 L^2}^2 \| v \|_{\ell^2 L^2}^2 \label{eqn-thm main bilinear Up-bilinear L2 estimate}
	\end{align}
where we used the definition of $\epsilon_m$ and $M$. To obtain the $L^2_t L^1_x$ bound, we start by noting that for any $1\les m \les M$, and $a_0\g2$ we have
    \begin{align}
    \bigg( \sum_{B \in \mc{Q}_{\frac{2R}{4^m}}(Q) } \big\| u^{(B)}_{m} \big\|_{L^\infty_t L^2_x}^2 \bigg)^\frac{1}{2}
                &\lesa   \bigg( \sum_{B \in \mc{Q}_{\frac{2R}{4^m}}(Q) } \big\| u^{(B)}_{m} \big\|_{L^\infty_t L^2_x}^2 \bigg)^\frac{1}{2} \notag \\
                &\lesa   4^{m (\frac{n+1}{2} - \frac{n+1}{a_0})} \bigg( \sum_{B \in \mc{Q}_{\frac{2R}{4^m}}(Q) } \big\| u^{(B)}_{m} \big\|_{\ell^{a_0} L^2_x}^{a_0} \bigg)^\frac{1}{a_0}  \notag \\
                &\lesa \mc{H}_2^{ (n+1)(\frac{1}{a_0}-\frac{1}{2})}  \| u \|_{\ell^{a_0} L^2_x} \label{eqn-thm main bilinear U p-L2 linear bound for u}
    \end{align}
where we applied the energy inequality for $u_m^{(B)}$. Therefore, an application of Holder's inequality gives for any $a_0 \g 2$
	\begin{align}
	\Big\| \Big( \big[ u^{(\cdot)}_{m-1}\big] -\big[ u^{(\cdot)}_{m} \big]\Big) v \Big\|_{L^2_t L^1_x(X[Q])}
	&\les \Big\| \Big( \big[ u^{(\cdot)}_{m-1}\big]- \big[ u^{(\cdot)}_{m} \big] \Big)\Big\|_{L^2_{t,x}(Q)}  \| v\|_{L^\infty_t L^2_x(Q)}\notag \\
	&\lesa \Big( \frac{R}{4^m}\Big)^\frac{1}{2} \Bigg(\sum_{B_{m-1}} \| u^{(B_{m-1})}_{m-1}\|_{L^\infty_t L^2_x(B_{m-1})}^2 + \sum_{B_m} \| u^{(B_m)}_m \|_{L^\infty_tL^2_x(B_m)}^2 \Bigg)^\frac{1}{2}  \| v \|_{\ell^{\infty} L^2} \notag \\
	&\lesa  \big( \mc{H}_2 R\big)^{\frac{1}{2}} \mc{H}_2^{ (n+1)(\frac{1}{a_0}-\frac{1}{2})}  4^{(M-m) \frac{1}{2}} \| u \|_{\ell^{a_0} L^2} \|v \|_{\ell^{\infty} L^2}. \label{eqn-thm main bilinear Up-bilinear L2L1}
	\end{align}
On the other, for the linear $L^2_{t,x}$ bound, we can apply a similar argument to deduce that
	\begin{align}
		\Big\| \Big( \big[ u^{(\cdot)}_{m-1}\big] -\big[ u^{(\cdot)}_{m} \big]\Big) v \Big\|_{L^2_{t,x}(X[Q])}
			&\lesa \Big( \frac{R}{4^m}\Big)^\frac{1}{2}  \Bigg(\sum_{B_{m-1}} \| u^{(B_{m-1})}_{m-1} v\|_{L^\infty_t L^2_x (B_{m-1})}^2 + \sum_{B_m} \| u^{(B_m)}_m v \|_{L^\infty_t L^2_x(B_m)}^2 \Bigg)^\frac{1}{2}   \notag \\
            &\lesa (\mu + \mb{d}_0)^{\frac{n}{2}} \Bigg(\sum_{B_{m-1}} \| u^{(B_{m-1})}_{m-1}\|_{L^\infty_t L^2_x}^2 + \sum_{B_m} \| u^{(B_m)}_m \|_{L^\infty_t L^2_x}^2 \Bigg)^\frac{1}{2} \| v\|_{L^\infty_t L^2_x}   \notag \\
            &\lesa (\mu + \mb{d}_0)^{\frac{n}{2}} \big( \mc{H}_2 R\big)^{\frac{1}{2}} \mc{H}_2^{ (n+1)(\frac{1}{a_0}-\frac{1}{2})}  4^{(M-m) \frac{1}{2}} \| u \|_{\ell^{a_0} L^2} \|v \|_{\ell^{\infty} L^2} \label{eqn-thm main bilinear Up-linear L2}
    \end{align}
where we used the fact that at least one of $u^{(B)}$ of $v^{(B')}$ has Fourier support contained in a set of diameter $\mb{d}_0+\mu$. Interpolating between (\ref{eqn-thm main bilinear Up-bilinear L2 estimate}), (\ref{eqn-thm main bilinear Up-bilinear L2L1}), and (\ref{eqn-thm main bilinear Up-linear L2}) then gives for any $\frac{1}{2}\les \frac{1}{r} \les 1$, $\frac{1}{b} \les 1-\frac{1}{r}$, and  $\frac{1}{b}\les \frac{1}{a} \les \frac{1}{2}$,
	\begin{align*} \Big\| \Big( \big[ u^{(\cdot)}_{m-1}\big]_{\frac{2R}{4^{m-1}}} &- \big[ u^{(\cdot)}_{m} \big]_{\frac{2R}{4^m}} \Big) v \Big\|_{L^2_t L^r_x (X[Q])}\\
			&\lesa \epsilon^{-C}  4^{-(M-m)\delta^*} \big( \mc{H}_2 R \big)^{\frac{1}{2} - \frac{n+1}{2b}} (\mu + \mb{d}_0)^{n(1-\frac{1}{r} - \frac{1}{b})}\mc{H}_2^{(n+1)(\frac{1}{a} - \frac{1}{2})} \| u \|_{\ell^a L^2_x} \| v \|_{\ell^b L^2_x}
	\end{align*}
where $\delta^*  = \frac{n+1}{2b} - \frac{1}{2} - 2 \delta C_n \frac{1}{b}$. Consequently, provided that $\frac{1}{n+1}<\frac{1}{b} \les 1 - \frac{1}{r}$, and we choose $\delta$ sufficiently small depending only on $C_n$, $b$, and $n$, we have $\delta^*>0$. Thus by telescoping the sum over $m$ and letting $ u_0^{(Q)} = u$, we deduce that
	\begin{align*} \big\| \big( |u| - [u^{(\cdot)}]\big) v \big\|_{L^2_t L^r_x(X[Q])} &\les \sum_{m=1}^M \big\| \big( [ u^{(\cdot)}_{m-1}] - [ u^{(\cdot)}_{m} ] \big) v \big\|_{L^2_t L^r_x(X[Q])} \\
	&\lesa \epsilon^{-C}\big( \mc{H}_2 R \big)^{\frac{1}{2} - \frac{n+1}{2b}} (\mu + \mb{d}_0)^{n(1-\frac{1}{r} - \frac{1}{b})}\mc{H}_2^{(n+1)(\frac{1}{a} - \frac{1}{2})} \| u \|_{\ell^a L^2_x} \| v \|_{\ell^b L^2_x}.
	\end{align*}
It only remains to estimate the second term on the left hand side of \eref{eqn-thm main bilinear Up est-temp bilinear est}. To this end, applying the definition of $v^{(B')}$ together with Theorem \ref{thm-prop of wave tables}, we have
		\begin{align*} \big\| [u^{(\cdot)}] \big( |v| - [v^{(\cdot)}] \big)\big\|_{L^2_{t,x} (X[Q])}^2 &\les\sum_{J \in \mc{J}} \big\| U \big( |v_J| - [v^{(\cdot)}_J] \big)\big\|_{L^2_{t,x} (I^{\epsilon_1, \frac{R}{2}}(Q))}^2  \lesa  \epsilon^{-2C_n}  \big( \mc{H}_2 R \big)^{-\frac{n-1}{2}} \| u \|_{\ell^2 L^2}^2 \| v \|_{\ell^2 L^2}^2.
        \end{align*}
On the other hand an application of Holder's inequality together with the energy estimates and (\ref{eqn-thm main bilinear U p-L2 linear bound for u}) gives for any $a_0 \g 2$ the bounds
	 \begin{align}  \big\| [u^{(\cdot)}] \big( |v| - [v^{(\cdot)}] \big)\big\|_{L^2_t L^1_x (X[Q])} &\les \bigg( \sum_{B \in \mc{Q}_{\frac{2R}{4^M}}} \big\| u^{(B)} \|_{L^2_{t,x}(B)}^2 \bigg)^\frac{1}{2} \sup_{J \in \mc{J}} \Big( \|v_J \|_{L^\infty_t L^2_x} + \| [v^{(\cdot)}_J] \|_{L^\infty_t L^2_x}\Big) \notag \\
	 &\lesa  \big( \mc{H}_2 R \big)^\frac{1}{2} \mc{H}_2^{ (n+1)(\frac{1}{a_0} - \frac{1}{2})}  \| u \|_{\ell^{a_0} L^2} \| v \|_{\ell^{\infty} L^2} \label{eqn-thm main bilinear Up est-v bilinear L2L1 bound}
	 \end{align}
and
	\begin{align}  \big\| [u^{(\cdot)}] \big( |v| - [v^{(\cdot)}] \big)\big\|_{L^2_{t,x} (X[Q])}
	&\lesa (\mu + \mb{d}_0)^{\frac{n}{2}} \bigg( \sum_{B \in \mc{Q}_{\frac{2R}{4^M}}} \big\| u^{(B)} \|_{L^2_{t,x}(B)}^2 \bigg)^\frac{1}{2} \sup_{J \in \mc{J}} \Big( \|v_J \|_{L^\infty_t L^2_x} + \| [v^{(\cdot)}_J] \|_{L^\infty_t L^2_x}\Big) \notag \\
	 &\lesa (\mu +  \mb{d}_0)^{\frac{n}{2}} \big( \mc{H}_2 R \big)^\frac{1}{2} \mc{H}_2^{ (n+1)(\frac{1}{a_0} - \frac{1}{2})}  \| u \|_{\ell^{a_0} L^2} \| v \|_{\ell^{\infty} L^2}.
	 \label{eqn-thm main bilinear Up est-v linear L2 bound}
	 \end{align}
Therefore, by interpolation, we deduce that for $\frac{1}{2}\les \frac{1}{r} \les 1$, $\frac{1}{b} \les 1-\frac{1}{r}$, and  $\frac{1}{b}\les \frac{1}{a} \les \frac{1}{2}$, we have
	$$ \big\| [u^{(\cdot)}] \big( |v| - [v^{(\cdot)}] \big)\big\|_{L^2_t L^r_x(X[Q])} \lesa \epsilon^{-C}\big( \mc{H}_2 R \big)^{\frac{1}{2} - \frac{n+1}{2b}} (\mu + \mb{d}_0)^{n(1-\frac{1}{r} - \frac{1}{b})}\mc{H}_2^{(n+1)(\frac{1}{a} - \frac{1}{2})}   \| u \|_{\ell^a L^2} \| v \|_{\ell^b L^2}$$
and consequently (\ref{eqn-thm main bilinear Up est-temp bilinear est}) follows.

\bibliographystyle{amsplain}

\providecommand{\bysame}{\leavevmode\hbox to3em{\hrulefill}\thinspace}
\providecommand{\MR}{\relax\ifhmode\unskip\space\fi MR }
\providecommand{\MRhref}[2]{%
  \href{http://www.ams.org/mathscinet-getitem?mr=#1}{#2}
}
\providecommand{\href}[2]{#2}

\end{document}